\documentclass[a4paper,10pt,reqno]{article}
\usepackage[utf8x]{inputenc}
\usepackage[leqno]{amsmath}
\usepackage{amsfonts}
\usepackage{amssymb}
\usepackage{amsthm}
\usepackage[T1]{fontenc}
\usepackage{graphicx,caption}

\newtheorem{thm}{Theorem}
\newtheorem{cor}{Corollary}
\newtheorem{lem}{Lemma}
\newtheorem{prop}{Proposition}

\theoremstyle{definition}

\theoremstyle{remark}
\newtheorem*{rmk}{Remark}



\newcommand{\R}{\mathbb R}

\newcommand{\N}{\mathbb N}
\newcommand{\C}{\mathbb C}

\newcommand{\T}{\mathbb T}

\newcommand{\ka}[1]{\mathcal{C}^{#1}}
\newcommand{\Lip}[1]{\operatorname{Lip}{#1}}

\newcommand{\Div}{\operatorname{div}}
\newcommand{\spt}{\operatorname{supt}}
\newcommand{\rot}{\operatorname{curl}}

\newcommand{\dmt}{\operatorname{diam}}
\title{On the analyticity of the trajectories of the particles in the patch problem for $2$D Euler and aggregation equations}
\author{J. M. Burgués \and J. Mateu}
\begin{document} 	

\maketitle

%\let\thefootnote\relax
%\footnotetext{Both authors are  partially supported by grants 2017-SGR-0395 (Generalitat de Catalunya) and MDM-2014-044 (MICINN, Spain). The second named author is also partially supported by MTM-2016-75390 (MINECO, Spain).}

\begin{abstract}We give a proof of the analyticity in time for the particle trajectories associated with the solutions of some transport equations when the initial datum is a patch. These results are obtained from a precise study of the Beurling transform, which provides estimates for the solutions of some equations satisfied by the lagrangian flow. 

\end{abstract}

% !TeX encoding = UTF-8
% !TeX spellcheck = ca_ES
% !TeX root = Analit4.tex 

\section{Introduction} 
Let us consider some classical transport equations in the plane. One of these is the $2$D Euler equation, in vorticity form, for an incompressible inviscid fluid:  
\begin{equation}\label{eq(E)}
\begin{cases}
\partial_t\omega+v\cdot\nabla\omega=0,\\
\omega(0)=\omega_0,
\end{cases}
\tag{$E$}
\end{equation}
where $\omega$ is a scalar function representing the vorticity, and $v$ is the velocity of the fluid. In this case the vorticity and the velocity of the fluid are related by the Biot--Savart law 
$$
v(t)=\omega(t)*\biggl(\frac{x^\perp}{|x|^2}\biggr),
$$ 
where $x\in\R^2$.

Another equation, related with biological systems, is the aggregation equation 
\begin{equation}\label{eqtildeA}
\begin{cases}
\partial_t\rho+v\cdot\nabla\rho=0,\\
\rho(0)=\rho_0,	
\end{cases}
\tag{$\tilde A$}
\end{equation}
$\rho$ representing the density of mass of an irrotational inviscid and compressible fluid, where 
$$
v(t)=\rho(t)*\biggl(\frac{x}{|x|^2}\biggr).
$$

Equation~\eqref{eqtildeA} is closely related to the continuity equation
\begin{equation}\label{eq(A)}
\begin{cases}
\partial_t\rho+\Div(\rho\, v)=0,\\
\rho(0)=\rho_0.
\end{cases}
\tag{$A$}
\end{equation}

See \cite{BeLaLe} for the relationship between equations~\eqref{eq(A)} and~\eqref{eqtildeA} and problems in Biology.

%\medskip

A classical theorem due to Yudovich asserts that for any initial datum $\omega_0$ in $(L^\infty\cap L^p)(\R^2)$, where $1<p<\infty$, there is a unique weak solution for the equation~\eqref{eq(E)}. Using similar arguments an analogous result is obtained for the equation~\eqref{eq(A)}, and in both cases the solution depends continuously on the time variable and is as regular in the space variables as the datum (see\ 
\cite{Yud}, \cite{Che}, \cite{MaBe}, \cite{BeLaLe}).

The aim of this paper is to prove real analyticity in time of the particle trajectories of these fluids. These trajectories are globally described by a flow associated to the velocity field $v$. 
The regularity in time has largely already been studied  (see \cite{Che2}, \cite{Her}, \cite{Shn}, \cite{FrZh}, \cite{Sue}, \cite{Ser}, \cite{Ser1}, \cite{Ser2}, \cite{BaKe}). We will use complex notation, owing to the presence of some intrinsec objects of complex analysis and also because it simplifies arguments and notation. 

%\bigskip
%\bigskip

\subsection{Basic notation} 

In $\C$ we consider the standard coordinate $z=x+i\, y$. Then 
$$
\frac{\partial}{\partial z}=\frac{1}{2} 
\biggl(\frac{\partial}{\partial x}-i\, \frac{\partial}{\partial y}\biggr),
\quad 
\frac{\partial}{\partial\bar z}=\frac{1}{2} 
\biggl(\frac{\partial}{\partial x}+i\, \frac{\partial}{\partial y}\biggr)
$$ 
and we have, identifying $v=(v_1,\, v_2)$ with $v=v_1+i\, v_2$, that 
$$
\rot v=2\, \Im \biggl(\frac{\partial v}{\partial z}\biggr)=\frac{1}{i} 
\biggl(\frac{\partial\, v}{\partial z}-\frac{\partial\bar v}{\partial\bar z}\biggr)
$$
and 
$$
\Div v=2\, \Re \biggl(\frac{\partial v}{\partial z}\biggr)=\frac{\partial\, v}{\partial z}+\frac{\partial\bar v}{\partial\bar z}.
$$

%\medskip

We will denote by $m$ the Lebesgue measure in $\R^2$, associated to the standard volume form in $\C$, $\frac{1}{2i}\, d\bar z\wedge dz$.

%\medskip

The {\it conjugate Cauchy transform}, inverse of the operator $\frac{\partial}{\partial z}$ and will be denoted by 
$$
\bar C[\varphi](z)=\frac{1}{2\pi i}\int_\C\frac{\varphi(\zeta)}{\bar\zeta-\bar z}\, d\bar z\wedge dz=\frac{1}{\pi}\int_\C\frac{\varphi(\zeta)}{\bar\zeta-\bar z}\, dm, 
$$ 
defined for suitable functions $\varphi$. Then 
$$
\frac{\partial}{\partial z}\bar C[\varphi]=\varphi
$$ 
and the derivative 
\begin{equation*}
\begin{split}
\bar B[\varphi](z)&=\biggl(\frac{\partial}{\partial\bar z}\bar C[\varphi]\biggr)(z)\\
&=\text{p.\ v.\ } \frac{i}{2\pi}\int_\C\frac{\varphi(\zeta)}{(\bar\zeta-\bar z)^2}\, d\bar z\wedge dz=\text{p.\ v.\ } \frac{-1}{\pi}\int_\C\frac{\varphi(\zeta)}{(\bar\zeta-\bar z)^2}\, dm, 
\end{split}
\end{equation*}
is the {\it conjugate Beurling transform}.

%\bigskip
%\bigskip

Let $\phi$ be a test function supported in the unit ball, whose integral is equal to $1$ and such that $\phi(0)=1$, we consider, for a distribution $T$ the limit 
\begin{equation}\label{test}
\lim_{\epsilon\to0} \langle T,\, \phi_{x_0,\, \epsilon}\rangle,
\end{equation} 
where $\phi_{x_0,\, \epsilon}(x)=\frac{1}{\epsilon^2}\, \phi(\frac{x-x_o}{\epsilon})$. 

If this limit exists at some point $x_0$ and is independent of the choice of $\phi$ we call it {\it density function of $T$ at $x_0$} and denote it by $\Theta(T,x_0)$.

%\medskip

In fact, we have

\begin{lem} 
Let $\Omega\subset\C$ be a bounded domain such that $\partial\Omega\in\ka{1}$. If $T_{\chi_\Omega}$ is the distribution given by $\chi_\Omega$, we have 
$$
\Theta(T_{\chi_\Omega},z_0)=\begin{cases} 
1  & \text{if }  z_0\in\Omega, \\ 
\frac{1}{2} & \text{if }  z_0\in\partial\Omega, \\	
0 &\text{if }  z_0\in\bar\Omega^c. 
\end{cases}
$$ 
\end{lem}

\begin{proof} Let $\rho$ be a conveniently chosen defining function for $\Omega$ (see \cite{Bur}). Take $x_0\in\partial\Omega$ and consider 
$$
\kappa(x;x_0)=(\nabla\rho(x_0),x-x_0)
$$ 
and
$$
\omega(x;x_0)=\rho(x)-\kappa(x;x_0).
$$ 
Then, taking the function $\phi_{x_0,\, \epsilon}$ used in the definition of \eqref{test},  we have 
$$
\int_{B_\epsilon(z_0)\cap\Omega} 
\phi_{z_0,\epsilon}(z)\, dm(z)=\int_{B_\epsilon(z_0)\cap\{\rho<0\}} 
\phi_{z_0,\epsilon}(z)\, dm(z)=(*),
$$ 
and since 
$$
\{\rho<0\}=(\{\rho<0\}\cap\{\kappa<0\})\cup(\{\rho<0\}\cap\{\kappa\geq0\}),
$$ 
then
\begin{equation*}
\begin{split}
(*)&=\int_{B_\epsilon(z_0)\cap\{\kappa<0\}} 
\phi_{z_0,\epsilon}(z)\, dm(z)\\
&\quad+\{\int_{B_\epsilon(z_0)\cap\{\rho<0\}\cap\{\kappa\geq0\}} 
\phi_{z_0,\epsilon}(z)\, dm(z)\\
&\quad-\int_{B_\epsilon(z_0)\cap\{\rho\geq0\}\cap\{\kappa<0\}} 
\phi_{z_0,\epsilon}(z)\, dm(z)\}=(I)+(II).
\end{split}
\end{equation*}

%\bigskip

So, after a rotation 
\begin{equation*}
\begin{split}
(I)&=
\lim_{\epsilon\to0}\int_{B_\epsilon(z_0)\cap\{\kappa<0\}} 
\phi_{z_0,\epsilon}(z)\, dm(z)\\
&=\lim_{\epsilon\to0}\int_{B_\epsilon(z_0)\cap\R^2_+}\phi_{z_0,\epsilon}(z)\, dm(z)=\frac{1}{2}.
\end{split}
\end{equation*}

On the other hand 
\begin{equation*}
\begin{split}
|(II)|&=
\biggl|\int_{\!B_\epsilon(z_0)\cap\{\rho<0\}\cap\{\kappa\geq0\}} \!
\phi_{z_0,\epsilon}(z)\, dm(z)-\!\!\int_{\!B_\epsilon(z_0)\cap\{\rho\geq0\}\cap\{\kappa<0\}} 
\!\phi_{z_0,\epsilon}(z)\, dm(z)\biggr|\!\\*[5pt]
&\le\frac{1}{\epsilon^2} \{ m(B_\epsilon(z_0)\cap\{\rho<0\}\cap\{0\le\kappa\})+m(B_\epsilon(z_0)\cap\{\rho\geq0\}\cap\{\kappa<0\})\},
\end{split}
\end{equation*}
and  since the set $\{\rho<0\}\cap B_\epsilon(z_0)$ can be described as the graph of a $\ka{1}$ function over $\{\kappa=0\}$ vanishing at $z_0$ at the order $1$, then both 
$$
m(B_\epsilon(z_0)\cap\{\rho<0\}\cap\{0\le\kappa\})
$$ 
and 
$$
m(B_\epsilon(z_0)\cap\{\rho\geq0\}\cap\{\kappa<0\})
$$ 
are of size $\epsilon\, o(\epsilon)
$, so $(II)\to_{\epsilon\to0}0$, and the lemma is proved. 
\end{proof}

%\bigskip
%\bigskip
%\bigskip
%\bigskip

We will use the space $\ka{k,\gamma}(U)$, where $U$ is an open subset of the plane,  $k$ is a non-negative integer and $0<\gamma<1$. 
It is the space of functions with continuous derivatives up to the order $k$ such that each derivative of order $k$ extends to a $\gamma$-H\"older function in the closure of $U$.

%\medskip

In this paper we will mainly use the spaces $\ka{k,\gamma}(U)$ for $k=0,1$, equipped with the norms 
$$
\|f\|_\gamma=\|f\|_\infty+\sup_{z\neq w;\, z,\, w\in U}\frac{|f(z)-f(w)|}{|z-w|^\gamma}
$$ 
and
$$
\|f\|_{1,\, \gamma}=\|f\|_\infty+\|\nabla f\|_\gamma=\|f\|_\infty+\biggl\|\frac{\partial f}{\partial\bar z}\biggr\|_\gamma+\biggl\|\frac{\partial f}{\partial z}\biggr\|_\gamma.
$$

%\bigskip

We let $\ka{\omega}$ stand for a space of real analytic functions and~$\mathcal{D}(A)$ for the compactly supported $\ka{\infty}$ functions whose support is contained in the closed set $A\subset\R^n$. 

%\bigskip
%\bigskip

\subsection{Statement of results}

Let $\Omega\subset\C$ be a bounded domain such that $\partial\Omega\in\ka{1,\gamma}$ for $\gamma\in(0,1)$. 

Let $f\in\ka{\infty}(\C)$ and consider the distribution $T=f\, T_{}\chi_\Omega$. Then it is an immediate fact that 
$T$ has a density function 
\begin{equation} \label{dns1} 
\Theta(T,z)=\begin{cases} 
f(z) &\text{if } z\in\Omega,\\
\frac{f(z)}{2} &\text{if } z\in\partial\Omega,\\
0  &\text{if } z\notin\bar\Omega.
\end{cases}
\end{equation}

We will consider 
$$
\varpi_0(z)=\Theta(T,z)
$$ 
as the initial datum of a Cauchy problem ({\it i.\ e.\ vorticity for the Euler equation and density for aggregation and transport equations})
and then we define 
$$
v_0=\bar C[\varpi_0]
$$ 
as the {\it initial velocity}. We have (cf.\ \cite[section~4.3.2]{AsIwMa}) that %\newline 
$v_0\in\Lip{(1,\, \C)}\cap\ka{\infty}(\C\setminus\partial\Omega).$

%\bigskip

We will establish the analyticity in time for the flow associated to the problems~\eqref{eq(E)}, \eqref{eq(A)} and 
\eqref{eqtildeA} in the case in which the initial datum of the corresponding {\it patch},~$\rho_0$ or~$\omega_0$, is a real multiple of the characteristic function of $\Omega$. In the current literature, this is the {\it patch problem}. 

We will prove (in Lemma \ref{jak}) that the trajectories of the particles associated to these problems satisfy the equation (\ref{Meq}) below. This equation relates the flow with the vorticity or the density for a fluid satisfying a transport or continuity equation.

Equation (\ref{Meq}) is used in an a priori treatment of the flow when the datum (the term in the right hand side of the equation) is analytic in time.

An inductive procedure provides a solution of (\ref{Meq}) that is analytic too and whose velocity field satisfies the equation mentioned above. The uniqueness of solution of the equations \eqref{eq(E)} and \eqref{eq(A)} concludes the argument.

This inductive procedure also prove the existence of an analytic solution for the problem (\ref{Meq}). In fact, we have 

\begin{thm}\label{Mth1}
Let $\varpi_0$ as in \eqref{dns1}, $T_0\in (0,\, \infty)$ and  
 $a\colon\C\times(-T_0,\, T_0)\rightarrow\C$ a function in $\ka{\infty}((\C\setminus \partial\Omega)\times(-T_0,\, T_0))$ and such that 
$$
t\rightarrow a(z,t)
$$ 
is analytic in $(-T_0,\, T_0)$ for every $z\in\C$. 
	
Then the problem  
\begin{equation}\label{Meq}  
\begin{cases}
(\frac{\partial^2\psi}{\partial t\partial z}
\overline{\frac{\partial\psi}{\partial z}}-\frac{\partial^2\psi}{\partial t\partial\bar z}
\overline{\frac{\partial\psi}{\partial\bar z}})(z,t)=\varpi_0(z)\, (1+a(z,\, t)\, t),\\*[5pt]
\psi(z,0)=z	,
\end{cases}
\end{equation} 
has a solution $\psi$ in $\C\times(-T_0,\, T_0)$, such that $t\rightarrow\psi(z,\, t)$ is in $\ka{\omega}((-T_0,\, T_0))$ for every~$z\in\C$.
\end{thm}

%\bigskip

\begin{rmk} It is worth to observe that for a general $a$ (as in the statement) the uniqueness of the solution of (\ref{Meq}) is far from being granted. 
	
\end{rmk}
	
%\bigskip
%\bigskip

The main consequence of this theorem concerns the solution of the {\it patch problem} for  Euler and aggregation equations. The solution $\omega$ or $\rho$ of these problems give rise to a (velocity) field $v$ that can be written in complex form. 

For the case  $\omega(t)*(\frac{x^\perp}{|x|^2})$, this is 
$$
v(t)=\omega(t)*\biggl(\frac{-y+ix}{|z|^2}\biggr)=(i\, \omega(t))*\biggl(\frac{1}{\bar z}\biggr)=\pi\, \bar C[i\, \omega(t)],
$$ 
and  for $\rho(t)*(\frac{x}{|x|^2})$ is 
$$
v(t)=\rho(t)*\biggl(\frac{x+iy}{|z|^2}\biggr)= \rho(t)*\biggl(\frac{1}{\bar z}\biggr)=\pi\, \bar C[\rho(t)].
$$ 

So $\frac{\partial\, v}{\partial z}=\omega \ \ \text{or}\ \  \rho$, where $\omega$ is the vorticity in Euler's equation ($\Div v=0$) or the density of mass in the aggregation equation ($\rot v=0)$.

In both cases there exists an associated flow, 
$$
\psi(z,t)=z+\int_0^tv(\psi(z,\tau),\, \tau)\, d\tau,
$$ 
 in $\C\times[-T,\, T]$. Moreover

\begin{thm}\label{Mth2} 
Let $\Omega\subset\C$ a bounded domain such that $\partial\Omega\in\ka{1,\gamma}$ for $\gamma\in(0,1)$. If $\psi$ is the flow corresponding to the solution of the equations \eqref{eq(E)}, \eqref{eq(A)} or 
\eqref{eqtildeA} with initial condition $\chi_\Omega$, then the function 
$$
t\rightarrow\psi(z,\, t), \ \ \ z\in\C
$$ In $t$
is in $\ka{\omega}(I)$, where $I$ is the interval of existence of the flow.	
\end{thm}

%\bigskip
%\bigskip
%\bigskip

\subsection{Plan of the paper}

The remaining of the paper is devoted to the proof of the two theorems above in three sections.

%\bigskip

In {\it section \ref{Mth2}} we prove Theorem \ref{Mth2} showing that, since both the vorticity in case of  Euler's equation $(E)$ and the density in the case of the continuity equation~$(A)$, are transported by the corresponding flow, then Theorem~\ref{Mth1}  \label{ProP} applies locally in time, providing a family of homeomorphisms of the whole plane parametrized by the time $t$. This family $\phi(z,\, t)$ is regular away of the boundary of $\Omega$, for any fixed $t$ and depends analytically on $t$ in a neighborhood of $t=0$. In each case, the time derivative of $\phi$ gives rise to a velocity field whose $z$-derivative in the regularity points provides a new vorticity or density that satisfies $(E)$ or $(A)$ respectively in these points. 

Since equation $(E)$ is self-dual, and $(\tilde A)$ is the dual equation of $(A)$, then the new vorticity or density satisfy $(E)$ or $(\tilde A)$ and it satisfies duality $(E)$ or $(A)$ in a weak sense. By uniqueness it is the solution of $(E)$ or $(\tilde A)$. As a consequence, since $\partial\Omega$ is regular, also the velocities coincide, and then the function obtained in {\it section \ref{Mth2}}, is in fact the flow corresponding to $(E)$ and $(\tilde A)$, proving that this flow is analytic in time, locally.

The persistence of the regularity, established in~\cite{Che1}, \cite{CoBe} and \cite{Ser3} in the Euler case and in \cite{BeGaLaVe} in the aggregation case allows the extension of this solution at all values of $t$. Then the uniqueness of the solution, established in~\cite{Yud} (see also [Che1] or \cite{CoBe} in the Euler case and in \cite{BeLaLe} in the aggregation case) shows that the flow is analytic in each case.  

\bigskip

In {\it section 3} we prove Theorem \ref{Mth1} following a standard a priori method based in power series developments for $\psi$ in the equation
$$
\frac{\partial^2\psi}{\partial t\partial z}
\overline{\frac{\partial\psi}{\partial z}}-\frac{\partial^2\psi}{\partial t\partial\bar z} 
\overline{\frac{\partial\psi}{\partial\bar z}}=\varpi_0(z)\, e^{2\, \Re\{\varpi_0(z)\}\, t},
$$ 
where $\varpi_0$ is essentially $i\, \omega_0$ or $\rho_0$, that changes the PDE problem to a system of functional equations.
Our procedure is inspired in the paper~\cite{FrZh}, showing analyticity of flows in a different context.

%\bigskip

{\it Section 4} is devoted to some technical facts. The proof of Theorem~\ref{Mth1} heavily relies on some formulas and precise bounds for the (conjugate) Beurling transform on domains with regular and bounded boundary. This technical result is the statement of Theorem \ref{Mth3} in section 4.  

% !TeX encoding = UTF-8
% !TeX spellcheck = ca_ES
% !TeX root = Analit4.tex 

\section{Proof of Theorem \ref{Mth2}}

First of all we show that we can reduce Euler and aggregation equations, to Theorem \ref{Mth1}, at least locally. 

\begin{lem}\label{jak}	
Let $V\subset\C$ an open subset, $\alpha\colon V\rightarrow(0,\infty)$. Let us consider 
$$
U=\{(z,\, t): z\in V,\, t\in(-\alpha(z),\, \alpha(z))\},
$$ 
and a complex valued function $a$, such that there exist functions $v$ and $\psi$ defined in~$U$ and regular enough such that 
$$
\frac{\partial v}{\partial z}(z,t)=a(z,\, t)
$$ 
and 
\begin{equation}\label{flux1}
\frac{\partial\psi}{\partial t}(z,\, t)=v(\psi(z,t),\, t),
\end{equation}
then 
\begin{multline*}
\biggl(\frac{\partial^2\psi}{\partial t\, \partial z} \overline{\frac{\partial\psi}{\partial z}}-
\frac{\partial^2\psi}{\partial t\, \partial\bar z} \overline{\frac{\partial\psi}{\partial\bar z}}\biggr)(z,\, t)\\
=a(\psi(z,t),t) \biggl(\biggl|\frac{\partial\psi}{\partial z}\biggr|^2-\biggl|\frac{\partial\psi}{\partial\bar z}\biggr|^2\biggr) (z,\, 0)\, e^{2\, \int_0^t\Re\{a(\psi(z,\tau),\tau)\}\, d\tau}.
\end{multline*}
\end{lem}

\begin{proof}
Taking derivatives with respect to $z$ in \eqref{flux1} we have
\begin{equation*}
\begin{split}	
\frac{\partial^2\psi}{\partial t\, \partial z}(z,\, t)&=\frac{\partial v}{\partial\zeta}(\psi(z,\, t),\, t) \frac{\partial\psi}{\partial z}(z,\, t)+\frac{\partial v}{\partial\bar\zeta}(z,\, t)\, \frac{\partial\bar\psi}{\partial z}(z,\, t)\\*[5pt]
&=a(\psi(z,t),t) \frac{\partial\psi}{\partial z}(z,\, t)+\frac{\partial v}{\partial\bar\zeta}(z,\, t) \frac{\partial\bar\psi}{\partial z}(z,\, t)
\end{split}
\end{equation*} 
and multiplying by $\overline{\frac{\partial\psi}{\partial z}}$ we have
\begin{equation}\label{det1}
\frac{\partial^2\psi}{\partial t\, \partial z}(z,\, t) \overline{\frac{\partial\psi}{\partial z}}(z,\, t)=a(\psi(z,t),t)
\biggl|\frac{\partial\psi}{\partial z}\biggr|^2(z,\, t)+\frac{\partial v}{\partial\bar\zeta}(z,\, t) \frac{\partial\bar\psi}{\partial z}(z,\, t)\overline{\frac{\partial\psi}{\partial z}}(z,\, t).
\end{equation} 
	
Also taking derivatives with respect to $\bar z$ in \eqref{flux1}, we have
\begin{equation*}
\begin{split}		
\frac{\partial^2\psi}{\partial t\, \partial\bar z}(z,\, t)&=\frac{\partial v}{\partial\zeta}(\psi(z,\, t),\, t) \frac{\partial\psi}{\partial\bar z}(z,\, t)+\frac{\partial v}{\partial\bar\zeta}(z,\, t) \frac{\partial\bar\psi}{\partial\bar z}(z,\, t)\\*[5pt]
&=a(\psi(z,t),t) \frac{\partial\psi}{\partial\bar z}(z,\, t)+\frac{\partial v}{\partial\bar\zeta}(z,\, t) \frac{\partial\bar\psi}{\partial\bar z}(z,\, t)
\end{split}
\end{equation*} 
and multiplying by $\overline{\frac{\partial\psi}{\partial\bar z}}$ we have
\begin{equation}\label{det2}\frac{\partial^2\psi}{\partial t\, \partial\bar z}(z,\, t) \overline{\frac{\partial\psi}{\partial\bar z}}(z,\, t)=a(\psi(z,t),t)\biggl|\frac{\partial\psi}{\partial\bar z}\biggr|^2(z,\, t)+\frac{\partial v}{\partial\bar\zeta}(z,\, t) \frac{\partial\bar\psi}{\partial\bar z}(z,\, t)\, \overline{\frac{\partial\psi}{\partial\bar z}}(z,\, t).
\end{equation}
	
Substracting \eqref{det2} from \eqref{det1}, and taking advantage of a cancellation we conclude that 
$$
\biggl(\frac{\partial^2\psi}{\partial t\, \partial z} \overline{\frac{\partial\psi}{\partial z}}-
\frac{\partial^2\psi}{\partial t\, \partial\bar z} \overline{\frac{\partial\psi}{\partial\bar z}}\biggr)(z,\, t)=a(\psi(z,t),t)\biggl(\biggl|\frac{\partial\psi}{\partial z}\biggr|^2-\biggl|\frac{\partial\psi}{\partial\bar z}\biggr|^2\biggr).
$$
	
Moreover, we have 
\begin{equation*}
\begin{split}
\frac{\partial}{\partial t}\biggl(\biggl|\frac{\partial\psi}{\partial z}\biggr|^2-\biggl|\frac{\partial\psi}{\partial\bar z}\biggr|^2\biggr)(z,\, t)&=2\, \Re\biggl\{\biggl(\frac{\partial^2\psi}{\partial t\, \partial z} \overline{\frac{\partial\psi}{\partial z}}-
\frac{\partial^2\psi}{\partial t\, \partial\bar z} \overline{\frac{\partial\psi}{\partial\bar z}}\biggr)(z,\, t)\biggr\}\\*[5pt]
&=2\, \Re\{a(\psi(z,t),t)\} \biggl(\biggl|\frac{\partial\psi}{\partial z}\biggr|^2-\biggl|\frac{\partial\psi}{\partial\bar z}\biggr|^2\biggr)(z,\, t),
\end{split}
\end{equation*}
and if $|\frac{\partial\psi}{\partial z}|^2-|\frac{\partial\psi}{\partial\bar z}|^2)(z,\, t)$ never vanishes on $U$, then 
$$
\biggl(\biggl|\frac{\partial\psi}{\partial z}\biggr|^2-\biggl|\frac{\partial\psi}{\partial\bar z}\biggr|^2\biggr)(z,\, t)=
\biggl(\biggl|\frac{\partial\psi}{\partial z}\biggr|^2-\biggl|\frac{\partial\psi}{\partial\bar z}\biggr|^2\biggr)(z,\, 0)\, e^{2\, \int_0^t\Re\{a(\psi(z,\tau),\tau)\}\, d\tau}.\rlap{\qedhere}
$$
\end{proof}

%\bigskip
%\bigskip	

%\begin{rmk}
%The jacobian of $\psi$ as considered a map from $\R^2$ into itself is, in complex %coordinates, 
%$$
%J_z\psi(z,\, t)=\biggl(\biggl|\frac{\partial\psi}{\partial %z}\biggr|^2-\biggl|\frac{\partial\psi}{\partial\bar z}\biggr|^2\biggr)(z,\, t).
%$$ 
%\end{rmk}

%\bigskip
%\bigskip

Now

\begin{enumerate}	
\item	In the case of {\it the vortex patch problem}, we start by considering, for the case of {\it the Euler equation}, the {\bf purely imaginary valued} function defined on~$\C$ by 
$$
\varpi_0=i\, c\, \biggl\{\chi_\Omega+\frac{1}{2}\, \chi_{\partial\Omega}\biggr\},
$$ 
where $c\in\R$. 
	
Then 
$$
v_0(z)=\bar C[\varpi_0](z)
$$ 
and by Yudovich's theorem there exist functions $\varpi$ and $v$ defined in $\C\times\R$ such that $\varpi\in 
L^{\infty}(\C\times\R)\cap L^{\infty}_{\text{loc}}(\R; L^p(\C))$  for $1<p<\infty$ and takes its values in $i\, \R$, and $v$ satisfies 
$$
\frac{\partial v}{\partial z}(\ ,t)=\frac{1}{2}\, \varpi(\ ,t)
$$ 
in the distributions sense, and the couple $(\varpi,\, v)$ satisfy the Euler equation~\eqref{eq(E)} in the weak sense. 

%\bigskip

Moreover, there exists a unique function $\psi\in\ka{}(\C\times\R;\C)$ such that 
$$
\psi(z,t)=z+\int_0^t v(\psi(z,\tau),\tau)\, d\tau,
$$ 
and there is a constant $C>0$ such that for any $t\in\R$, 
$$
\psi(\ ,t)-I\in\ka{e^{-Ct\|\varpi_0\|_{(L^p\cap L^\infty)(\C)}}}(\C).
$$

%\bigskip

Then, from the particular shape of $\varpi_0$, we have that if $U\subset\Omega$ or $U$ is a bounded subset of $\C\setminus\bar\Omega$, then $v_0\in\ka{\infty}(U)$  and then, using for instance Proposition~8.3 in \cite{MaBe},  we can conclude that for any $t\in\R$, %\newline
$
\varpi(\ ,t)\in\ka{\infty}(\psi(U,\, t))
$
and 
$
v(\ ,t)\in\ka{\infty}(\psi(U,\, t)).
$ 

Moreover, both the velocity and the {flow} given by the function $\psi$ in the theorem are globally defined with respect to the time variable and in general are regular beyond the continuity in the $z$ variable.

%\bigskip

The flow $\psi$ also inherites the local space regularity (after derivation under the integral sign). Moreover, using 
Theorem~1.3.1 of Chapter~1 in \cite{Hor}, we have that $\psi\in\ka{1}(U\times [0,\infty))$ and $\frac{\partial\psi}{\partial t}(\ ,t)\in\ka{1}(U)$.

%\medskip 

Then the incompressibility of the fluid can be written in terms of the jacobian of $\psi$ as
$$J_z\psi=
\biggl|\frac{\partial\psi}{\partial z}\biggr|^2-\biggl|\frac{\partial\psi}{\partial\bar z}\biggr|^2\equiv 1,
$$ 
also the vorticity is constant along the flow lines.  

From these facts and Lemma \ref{jak} we get the relationship between the flow and the vorticity 
\begin{equation}\label{PE}
\biggl(\frac{\partial^2\psi}{\partial t\partial z}
\overline{\frac{\partial\psi}{\partial z}}-\frac{\partial^2\psi}{\partial t\partial\bar z}
\overline{\frac{\partial\psi}{\partial\bar z}}\biggr)(z,t)=\frac{1}{2}\, \varpi_0(z),
\end{equation}
in the region of $\C\times\R$ where it makes sense. 

%\bigskip

The formula above falls in the so called lagrangian approach. This approach was introduced by A.~Cauchy in~\cite{Cau} and used by several authors since then (e.\ g.\ \cite{FrZh} or \cite{Shn}). 

%\bigskip
%\bigskip
%\bigskip
%\bigskip

\item For the case of {\it the aggregation equation}, we consider the {\bf real valued} function defined on $\C$ by  
$$
\varpi_0=c\biggl\{\chi_\Omega+\frac{1}{2}\, \chi_{\partial\Omega}\biggr\},
$$ 
where $c\in\R$. 

Again 
$$
v_0(z)=\bar C[\varpi_0](z)
$$ 
and, as proven in \cite[Theorems 2.3,  2.4 and 3.1]{BeLaLe}, there exists a constant $T=T(c)$ and functions $\varpi$ and $v$ defined in $\C\times[0,T)$ such that $\varpi\in 
L^{\infty}(\C\times[0,T))\cap \ka{}([0,\, T); L^1(\C))$, also $\varpi(\ , t)$ has bounded support for each $t\in[0,\, T)$, and the function 
$$
\frac{\partial v}{\partial z}(z,t)=\frac{1}{2}\, \varpi(z,t)
$$ 
in weak sense. The functions $\varpi$ and $v$ are unique, solving the equation 
\begin{equation}\label{eq(A)bis}
\begin{cases}
\frac{\partial\varpi}{\partial t}+2\, \Re(\frac{\partial(\varpi\, v)}{\partial z})=0,\\
\varpi(\ , 0)=\varpi_0.
\end{cases}
\tag{$A$}
\end{equation}

%\medskip

All this implies that $v_0\in \Lip{(1,\C)}\cap\ka{\infty}(\C\setminus\partial\Omega),$ and since 
$$
v(z,t)=\bar C[\varpi(\ , t)](z)
$$ 
also (cf.\ \cite{AsIwMa})  
$$
v\in  L^{\infty}(\C\times[0,T))\cap \ka{}([0,\, T); \Lip{(1,\,\C)}).
$$

Then (cf.\ \cite[Theorem~5.2.1]{Che}), there exists a unique function $\psi\in\ka{}(\C\times[0,T);\C)$ such that 
$$
\psi(z,t)=z+\int_0^t v(\psi(z,\tau),\tau)\, d\tau,
$$ 
and  $(\varpi,\, v)$ and there is a constant $C>0$ such that for any $t\in\R$, 
$$
\psi(\ ,t)-I\in\ka{e^{-Ct\|w_0\|_{(L^p\cap L^\infty)(\C)}}}(\C).
$$

It is then clear that 
$$
\biggl(\biggl|\frac{\partial\psi}{\partial z}\biggr|^2-\biggl|\frac{\partial\psi}{\partial\bar z}\biggr|^2\biggr)(z,\, 0)=1.
$$

\bigskip

As in the case of the Euler equation, for each $t\in[0,T)$ and $U\subset\C\setminus\partial\Omega$ and bounded, we have that $\psi(U,\, t)\subset\C\setminus\psi(\partial\Omega,\, t)$ and bounded. Then, from the particular shape of $\varpi_0$, we have that if $U\subset\Omega$ or $U\subset\C\setminus\bar\Omega$ and bounded, then $v_0\in\ka{\infty}(U)$  and, using Proposition~8.3 in \cite{MaBe},  we can conclude that for any $t\in[0,\, T)$, 
$
\varpi(\ ,t)\in\ka{\infty}(\psi(U,\, t))
$
and 
$
v(\ ,t)\in\ka{\infty}(\psi(U,\, t)).
$ 

Moreover, both the velocity and the {flow} given by the function $\psi$ in the theorem are globally defined with respect to the time variable and in general are regular beyond the continuity in the $z$ variable. This is because the local regularity is improved by the Cauchy transform, and so the regularity of initial conditions is propagated by solutions of ordinary differential equations.

%\bigskip

The flow $\psi$ also inherites the local space regularity (after derivation under the integral sign). Moreover, using 
Theorem~1.3.1 of Chapter~1 in \cite{Hor}, we have that $\psi\in\ka{1}(U\times[0, T))$ and $\frac{\partial\psi}{\partial t}(\ ,t)\in\ka{1}(U)$.

%\medskip 

Then we have that (cf.\ \cite{BeLaLe})   
 the transport of the density by the flow satisfies 
 $$
 \varpi(\psi(z,t),t)=\frac{\varpi_0(z)}{1-\varpi_0(z)\, t}
 $$ 
 in $\C\setminus\partial\Omega$. 	
 
These facts allow us to formulate the relationship between the flow and the density, in the region of $\C\times[0,T)$ where it makes sense, as 
\begin{equation}\label{PA}
\begin{split}
\biggl(\frac{\partial^2\psi}{\partial t\partial z}
\overline{\frac{\partial\psi}{\partial z}}-\frac{\partial^2\psi}{\partial t\partial\bar z}
\overline{\frac{\partial\psi}{\partial\bar z}}\biggr)(z,t)&=\frac{1}{2}\frac{\varpi_0(z)}{1-\varpi_0(z)\, t}\, e^{2\, \int_0^t\frac{\varpi_0(z)}{1-\varpi_0(z)\, \tau}\, d\tau}\\*[5pt] 
&=\frac{c}{2\, (1-c\, t)^3}\biggl\{\chi_\Omega+\frac{1}{2}\, \chi_{\partial\Omega}\biggr\}.
\end{split}
\end{equation}
\bigskip
%\bigskip

Let $\phi$ be the solution of the equation \eqref{PE} or \eqref{PA} provided by Theorem \ref{Mth1}. We have that $\phi(z,\, )\in\ka{\omega}(-T_0,\, T_0)$ for $z\in\C$.

Moreover

\begin{prop}\label{ProP} The solution $\phi$ of the equation \eqref{PE} or \eqref{PA} obtained using Theorem \ref{Mth1} satisfies, for $t\in[-T_0,\, T_0]$, that
	
	\begin{enumerate}
		\item $\phi(\ ,\, t)\colon\C\rightarrow\C$ is an homeomorphism.
		
		\medskip
		
		\item $\phi(\ ,\, t)\in\ka{\infty}(\C\setminus\partial\Omega)$.
	\end{enumerate}
	\end{prop}  

\medskip

This result will be proved in section 3.2. 

\bigskip
%\bigskip

 As we have seen, the (unique) flow $\psi$ corresponding to the equation \eqref{eq(E)} or \eqref{eq(A)} with the initial condition $\omega_0(z)$ or $\rho_0(z)$, also satisfies equation \eqref{PE} or \eqref{PA}. Now we prove that $\phi(\ ,\, t)\equiv
\psi(\ ,\, t)$ for $t\in[-T_0,\, T_0]$:

We have that $\psi$ and $\phi$ generate vector fields $v$ and $u$, namely
$$v(\zeta,\, t)=\frac{\partial\psi}{\partial t}(\psi(\ ,\, t)^{-1}(\zeta),\, t)$$ and  
$$u(\zeta,\, t)\stackrel{\text{def}}{=}\frac{\partial\phi}{\partial t}(\phi(\ ,\, t)^{-1}(\zeta),\, t).$$

On the other hand, we have that 
$$\frac{\partial v}{\partial\zeta}(\zeta,\, t)=\varpi(\zeta,\, t)=\frac{1}{2}\, \varpi_0(\psi(\ ,\, t)^{-1}(\zeta))\, b(t)$$ and  $$\frac{\partial u}{\partial\zeta}(\zeta,\, t)\stackrel{\text{def}}{=}\small{\Xi}(\zeta,\, t)=\frac{1}{2}\, \varpi_0(\phi(\ ,\, t)^{-1}(\zeta))\, b(t)$$  in $\C\setminus\partial\Omega$, where $b(t)=1$ in the case of equation~\eqref{eq(E)} and $b(t)=\frac{1}{(1-c\, t)^3}$ in the case of equation~\eqref{eq(A)}.

Now, we will prove that $\small{\Xi}$ is a weak solution of equation~\eqref{eq(E)} or equation~\eqref{eq(A)}, and then, by uniqueness $\varpi\equiv\small{\Xi}$ in the weak sense, so $\varpi=\small{\Xi}$ almost everywhere.

Then $$\frac{\partial (u-v)}{\partial\zeta}(\zeta,\, t)=0$$ for $\zeta$ outside a subset of $\C$ whose  continuous analytic capacity (see \cite[pg~38]{Tol}, also \cite{Gam} or \cite{Gar}) is equal to $0$. Since $u$ and $v$ are bounded and vanish at infinity, by Liouville's theorem we conclude that $u\equiv v$. That  concludes the argument.

\medskip

Let us prove that $\small{\Xi}$ is a weak solution of equation~\eqref{eq(E)} or equation~\eqref{eq(A)}. It is worth to remark that the equation $(E)$ is self-dual, but the dual equation of $(A)$ is $(\tilde A)$. 

\newpage
For $\varphi\in\ka{1}([0,T];\, \ka{1}_0(\C))$
\begin{equation}\begin{split} 
&\int_0^T\, \int_\C \small{\Xi}(\zeta,\, t)\, \frac{D\varphi}{Dt}(\zeta,\, t)\, dm(\zeta)\, dt\\
=&\int_0^T\, \int_\C \varpi_0(\phi(\ ,\, t)^{-1}(\zeta))\, b(t)\, \frac{D\varphi}{Dt}(\zeta,\, t)\, dm(\zeta)\, dt\\
=&\int_0^T\, \int_\C \varpi_0(z)\, b(t)\, \frac{D\varphi}{Dt}(\phi(z,\, t),\, t)\, J_z\phi(z,\, t)\, dm(z)\, dt\\
=&\int_\C \varpi_0(z)\, dm(z)\int_0^T \frac{d\ }{dt}(\varphi(\phi(z,\, t),\, t)\, dt\\
=&\int_\C \varpi_0(z)\, \{\varphi(\phi(z,\, T),\, T)-\varphi(z,\, 0)\}\, dm(z),
\end{split}
\end{equation}
in the second equality we have used that $b(t)\, J_z\phi(z,\, t)\equiv1$.

\bigskip

Finally

\begin{enumerate}
\item In the Euler's case there exists a number $T_1>0$ and a function $\psi\colon\C\rightarrow\C$, such that $t\rightarrow\psi(z,\, t)$ is real analytic in $(-T_1,\, T_1)$.
	
%\medskip
	
If $T_1=+\infty$ there is nothing else to say. Otherwise, at time $\pm T_1$, from the theorem of persistence of regularity (\cite{Che1}, \cite{CoBe}, \cite{Ser3}), we have that $\partial\Omega_{\pm T_1}\in\ka{1,\, \gamma}$ and we can iterate the procedure and use the uniqueness of the solution to obtain $t$-analyticity in $(-T_2,\, T_2)$ for~$T_2>T_1$. This implies that the set of analyticity is open and closed in $\R$, so it is $\R$.
	
%\bigskip
	
\item For the equations~\ref{eq(A)} and~\ref{eqtildeA}, the procedure is similar. We only have to take in account that now $t>0$. For the equation~\ref{eqtildeA}, and for~\ref{eq(A)} in the case of $c<0$, the analyticity will occur for $t\in[0,\, \infty)$. For the equation~\ref{eq(A)} and $c>0$, the analyticity will occur for $t\in[0,\, \frac{1}{c})$. The only necessary ingredient in the proof is the persistence of the regularity.  
	
%\medskip
	
In fact it is enough to have persistence in the case~\eqref{eq(A)bis}. The rescaling 
$$
s(t)=\ln\biggl(\frac{1}{1-ct}\biggr)
$$ 
and 
$$
\tilde\varpi(z,\, s)=\frac{1-c\, t(s)}{c}\, \varpi(z,\, t(s))
$$ 
transforms the problem 
\begin{equation}
\begin{cases}
\frac{\partial\varpi}{\partial t}+2\Re(\frac{\partial(\varpi\, v)}{\partial z})=0,\\
v(z,\, t)=-\bar C[\varpi(\ , t)](z),\\
\varpi(\ , 0)=\varpi_0,
\end{cases}
\tag{$A$}
\end{equation}
in the problem 
\begin{equation}\label{eqtildeAbis}
\begin{cases}
\frac{\partial\tilde\varpi}{\partial s}+\Re(\tilde v\, \frac{\partial\tilde\varpi}{\partial z})=0,\\
\tilde v(z,\, s)=-\bar C[\tilde\varpi(\ , s)](z),\\
\tilde\varpi(\ , 0)=\chi_\Omega,
\end{cases}
\tag{$\tilde A$}
\end{equation}
which is a transport equation with initial datum the indicator function of $\Omega$, a region with $\ka{1,\gamma}$ boundary.  

For the problem $(\tilde A)$ it is proven in \cite{BeGaLaVe} that for every %\newline 
$s\!\in\![0,\infty)$, if $\tilde\psi$ is the corresponding flow, then $\tilde\psi(\partial\Omega,\, s)$ is a $\ka{1,\gamma}$ embedded submanifold of $\C$ of real dimension $1$. Since the rescalings above do not affect the $z$ variable, the same regularity is true in the case of~\eqref{eq(A)bis}, for $t\in[0,\frac{1}{c})$.

%\medskip

From Theorem \ref{Mth1} there exists a number $T_1>0$ and a function $\psi\colon\C\rightarrow\C$, such that $t\rightarrow\psi(z,\, t)$ is real analytic in $(-T_1,\, T_1)$.

%\medskip

If $T_1=+\infty$ there is nothing else to say. Otherwise, after a time $T_1$, we have 
$$
\varpi(z,\, T_1)=\frac{c}{(1-c\, T_1)^3}\, \varpi_0(z),
$$
that must be used as initial density to iterate the procedure, because the boundary $\partial\Omega_{T_1}$ is $\ka{1,\, \gamma}$, getting a new $T_2>T_1$ and analyticity in~$(-T_1,\, T_2)$.

Again an argument of connectivity and the uniqueness conclude that the flow is analytic in $[0,\frac{1}{c})$.    
\end{enumerate}  
\end{enumerate}

%\bigskip

% !TeX encoding = UTF-8
% !TeX spellcheck = ca_ES
% !TeX root = Analit4.tex 
\section{Proof of Theorem \ref{Mth1} and consequences} We divide this section in two parts.

\subsection{Construction of the solution} 

The analyticity in $t$ of $a$ implies the local existence of functions  
$$
a^{(s)}\colon\C\rightarrow\C,
$$ 
for each $s\in\N\setminus\{0\}$ such that 
$$
\varpi_0(z)\, (1+a(z,\, t)\, t)=\varpi_0(z)+\sum_{s=1}^\infty a^{(s)}(z)\, t^s,
$$ 
for $t\in(-T_0,\, T_0)$. 

It is worth to remark that $a^{(s)}(z)=0$ for every $z\in\C\setminus\bar\Omega$.

%\bigskip

Assume that there exists a family of functions 
$$
\xi^{(s)}\colon\C\rightarrow\C,\quad s\in\N,
$$ 
having first order derivatives with respect to $z$ and $\bar z$ at each point of $\C\setminus\partial\Omega$ and such that 
$$
\psi(z,t)=\sum_{s=0}^\infty\xi^{(s)}(z)\, t^s,
$$ 
then
\begin{multline*} 
\biggl(\frac{\partial^2\psi}{\partial t\partial z}
\overline{\frac{\partial\psi}{\partial z}}-\frac{\partial^2\psi}{\partial t\partial\bar z}
\overline{\frac{\partial\psi}{\partial\bar z}}\biggr)(z,t)=\sum_{s=0}^\infty \biggl\{\sum_{k=0}^s(k+1) \biggl(\frac{\partial\xi^{(k+1)}}{\partial z}(z)
\overline{\frac{\partial\xi^{(s-k)}}{\partial z}(z)}\\*[5pt]
-\frac{\partial\xi^{(k+1)}}{\partial\bar z}(z) \overline{\frac{\partial\xi^{(s-k)}}{\partial\bar z}(z)}\biggr)\biggr\}\, t^s.
\end{multline*}

%\bigskip
%\bigskip

If $\psi$ is a solution analytic in $t$ in a neighborhood of $\C\times\{0\}$ of the problem~\eqref{Meq}, then 
$$
\xi^{(0)}(z)\equiv z
$$ 
and if $z\in\C\setminus\partial\Omega$, then
\begin{equation}\label{eqfunct1} 
\!\!\begin{cases}\frac{\partial\xi^{(1)}}{\partial z}(z)=\varpi_0(z)\overset{\text{def}}{=}a^{(0)}(z),\\*[5pt]
\frac{\partial\xi^{(s+1)}}{\partial z}(z)\!=\!\frac{a^{(s)}(z)}{s+1}-\frac{1}{s+1} \sum_{k=0}^{s-1}(k\!+\!1)(\frac{\partial\xi^{(k+1)}}{\partial z} 
\overline{\frac{\partial\xi^{(s-k)}}{\partial z}}
-\frac{\partial\xi^{(k+1)}}{\partial\bar z} \overline{\frac{\partial\xi^{(s-k)}}{\partial\bar z}})(z),\!\!\!\!\!\\*[5pt] 
s\in\N\setminus\{0\}.
\end{cases}\!\!\!\!
\end{equation}

%\bigskip
%\bigskip

The system $\eqref{eqfunct1}$ provides for each $s\in\N\setminus\{0\}$ the $z$-derivative of $\xi^{(s)}$ in terms of the first order derivatives of the functions $\xi^{(l)}$, where $1\le l<s$, if $s>1$, and in terms of $\varpi_0$ in the case $s=1$.

%\medskip

Let us choose  
$$
\xi^{(1)}(z)=\bar C[\varpi_0](z).
$$ 

We have

\begin{prop} 
The function $\bar C[\varpi_0](z)$ is in $\ka{\infty}(\C\setminus \bar\Omega)\cap \Lip{(1,\, \C)}$ and has a decay at the infinity of the type 
$$
\frac{1}{\max\{R_0,\, d(z,\partial\Omega)\}}.
$$

%\bigskip

Moreover, if $T^{(1)}$ is the distribution corresponding to $\frac{\partial\xi^{(1)}}{\partial z}$, then 
$$
\Theta(T^{(1)},\, z)=\varpi_0(z)
$$ 
for all $z\in\C$. 
\end{prop}

\begin{proof} Of the last part
\begin{equation*}
\begin{split}
\biggl\langle \frac{\partial T_{\bar C[\varpi_0]}}{\partial z},\phi_{z_0,\epsilon}\biggr\rangle&=-
\biggl\langle T_{\bar C[\varpi_0]},\frac{\partial \phi_{z_0,\epsilon}}{\partial z}\biggr\rangle\\*[5pt]
&=\frac{i}{2\pi}\int_\C\biggl\{\int_\C\frac{\varpi_0(\zeta)}{\bar\zeta-\bar z}\, dm_2(\zeta)\biggr\} 
	\frac{\partial \phi_{z_0,\epsilon}}{\partial z}(z)\, dm_2(z)\\*[5pt]
&=\frac{i}{2\pi}\int_\C\biggl\{\int_\C\frac{\frac{\partial \phi_{z_0,\epsilon}}{\partial z}(z)}{\bar\zeta-\bar z}\, 
	dm_2(z)\biggr\}\, \varpi_0(\zeta)\, dm_2(\zeta)\\*[5pt]
&=\int_\C\phi_{z_0,\epsilon}(\zeta)\, \varpi_0(\zeta)\, dm_2(\zeta).\qedhere
\end{split}
\end{equation*}
\end{proof}

%\bigskip
%\bigskip

In the remaining we use the following theorem on the boundedness of the conjugate Beurling transform

\begin{thm}\label{Mth3} 
Let $\Omega$ be a domain such that $\partial\Omega\in\ka{1,\gamma}$, where $\gamma\in(0,1)$ and $g$ a function defined at every $z\in\C$, $g\in\ka{\infty}(\C\setminus\partial\Omega)$ and $\chi_\Omega g$ extends to a (unique) function $g_-\in \Lip{(\gamma,\, \bar\Omega)}$, $\chi_{\bar\Omega^c} g$ extends to a (unique) function $g_+\in \Lip{(\gamma,\, \C\setminus\Omega)}$ and for $z\in\partial\Omega$, $g(z)=\frac{1}{2}\, \{g_+(z)+g_-(z)\}$.

%\medskip

We also assume that there is a constant $R_0>0$ depending only on $\Omega$ such that if 
$$
U_{R_0}'=\{z\in\C: d(z,\Omega)<R_0\},
$$ 
then  for $g$ there is a constant $C(g)>0$ such that for $z\in\C\setminus U_{\R_0}'$, 
$$
|g(z)|\le\frac{C(g)}{\max\{R_0^2,\, d(z,\partial\Omega)^2\}}.
$$ 

Then  $\bar B[g](z)$ is well defined for each $z\in\C$ and 
\begin{align*}
\chi_\Omega\, \bar B[g]\in\ka{\infty}(\Omega)&\cap \Lip{(\gamma,\, \bar\Omega)},\\
\chi_{\C\setminus \bar\Omega}\, \bar  B[g]\in\ka{\infty}(\C\setminus \bar\Omega)&\cap \Lip{(\gamma,\, \C\setminus \Omega)},
\end{align*}
and for any $z\in\partial\Omega$ we have 
\begin{equation}\label{jump}
\bar B[g](z)=\frac{1}{2}\biggl\{\lim_{w\to z;\, w\in\Omega} \chi_\Omega\,  \bar B[g](w)+\lim_{w\bar\to z;\, w\in\C\setminus\Omega} \chi_{\C\setminus \bar\Omega}\, \bar B[g](w)\biggr\}.
\end{equation}

\medskip

Moreover, there exists a constant $K=K(\gamma,\, \Omega,\, R_0)>0$ such that  
\begin{align*}
\|\bar B[g]\|_{L^\infty(\C)}&\le K\, \|g\|_\gamma,\\*[5pt]
\|\bar B[g]\|_{\gamma,\, \bar\Omega} &\le K\, \|g\|_\gamma,
\intertext{and} 
\|\bar B[g]\|_{\gamma,\, \C\setminus\Omega}&\le K\, \|g\|_\gamma.
\end{align*}

%\medskip

Moreover, for any $z\in (U_{R_0}\cup\Omega)^c$, 
\begin{equation} \label{decay}
 |\bar B[g](z)|\le K\, (1+C(g))\, \frac{1+\ln d(z,\, \Omega)}{\max\{R_0^2,\, d(z,\partial\Omega)^2\}}\, \|g\|_\gamma.
 \end{equation} 
\end{thm}

%
%\bigskip
%\bigskip

Using Theorem \ref{Mth3} we have that

\begin{prop} If $\tilde T^{(1)}$ is the distribution corresponding to $\frac{\partial\xi^{(1)}}{\partial\bar z}$, then 
$$
\Theta(\tilde T^{(1)},\, z)=\bar B[\varpi_0](z)
$$ 
for all $z\in\C$.

%\bigskip

So $\tilde T^{(1)}$ is given by a function whose decay at $\infty$ is 
of type 
\begin{equation}\label{decay1}
\frac{1+\ln d(z,\, \Omega)}{\max\{R_0^2,\, d(z,\partial\Omega)^2\}}.
\end{equation}	
\end{prop}

\begin{proof} 
The density part.  
\begin{equation*}
\begin{split}
\left\langle \frac{\partial T_{\bar C[\varpi_0]}}{\partial\bar z},\phi_{z_0,\epsilon}\right\rangle &=-\left\langle T_{\bar C[\varpi_0]},\frac{\partial \phi_{z_0,\epsilon}}{\partial\bar z}\right\rangle\\*[5pt]
&=\frac{i}{2\pi}\int_\C\biggl\{\int_\C\frac{\varpi_0(\zeta)}{\bar\zeta-\bar z}\, dm_2(\zeta)\biggr\}
\frac{\partial \phi_{z_0,\epsilon}}{\partial\bar z}(z)\, dm_2(z)\\*[5pt]
&=\frac{i}{2\pi}\int_\C\biggl\{\int_\C\frac{\frac{\partial \phi_{z_0,\epsilon}}{\partial \bar z}(z)}{\bar\zeta-\bar z}\, 
dm_2(z)\biggr\}\, \varpi_0(\zeta)\, dm_2(\zeta)=(*)
\end{split}
\end{equation*}
and using the Stokes formula and Theorem \ref{Mth3}, we have 
$$
\int_\C\frac{\frac{\partial \phi_{z_0,\epsilon}}{\partial\bar z}(z)}{\bar\zeta-\bar z}\, 
dm_2(z)=\lim_{\epsilon\to0}\int_{\C\setminus B_\epsilon(\zeta)}\frac{\frac{\partial \phi_{z_0,\epsilon}}{\partial\bar z}(z)}{\bar\zeta-\bar z}\, 
dm_2(z)=(**)
$$ 
and 
\begin{equation*}
\begin{split}
\int_{\C\setminus B_\epsilon(\zeta)}\!\frac{\frac{\partial \phi_{z_0,\epsilon}}{\partial\bar z}(z)}{\bar\zeta-\bar z}\, 
dm_2(z)&=\frac{1}{2\, i} \int_{\C\setminus B_\epsilon(\zeta)}\frac{\frac{\partial \phi_{z_0,\epsilon}}{\partial\bar z}(z)}{\bar\zeta-\bar z}\, d\bar z\wedge dz\\*[5pt]
&=\frac{1}{2\, i} \int_{\partial B_\epsilon(\zeta)} \!\frac{\phi_{z_0,\epsilon}(z)}{\bar\zeta-\bar z}\, 
dz-\frac{1}{2\, i} \int_{\C\setminus B_\epsilon(\zeta)}\!\frac{\phi_{z_0,\epsilon}(z)}{(\bar\zeta-\bar z)^2}\, d\bar z\wedge dz,\!\!\!
\end{split}
\end{equation*}
so 
$$
(**)=-\frac{1}{2\, i} \lim_{\epsilon\to0} \int_{\C\setminus B_\epsilon(\zeta)}\frac{\phi_{z_0,\epsilon}(z)}{(\bar\zeta-\bar z)^2}\, d\bar z\wedge dz,
$$ 
and then 
$$
(*)=\int_\C\phi_{z_0,\epsilon}(\zeta)\{\frac{1}{\pi}\text{ p.\ v.\ }  \int_\C \frac{\varpi_0(\zeta)}{(\bar\zeta-\bar z)^2}\,  dm_2(\zeta)\}\, dm(z).
$$

\bigskip

The decay estimate \eqref{decay1} is a consequence of \eqref{decay} in Theorem \ref{Mth3}. 				
\end{proof}

%\bigskip
%\bigskip

Now, for $s>1$ we will define $T^{(s)}$ and  $\tilde T^{(s)}$, the distributions corresponding to $\frac{\partial\xi^{(s)}}{\partial z}$
and $\frac{\partial\xi^{(s)}}{\partial\bar z}$ respectively, and then $\Theta(T^{(s)},\, z)=\theta^{(s)}(z)$ and $\Theta(\tilde T^{(s)},\, z)=\eta^{(s)}(z)$. All these objects exist, by a direct application of the previous proposition to the cases $s'<s$.

%\medskip

So, taking densities, the system \eqref{eqfunct1}  becomes
\begin{equation}\label{eqfunct2} 
\begin{cases}\theta^{(1)}(z)=\varpi_0(z),\\*[5pt]
\eta^{(1)}(z)=\bar B[\varpi_0](z),\\*[5pt]
\theta^{(s+1)}(z)=\frac{a^{(s)}(z)}{s+1}-\frac{1}{s+1}\sum_{k=0}^{s-1}(k+1)\\*[5pt]
\hphantom{\theta^{(s+1)}(z)=}
\times \{\theta^{(k+1)} 
\overline{\theta^{(s-k)}}
-\eta^{(k+1)} \overline{\eta^{(s-k)}}\}(z), & \text{for }  s\in\N\setminus\{0\}.
\end{cases}
\end{equation}

%\bigskip
%\bigskip

Then we consider,  in $\C$, the decomposition 
\begin{equation} \label{decomp}
\theta^{(j)}=\chi_{\Omega}\, \theta^{(j)}+\chi_{\bar\Omega^c}\, \theta^{(j)}+\chi_{\partial\Omega}\, \theta^{(j)}=\phi^{(j)}+\psi^{(j)}+
\i^{(j)}.
\end{equation}

So 
\begin{equation*}
\begin{split}
&(s+1)\, (\phi^{(s+1)}+\psi^{(s+1)}+\i^{(s+1)})(z)\\*[5pt]
&\quad=a^{(s)}(z)
-\!\sum_{k=0}^{s-1}(k+1) \Bigl\{(\phi^{(k+1)}+\psi^{(k+1)}+\i^{(k+1)}) (\overline{\phi^{(s-k)}+\psi^{(s-k)}}+\i^{(k+1)})\!\\*[5pt]
&\hspace*{2.6cm}-\bar B[\phi^{(k+1)}+\psi^{(k+1)}+\i^{(s-k)}] \overline{\bar B[\phi^{(s-k)}+\psi^{(s-k)}+\i^{(s-k)}]}\Bigr\}(z)\\*[5pt]
&\quad=a^{(s)}(z)-\!\sum_{k=0}^{s-1}(k+1) \Bigl\{(\phi^{(k+1)}+\psi^{(k+1)}+\i^{(k+1)}) (\overline{\phi^{(s-k)}+\psi^{(s-k)}+\i^{(k+1)}})\\*[5pt]
&\hspace*{5.2cm}-\bar B[\phi^{(k+1)}+\psi^{(k+1)}] \overline{\bar B[\phi^{(s-k)}+\psi^{(s-k)}]}\Bigr\}(z),
\end{split}
\end{equation*}
consequently, if we define 
$$
\Phi[\ ] =\chi_\Omega\, \bar B[\ ],\quad  
\Psi[\ ] =\chi_{\C\setminus\Omega}\, \bar B[\ ],\quad  \text{and}\quad 
\Gamma[\ ] =\chi_{\partial\Omega}\, \bar B[\ ],$$ we have
\begin{equation*}
\begin{split}
(s+1)\, \phi^{(s+1)}&=a^{(s)}-\sum_{k=0}^{s-1}(k+1)\Bigl \{\phi^{(k+1)}\overline{\phi^{(s-k)}}-\Phi[\phi^{(k+1)}]\overline{\Phi[\phi^{(s-k)}]}\Bigr\}\\*[5pt]
&\quad+\sum_{k=0}^{s-1}(k+1)\Bigl \{\Phi[\phi^{(k+1)}] \overline{\Phi[\psi^{(s-k)}]}+\Phi[\psi^{(k+1)}] \overline{\Phi[\phi^{(s-k)}]}\\*[5pt]
&\hspace*{5.4cm}+\Phi[\psi^{(k+1)}] \overline{\Phi[\psi^{(s-k)}]}\Bigr\},
 \end{split}
\end{equation*}
in $\Omega$, 
\begin{equation*}
\begin{split}
(s+1)\, \psi^{(s+1)}&=- \sum_{k=0}^{s-1}(k+1)\Bigl\{\psi^{(k+1)}) \overline{\psi^{(s-k)}}-\Psi[\psi^{(k+1)}] \overline{\Psi[\psi^{(s-k)}]}\Bigr\}\\*[5pt]
&\quad+\sum_{k=0}^{s-1}(k+1)\Bigl\{\Psi[\phi^{(k+1)}] \overline{\Psi[\psi^{(s-k)}]}+\Psi[\psi^{(k+1)}] \overline{\Psi[\phi^{(s-k)}]}\\*[5pt]
&\hspace*{5.1cm}+\Psi[\phi^{(k+1)}] \overline{\Psi[\phi^{(s-k)}]}\Bigr\}
 \end{split}
\end{equation*}
in $\Omega^c$, and
\begin{equation*}
\begin{split}
(s+1)\, \i^{(s+1)}&=a^{(s)}-\sum_{k=0}^{s-1}(k+1)\Bigl \{\i^{(k+1)}\overline{\i^{(s-k)}}-\Gamma[\phi^{(k+1)}] \overline{\Gamma[\phi^{(s-k)}]}\Bigr\}\\*[5pt]
&\quad+\sum_{k=0}^{s-1}(k+1)\Bigl \{\Gamma[\phi^{(k+1)}] \overline{\Gamma[\psi^{(s-k)}]}+\Gamma[\psi^{(k+1)}] \overline{\Gamma[\phi^{(s-k)}]}\\*[5pt]
&\hspace*{5.4cm}+\Gamma[\psi^{(k+1)}] \overline{\Gamma[\psi^{(s-k)}]}\Bigr\},
 \end{split}
\end{equation*}
in $\partial\Omega$.

%\bigskip
%\bigskip

And now we can determine inductively the functions $\phi$, $\psi$ and $\i$.

\begin{prop} 
If $\phi^{(1)}=\chi_\Omega\, \varpi_0$, $\psi^{(1)}=0$ and $\i^{(1)}=\chi_{\partial\Omega}\, \varpi_0$, then for any $s\in\N\setminus\{0,\, 1\}$, we have that $\phi^{(s)}\in \Lip{(\gamma,\, \bar\Omega)}$, $\psi^{(s)}\in \Lip{(\gamma,\, \Omega^c)}$ and $\i^{(s)}\in \Lip{(\gamma,\, \partial\Omega)}$.  
\end{prop}

\begin{proof}
The formulas above allow, for any $s\geq 2$, the control of $\phi^{(s+1)}$ and  $\psi^{(s+1)}$ in terms of $\phi^{(l)}$ and  $\psi^{(r)}$, for all $1\le r,l\le s$. 

In the case of $s=1$, Theorem \ref{Mth3} gives the desired estimate. 

For $s>1$, we need the existence and estimates of terms of type 
$$
B[B[f^{(l)}]\, B[f^{(r)}]]
$$ 
for $1\le r,l\le s$, where $f^{(s)}$ states for any of the functions in the statement.

The estimates in $\{d(z,\, \partial\Omega)>1\}$ provided by \eqref{decay} in Theorem \ref{Mth3} imply that 
$$
B[f^{(l)}]\, B[f^{(r)}]=O\biggl(\frac{\ln^2 d(z,\partial\Omega)}{d(z,\partial\Omega)^4}\biggr),
$$ 
and then the product is integrable and satisfies the hypotheses  Theorem \ref{Mth3} with the corresponding estimates and we can perform iteration.

\medskip

Also we have that $$|\psi^{(s)}(z)|\le C \frac{\ln d(z,\partial\Omega)}{\max\{1,\, d(z,\partial\Omega)^2\}},$$ with $C$ independent of $s$. 
\end{proof}
%\bigskip
%\bigskip

Let us control now the H\"older and some related integrability sizes of $\phi$, $\psi$ and $\i$.

%\medskip

If $\|\ \|$ denotes $$\|\ \|_{\gamma,\, \bar\Omega}+\|\ \|_{L^p(\bar\Omega)}+\|\ \|_{L^q(\bar\Omega)}$$ or $$\|\ \|_{\gamma,\, \Omega^c}+\|\ \|_{L^p(\Omega^c)}+\|\ \|_{L^q(\Omega^c)},$$ where $p=\frac{2}{1-\gamma}$ and $q=\frac{2}{1+\gamma}$ is the dual exponent, we have 
$$
\| \phi^{(1)} \|=\| \varpi_0 \|\overset{\text{def}}{=}\frac{1}{2}\, \alpha
$$ 
and 
$$
\| \psi^{(1)} \|=0.
$$

Using that $\|\ \|$ is a multiplicative norm, the remaining terms have the control 
\begin{equation*}
\begin{split}
\|\phi^{(s+1)}\|&\le\frac{1}{s+1}\, \|\varpi_0\|\, \frac{(2\|\varpi_0\|)^s}{s!}\\*[5pt] 
&\quad+\frac{1}{s+1}\,\sum_{k=0}^{s-1}(k+1)\Bigl\{\|\phi^{(k+1)}\|\, \|\phi^{(s-k)}\|+\|\Phi[\phi^{(k+1)}]\| \|\Phi[\phi^{(s-k)}]\|\Bigr\}\\*[5pt]
&\quad+\frac{1}{s+1}\, \sum_{k=0}^{s-1}(k+1)\Bigl \{\|\Phi[\phi^{(k+1)}]\|\,  \|\Phi[\psi^{(s-k)}]\|\\*[5pt]
&\hspace*{2.5cm}+\|\Phi[\psi^{(k+1)}]\|\  \|\Phi[\phi^{(s-k)}]\|
+\|\Phi[\psi^{(k+1)}]\ \|\Phi[\psi^{(s-k)}]\|\Bigr\},\!
 \end{split}
\end{equation*}
and 
\begin{equation*}
\begin{split}
\|\psi^{(s+1)}\|&\le\frac{1}{s+1}\, \sum_{k=0}^{s-1}(k+1)\Bigl \{\|\psi^{(k+1)}\|\ \|\psi^{(s-k)}\|+\|\Psi[\psi^{(k+1)}]\|\ \|\Psi[\psi^{(s-k)}]\|\Bigr\}\\*[5pt]
&\quad+\frac{1}{s+1}\, \sum_{k=0}^{s-1}(k+1)\Bigl \{\|\Psi[\phi^{(k+1)}]\|\ \|\Psi[\psi^{(s-k)}]\|\\*[5pt]
&\hspace*{2cm}+\|\Psi[\psi^{(k+1)}]\|\ \|\Psi[\phi^{(s-k)}]\|
+\|\Psi[\phi^{(k+1)}]\|\ \|\Psi[\phi^{(s-k)}]\|\Bigr\}.
 \end{split}
\end{equation*}

If $f=\phi^{(s)}$ or $\psi^{(s)}$, using Theorem \ref{Mth3}, we have that
$\|\Phi[f]\|$ and $\|\Psi[f]\|\le K\, \|f\|$, and then, using the notation $\alpha_p=\|\phi^{(p)}\|$ and $\beta_p=\|\psi^{(p)}\|$, 
\begin{equation*}
\begin{split}
\alpha_{s+1}&\le\frac{1}{s+1}\, \frac{\alpha}{2}\, \frac{\alpha^s}{s!}+(1+K^2)\, \frac{1}{s+1}\, \sum_{k=0}^{s-1}(k+1)\, \alpha_{k+1}\, \alpha_{s-k}\\*[5pt]
&\quad+K^2\, \frac{1}{s+1}\, \sum_{k=0}^{s-1}(k+1)\, \{\alpha_{k+1}\,  \beta_{s-k}+\beta_{k+1}\,  \alpha_{s-k}+\beta_{k+1}\, \beta_{s-k}\},
 \end{split}
\end{equation*}
and
\begin{equation*}
\begin{split}
\beta_{s+1}&\le(1+K^2)\, \frac{1}{s+1}\, \sum_{k=0}^{s-1}(k+1)\, \beta_{k+1}\, \beta_{s-k}\\*[5pt]
&\quad+K^2\, \frac{1}{s+1}\, \sum_{k=0}^{s-1}(k+1)\, \{\alpha_{k+1}\,  \beta_{s-k}+\beta_{k+1}\,  \alpha_{s-k}+\alpha_{k+1}\,  \alpha_{s-k}\}.
 \end{split}
\end{equation*}

Now, we have that if we define 
\begin{align*}
A_s&=\sum_{k=0}^{s-1}(k+1)\, \alpha_{k+1}\, \alpha_{s-k},\\
B_s&=\sum_{k=0}^{s-1}(k+1)\, \beta_{k+1}\, \beta_{s-k} 
\intertext{and}
C_s&=\sum_{k=0}^{s-1}(k+1)\, \{\alpha_{k+1}\,  \beta_{s-k}+\beta_{k+1}\,  \alpha_{s-k}\},
\end{align*}
then we have 
$$(s+1)\, \alpha_{s+1}\le\frac{1}{2}\, \frac{\alpha^{s+1}}{s!}+(1+K^2)\, A_s+K^2\, (C_s+B_s)$$
and
$$(s+1)\, \beta_{s+1}\le(1+K^2)\, B_s+K^2\, (C_s+A_s).$$

%\bigskip

%\bigskip
%\bigskip

Let us consider, next, the polynomic functions 
\begin{align*}
f_1(\xi)&=\frac{\alpha}{2}\, \xi,\\*[5pt]
g_1(\xi)&=0,\\*[5pt]
f_N(\xi)&=\sum_{p=0}^N\alpha_p\, \xi^p,
\intertext{and}
g_N(\xi)&=\sum_{q=0}^N\beta_q\, \xi^q,
\end{align*} 
where $\alpha_0=\beta_0=0$.

%\bigskip

Then the previous inequalities imply that, for $\xi\in[0,\infty)$, 
\begin{equation*}
\begin{split}
f'_{N+1}(\xi)&=\sum_{s=0}^N(s+1)\, \alpha_{s+1}\, \xi^s\\*[5pt]
&\le\frac{1}{2}\, \sum_{s=0}^N\, \frac{\alpha^{s+1}}{s!}\, \xi^s+(1+K^2)\, \sum_{s=0}^NA_s\, \xi^s+K^2\, \sum_{s=0}^N(C_s+B_s)\, \xi^s
\end{split}
\end{equation*}
and 
$$
g'_{N+1}(\xi)=\sum_{s=0}^N(s+1)\, \beta_{s+1}\, \xi^s\le(1+K^2)\, \sum_{s=0}^NB_s\, \xi^s+K^2\, \sum_{s=0}^N(C_s+A_s)\xi^s.
$$

%\bigskip

Then we have that $$\sum_{s=0}^NA_s\, \xi^s\le f_N'\, f_N,$$ because   
\begin{equation*}
\begin{split}
f'_N(\xi)\, f_N(\xi)&=\biggl(\sum_{q=0}^{N-1}(q+1)\, \alpha_{q+1}\, \xi^q\biggr) \biggl(\sum_{p=0}^{N}\alpha_p\, \xi^p\biggr)\!=\!
\sum_{q=0}^{N-1}\, \sum_{p=0}^{N}(q+1)\, \alpha_{q+1}\, \alpha_p\, \xi^{q+p}\\
&=\sum_{s=0}^{2\, N-1} \biggl(\sum_{q=0}^{N-1}(q+1)\, \alpha_{q+1}\, \alpha_{s-q}\biggr) \xi^s\!\geq\!\sum_{s=0}^{ N} \biggl(\sum_{k=0}^{s-1}(k+1)\, \alpha_{k+1}\, \alpha_{s-k}\biggr) \xi^s\\
&=\sum_{s=0}^{ N}\, A_s\, \xi^s
\end{split}
\end{equation*}
and similar developments imply that also 
$$\sum_{s=0}^NB_s\, \xi^s\le g_N'\, g_N$$ and
$$\sum_{s=0}^NC_s\, \xi^s\le f_N'\, g_N+f_N\, g_N',$$
so 
\begin{equation*}
f'_{N+1}(\xi)\le\frac{1}{2}\, \sum_{s=0}^N\, \frac{\alpha^{s+1}}{s!}\, \xi^s+(1+K^2)\biggl(\frac{1}{2}\, f_N^2\biggr)'+K^2 \biggl(\biggl(\frac{1}{2}\, g_N^2\biggr)'+(f_N\, g_N)'\biggr)
\end{equation*} 
and 
$$
g'_{N+1}(\xi)\le(1+K^2)\biggl(\frac{1}{2}\, g_N^2\biggr)'+K^2 \biggl(\biggl(\frac{1}{2}\, f_N^2\biggr)'+(f_N\, g_N)'\biggr),
$$ 
and integrating
$$
f_{N+1}(\xi)\le  c_N+\frac{1}{2}\, \sum_{s=0}^N\, \frac{\alpha^{s+1}}{(s+1)!}\, \xi^{s+1}+\frac{1+K^2}{2}\, f_N^2+K^2\biggl(\frac{1}{2}\, g_N^2+f_N\, g_N\biggr)
$$ 
and 
$$
g_{N+1}(\xi)\le d_N+\frac{1+K^2}{2}\, g_N^2+K^2\biggl(\frac{1}{2}\, f_N^2+f_N\, g_N\biggr).
$$ 

Adding these two inequalities and considering that $f_N(0)=g_N(0)=0$, we have that 
$$
(f_{N+1}+g_{N+1)}(\xi)\le  \frac{\alpha\, \xi}{2}\, e^{\alpha\, \xi}+\frac{1+2\, K^2}{2}\, (f_N+g_N)^2.
$$
Then, if $h_N=f_N+g_N$, we have 
\begin{equation} \label{induc}
h_{N+1}\le  \frac{\alpha\, \xi}{2}\, e^{\alpha\, \xi}+\frac{1+2\, K^2}{2}\, h_N^2,
\end{equation}
and then 

\begin{prop}\label{|||}
For any $N\in\N\setminus\{0\}$, let $h_N$ a sequence of functions satisfying~\eqref{induc}  we have 
$$
h_N(\xi)\le\frac{2}{1+2\, K^2}
$$ for $\xi\in(0,\, \frac{1}{\alpha\, (1+2\, K^2)}\, e^{-\frac{4}{1+2\, K^2}})$.
\end{prop}          

\begin{proof}
We want to prove that there exists $L>0$ such that if 
$
h_1(\xi)= \frac{\alpha}{2}\, \xi\le L
$ 
and for any $N\in\N\setminus\{0\}$, we have $h_N\le L$, then  $h_{N+1}\le L$. 

The inequality (\ref{induc}), shows that it is enough to find $L$ satisfying 
$$
\frac{\alpha\, \xi}{2}\, e^{\alpha\, \xi}+\frac{1+2\, K^2}{2}\, L^2-L\le0,
$$ for $
\xi\le \frac{2\, L}{\alpha}
$, so let us fix $\xi$ for a moment.

%\bigskip

Using the notation $c=\frac{\alpha\, \xi}{2}\, e^{\alpha\, \xi}$ and $R=\frac{1+2\, K^2}{2}$, the polynomial 
$$
Rx^2-x+c=R\biggl(x-\frac{1-\sqrt{1-4\, R\, c}}{2R}\biggr) \biggl(x-\frac{1+\sqrt{1-4\, R\, c}}{2R}\biggr)
$$ 
has two real positive different roots whenever $0<c\le\frac{1}{4R}$, and, independently of $c\in(0,\frac{1}{4R})$, the value $\frac{1}{R}$ is the mean value of the roots. 

%\bigskip

In this case $L$ must be in the interval determined by the roots, so if $c=\frac{\mu}{4\, R}$ and $0<\mu<1$, the smallest root is $\frac{1-\sqrt{1-\mu}}{2\, R}$ an it is larger than $c$, so any number in the interval 
$$
\biggl(\frac{1-\sqrt{1-4\, R\, c}}{2R},\,  \frac{1+\sqrt{1-4\, R\, c}}{2R}\biggr)
$$ 
can be chosen as $L$, and in particular $h_N(\xi)\le \frac{1}{R}$ for any $N$ and any $c\in(0,\frac{1}{4R})$.        	
\end{proof}
%\bigskip
%\bigskip

Back to the decomposition \eqref{decomp} and the fact that $\alpha=2\, \|\varpi_0 \|$, the previous proposition implies that 
$$
\sum_{s=0}^N\biggl\|\frac{\partial\xi^{(s)}}{\partial z}\biggr\|\, \tau^s\le h_N(\tau)\le \frac{2}{1+2\, K^2}
$$ 
for every $N$ and $\tau\in[0,\, \frac{1}{\alpha\, (1+2\, K^2)}\, e^{-\frac{4}{1+2\, K^2}})$. 

%\medskip

Since the functions $\frac{\partial\xi^{(s)}}{\partial z}$ are in $\ka{1,\, \gamma}$ and, as the recurrence shows, they exhibit a decay at $\infty$ of order $\frac{1}{|z|^2}$, then  
$$
\xi^{(s)}=\bar C\biggl[\frac{\partial\xi^{(s)}}{\partial z}\biggr]
$$ 
is in $\ka{1,\, \gamma}$ (cf.\ Theorems 4.3.11 and 4.3.12 in \cite{AsIwMa}). Moreover 
$$
\|\xi^{(s)}\|_{L^\infty}\le K_0 \biggl\|\frac{\partial\xi^{(s)}}{\partial z}\biggr\|_{\gamma}.
$$

%\medskip

Consequently, for any $N$, 
$$
\biggl|\sum_{s=0}^N \xi^{(s)}(z)\, t^s\biggr|\le\sum_{s=0}^N \|\xi^{(s)}\|_{L^\infty}\, |t|^s\le K_0\, \sum_{s=0}^N 
\biggl\|\frac{\partial\xi^{(s)}}{\partial z}\biggr\|\, |t|^s
$$ 
and the last term is uniformly bounded in $N$, for $|t|\le \frac{1}{\|\varpi_0 \|\, (1+2\, K^2)}\, e^{-\frac{4}{1+2\, K^2}}$.

This implies the existence of an analytic solution of the equation~\eqref{Meq}, providing the proof of Theorem \ref{Mth1}.

%\bigskip
\bigskip

\subsection{Proof of Proposition \ref{ProP}}

We have seen that $$\phi(z,\, t)=z+\sum_{s=1}^\infty \xi^{(s)}(z)\, t^s\stackrel{\text{def}}{=}z+A(z,\, t)$$ and $$\xi^{(s)}(z)=\bar C[\theta^{(s)}](z).$$ 

Let us prove first an auxiliary lemma.

\begin{lem} There exists  $T_0>0$ and $K>0$ such that for $|t|\le T_0$, we have $$\|A(\ ,\, t)\|_{L^\infty(\C)}\le K$$ and $$\|\nabla_zA(\ ,\, t)\|_{{L^\infty(\C)}\le\frac{1}{2}}.$$

\end{lem}
\begin{proof} Since $$\sum_{s=1}^N\|\xi^{(s)}\|_{L^\infty(\C)}\, \tau^s=\sum_{s=1}^N\|\bar C[\theta^{(s)}]\, \|_{L^\infty(\C)}\, \tau^s$$ and $$\|\bar C[\theta^{(s)}]\, \|_{L^\infty(\C)}\, \le C\{\|\theta^{(s)}\, \|_{L^\infty(\C)}+\|\theta^{(s)}\, \|_{L^p(\C)}\}$$ with $C$ independent of $s$ and $N$, and $1<p<2$, due to the Proposition \ref{|||}, we obtain the first inequality.

\bigskip

On the other hand 
\begin{equation*}
\begin{split}
\sum_{s=1}^N\|\nabla_z\xi^{(s)}\|_{L^\infty(\C)}\, \tau^s&\le\sum_{s=1}^N\{\|\theta^{(s)}\|_{L^\infty(\C)}+\|\bar B[\theta^{(s)}]\|_{L^\infty(\C)}\}\, \tau^s\\
&\le K\, \sum_{s=1}^N\|\theta^{(s)}\|_{\gamma}\, \tau^s
\end{split}
\end{equation*}
 where $K$ is independent of $s$ and $N$. We have used the part of Theorem 3.
\end{proof}
%\bigskip
\bigskip

\subsubsection{$\phi$ is a flow.}

\begin{prop} \begin{enumerate}
	\item $\phi(\C,\, t)$ is closed for every $t$.
	
	\bigskip
	
	\item $\phi(\ ,\, t)$ is one to one. 
	
	\bigskip
	
	\item $\phi(\C ,\, t)=\C$. 
	
	\bigskip
	
	\item $\phi(\ ,\, t)$ is an homeomorphism and differentiable in $\C\setminus\partial\Omega$. 
\end{enumerate}

\end{prop}

\begin{proof} \begin{enumerate}
	\item If $w_l=\phi(z_l,\, t)=z_l+A(z_l,\, t)$ is a sequence such that $w_l\rightarrow w_0$, then $|z_l|\le |w_l|+|A(z_l,\, t)|\le M+\|A(\ ,\, t)\|_{L^\infty},$ where $M$ is a bound for~$w_l$. Then there exists a partial sequence $z_{l_s}\rightarrow z_0$
	and since $\phi(\ ,\, t)$ is a continuous map then $w_0=\phi(z_0,\, t)$.

	\bigskip
	
	\item Since $$\phi(z,\, t)-\phi(z',\, t)=z-z'+A(z,\, t)-A(z',\, t),$$ and $\nabla A(\ ,\, t)$ exists almost everywhere in $\C$ and $\|\nabla A(\ ,\, t)\|_{L^\infty}<1$ if $|t|\le T'''_0$, then by Rademacher's theorem $A(\ ,\, t)\in \operatorname{Lip}(1,\, \C)$, so if %\newline 
	$\phi(z,\, t)-\phi(z',\, t)=0$, then $|z-z'|<|z-z'|$. That's a contradiction. 
	
	\bigskip
	
	\item Since for some $\tilde T\in(0,\, T'''_0]$ we have $|J_zA(z,\, t)|<\frac{1}{2}$ for $(z,\, t)\in(\C\setminus\partial\Omega)\times[-\tilde T,\, \tilde T]$, and $\phi(\ ,\, t)\in\ka{\infty}(\C\setminus\partial\Omega)$, then  $\phi(\C\setminus\partial\Omega ,\, t)$ is open and $\phi(\ ,\, t)$ is a diffeomorphism onto the image, because of the e inverse function theorem and the injectivity.
	
	\medskip 
	
	Let us start proving that the injectivity of $\phi(\ ,\, t)$ also implies that $$\phi(\partial\Omega,\, t)=\partial\phi(\Omega,\, t)\cap\partial\phi(\bar\Omega^c,\, t).$$
	
	Since for $w_0=\phi(z_0,\, t)$, with $z_0\in\partial\Omega$, we have that $z_0=\lim z_l=\lim z'_l$ where $z_l\in\Omega$ and $z'_l\in\bar\Omega^c$, then $w_0=\lim\phi(z_l,\, t)=\lim\phi(z'_l,\, t)$, as $\partial\Omega$ is regular. This implies the inclusion of the left hand side.
	
	On the other hand, for $w_0\in\partial\phi(\Omega,\, t)\cap\partial\phi(\Omega^c,\, t)$ we have that $w_0=\lim\phi(z_l,\, t)=\lim\phi(z'_l,\, t)$, where $z_l\in\Omega$ and $z'_l\in\bar\Omega^c$, so $z_l=w_0-A(z_l,\, t)$ and $z'_l=w_0-A(z'_l,\, t)$ and there are partial subsequences $z_{l_j}\to\alpha\in\bar\Omega$ and $z'_{l_k}\to\beta\in\Omega^c$. Then $\phi(\alpha,\, t)=\phi(\beta,\, t)=w_0$. The injectivity of $\phi(\ ,\, t)$ implies $\alpha=\beta\in\partial\Omega$. 
	
	\medskip
	
	Let us see that, in fact, $\phi(\ ,\, t)$ maps boundaries to boundaries. 
	
	We have
	$$\phi(\partial\Omega ,\, t)=\partial\phi(\Omega,\, t)=\partial\phi(\bar\Omega^c ,\, t),$$ because $\phi(\partial\Omega ,\, t)\subset\partial\phi(\Omega,\, t)$ by continuity, and if $w_0\in\partial\phi(\Omega,\, t)\subset\overline{\phi(\C,\, t)}=
	\phi(\C,\, t)$, then $w_0=\phi(z_0,\, t)$ for a unique $z_0\in\C$. 
	
	We have that $z_0\notin\Omega$ because otherwise $w_0\in\phi(\Omega,\, t),$ an open set, and similarly $z_0\notin\bar\Omega^c$. So $z_0\in\partial\Omega$.
	
	The same argument can be performed for $\partial\phi(\bar\Omega^c,\, t)$.
	
	\medskip
	
    Then, as $\Omega$ is a bounded subset of $\C$, we finally conclude that $$\partial\phi(\bar\Omega,\, t)=\phi(\partial\Omega,\, t),$$ because if $w\in\partial\phi(\bar\Omega,\, t),$ then $w\in\phi(\bar\Omega,\, t),$ so $w=\phi(z,\, t),$ for $z\in\bar\Omega$, and in fact $z\in\partial\Omega$ because $\phi(\Omega,\, t)$ is open. Assume now that $w\in\phi(\partial\Omega,\, t)\subset\phi(\bar\Omega,\, t)$ then $w=\phi(z,\, t),$ for $z\in\bar\Omega$, and if $w$ is an interior point of $\phi(\bar\Omega,\, t ),$ then, by continuity of $\phi(\ ,\, t)$, we have $z$ is an interior point of $\bar\Omega$, and since $\partial\Omega$ is a regular submanifold of $\C$, then $z\notin\partial\Omega$. That's a contradiction. 
	
	\medskip
	
	As a consequence of all these facts, $$\phi(\C,\, t)=\C.$$
	
	If $w_0\in\C\setminus\phi(\C,\, t)$, let $w_1=\phi(z_1,\, t)$ a point satisfiying $|w_0-w_1|=d(w_0,\, \phi(\C,\, t))$. Clearly $z_1\in\partial\Omega$ because $\phi(\Omega,\, t)$ and $\phi(\bar\Omega^c,\, t)$ are open sets. In fact, $$w_1\in\partial\phi(\Omega,\, t)\cap\partial\phi(\bar\Omega^c,\, t)\cap\partial\phi(\C,\, t).$$
	
	On the other hand, $\phi(\partial\Omega,\, t)$ is a Jordan curve in $\C$, and then it divides $\C$ in two disjoint regions, say $\Pi_+$ and $\Pi_-$ and the curve itself: $$\C=\Pi_+\cup\phi(\partial\Omega,\, t)\cup\Pi_-,$$ and $\Pi_-$ is bounded.
	
	Let us prove now that $\phi(\Omega,\, t)\subset\Pi_-$. Since $\phi(\bar\Omega,\, t)$ is a compact set and then there is a large ball $B_R(0)$ such that $$\phi(\bar\Omega,\, t)\cup\Pi_-\subset\overline{B_R(0)}.$$ For a point $u_0\notin B_{2R}(0)$, then $u_0\in\Pi_+$ and if we assume that $\phi(\Omega,\, t)\subset\Pi_+$, then since $\Pi_+$ is arc-connected there must exist a continuous curve $$p\colon [0,1]\rightarrow\Pi_+$$ and a point $z_2\in\Omega$ such that $p(0)=u_0$, $p(1)=\phi(z_2,\, t)$. By connection arguments this implies that $p([0,\, 1])\cap\partial\phi(\bar\Omega,\, t)\neq\emptyset$, so %\newline 
	$\partial\phi(\bar\Omega,\, t)\setminus\phi(\partial\Omega,\, t)\neq\emptyset$. That's a contradiction. 
	
	\medskip
	
	Using the same arguments, we also prove that $\phi(\bar\Omega^c,\, t)\subset\Pi_+$. Let us show now that in fact $\phi(\Omega,\, t)=\Pi_-$ and $\phi(\bar\Omega^c,\, t)=\Pi_+$. If $w_0\in\Pi_+\setminus\phi(\bar\Omega^c,\, t)$ There exists an arc in $\Pi_+$, $p$,  joining $w_0$ with a point in $\phi(\bar\Omega^c,\, t)$. Again this arc must cross the curve $\partial\phi(\bar\Omega^c,\, t)=\partial\phi(\Omega,\, t)$. that's a contradiction.  
	
	\bigskip
	
	\item Follows directly from the previous considerations. \qedhere
\end{enumerate}
\end{proof}
%\bigskip
%\bigskip

% !TeX encoding = UTF-8
% !TeX spellcheck = ca_ES
% !TeX root = Analit4.tex 
\section{Proof of Theorem \ref{Mth3}}
	
This is the most technical section of the paper. We divide it in several subsections, where we exhibit some structural and geometric facts about the Beurling transform, necessary for our purposes, and then we perform the uniform and Lipschitz estimates, necessary for the proof of Theorem \ref{Mth1}.
	
%	\bigskip
%	\bigskip
	
Let us consider a function $f$ defined on a domain $W\subset\C$ with bounded $\ka{1,\gamma}$-regular boundary. We will also denote by $f$ the extension of $f$ by $0$ to $\bar W^c$, and $f(z)=\lim_{w\in W\to z}f(w)$ if $z\in\partial W$. 
	
Whenever it be  necessary, we will specialize $f=\phi$ and then $W=\spt(\phi)$, the support of $\phi$, or $f=\psi$ and $W=\spt(\psi)$, respectively.
	
%	\bigskip
%	\bigskip
	
\subsection{The geometric lemma} 
The metric and geometric properties of $\partial\Omega$ play a crucial role in the behavior of the Beurling transform. The next lemma is a synthesis of these properties.
	
\begin{lem}\label{geom} 
Let $\Omega\subset\C$ be a bounded domain such that $\partial\Omega\in\ka{1,\gamma}$, defined by a function $\rho$. 
		
There exists $0<R_0<1$ such that if 
$$
U_{R_0}(\partial\Omega):=\cup_{z\in\partial\Omega}B_{R_0}(z),
$$ 
then there exists $R_1$ such that for any $z_0\in U_{R_0}$ the level set $\{\rho=\rho(z_0)\}\cap B_{R_1}(z_0)$ coincides with the graph of a function $\varphi_{z_0}$.

Moreover $\varphi_{z_0}$ is $\ka{1,\gamma}$ and $\varphi_{z_0}(0)=\varphi'_{z_0}(0)=0$.
\end{lem}
\medskip

\begin{rmk} The function $\varphi_{z_0}$ is defined on a segment of the tangent line to the level set of $\rho$ across $z_0$.  
	
	Then
$$	\rho(z_0+s\,\eta(z_0)^\perp+\varphi_{z_0}(s)\,\eta(z_0))=\rho(z_0).
	$$ 
	
	\end{rmk}

%\medskip
	
\begin{proof} 
Let $\rho\in\ka{1,\gamma}$ be the defining  function for $\Omega$, as in \cite{Bur}. Since $\partial\Omega$ is a compact set, let 
$$
m=\min\{\|\nabla\rho(z)\|: z\in\partial\Omega\}.
$$
		
Since $\nabla\rho\in \operatorname{Lip}(\gamma)$, then for  
$$
R=\biggl(\frac{ m}{2\sqrt{2}\, \|\rho\|_\gamma}\biggr)^\frac{1}{\gamma},
$$ 
we define $U_R=\cup_{\zeta\in\partial\Omega}B_R(\zeta)$. It is an open neighborhood of $\partial\Omega$, and for every~$z\in U_R$, 
$$
\|\nabla\rho(z)\|\geq\frac{m}{2}.
$$
		
%\bigskip
		
Fix $z_0\in U_{\frac{R}{2}}$. The vectors $\eta(z_0)=\frac{\nabla\rho(z_0)}{\|\nabla\rho(z_0)\|}$ and $\tau(z_0)=\eta(z_0)^\perp$ form an orthonormal basis of $T_{z_0}(\R^2)$, and the map 
$$
\Xi\colon \R^2\rightarrow\C
$$ 
given by 
$$
\Xi(s,t)=z_0+s\, \tau(z_0)+t\, \eta(z_0)
$$ 
is a rigid movement.
		
The function 
$$
r(s,t)=\rho\circ\Xi(s,t)
$$ 
satisfies 
$$
\nabla r(0,0)=\|\nabla\rho(z_0)\|\, e_2.
$$
		
Let us consider the map 
$$
\Psi(s,t)=\biggl(s,\, \frac{r(s,t)-r(0,0)}{\|\nabla\rho(z_0)\|}\biggr).
$$ 
We have that $\Psi(0,\, 0)=(0,0)$ and 
$$
\det J\Psi(0,0)=\frac{\frac{\partial\rho}{\partial y}(z_0)}{\|\nabla\rho(z_0)\|}.
$$
		
%\bigskip
		
Next, $\Psi(s,t)=\Psi(s',t')$ iff $s=s'$ and $r(s,t)=r(s',t')$, so, for some $\xi$ in the open interval delimited by $t$ and $t'$, we have  
\begin{equation*}
\begin{split}
0&=r(s,t)-r(s,t')=\frac{\partial r}{\partial t}(s,\xi)\, (t-t')\\*[5pt]
&=\frac{\partial r}{\partial t}(0,0)\, (t-t')+\biggl(\frac{\partial r}{\partial t}(s,\xi)-\frac{\partial r}{\partial t}(0,0)\biggr)\, (t-t'),
\end{split}
\end{equation*} 
so 
\begin{equation*}
\begin{split}
0&\geq\biggl|\frac{\partial r}{\partial t}(0,0)\biggr|\, |t-t'|-\biggl\|\frac{\partial r}{\partial t}\biggr\|_{\Lip{(\gamma)}}\, |t-t'|^{1+\gamma}\\*[5pt]
&=\biggl\{\|\nabla\rho(z_0)\|-\biggl\|\frac{\partial r}{\partial t}\biggr\|_{\Lip{(\gamma)}}\, |t-t'|^{\gamma}\biggr\}\, |t-t'|,
\end{split}
\end{equation*} 
and if we choose $R_1$ in such a way that if $z_0\in\overline U_{R_1}$, then $B_{R_1}(z)\subset U_{\frac{R}{2}}$ and $R_1\le(\frac{\|\nabla\rho(z_0)\|}{4\, \|\frac{\partial r}{\partial t}\|_{\Lip{(\gamma)}}})^\frac{1}{\gamma}$, then 
$$
0\geq\frac{1}{2}\, \|\nabla\rho(z_0)\|\, |t-t'|.
$$  
		
So $\Psi$ is one-to-one in $B_{R_1}(0,0)$.
		
%\bigskip
		
Choose, now, $0<R_2<R_1$. Then, if $s^2+t^2=R_2^2$, we have 
\begin{equation*}
\begin{split}
\|\Psi(s,t)\|^2&=s^2+ \frac{(r(s,t)-r(0,0))^2}{\|\nabla\rho(z_0)\|^2}\\*[5pt]
&=s^2+ \biggl(t+\frac{(\nabla r(s',t')-\nabla r(0,0),\, (s,t))}{\|\nabla\rho(z_0)\|}\biggr)^2\\*[5pt]
&\geq R_2^2\biggl\{1+\biggl(\frac{(\nabla r(s',t')-\nabla r(0,0),\, u)}{\|\nabla\rho(z_0)\|}\biggr)^2\\*[5pt]
&\hspace*{1.5cm}-2\, |(u,\, e_2)|
\biggl|\frac{(\nabla r(s',t')-\nabla r(0,0),\, u)}{\|\nabla\rho(z_0)\|}\biggr|\biggr\}=(*),
\end{split}
\end{equation*} 
where $u$ is a unitary vector. Moreover 
$$
(*)\geq R_2^2\biggl\{1-2\, \frac{\|\nabla r\|_{\Lip{(\gamma)}}\, R_2^\gamma}{\|\nabla\rho(z_0)\|}|\biggr\},
$$ 
and if 
$$
R_2\le\biggl(\frac{\|\nabla\rho(z_0)\|}{4\, \|\nabla r\|_{\Lip{(\gamma)}}}\biggr)^{\frac{1}{\gamma}},
$$ 
then 
$$
d(\Psi(\partial B_{R_2}(0,0)),\, (0,0))\geq\frac{R_2}{\sqrt{2}}.
 $$ 
This implies that for $(p,q)\in B_{\frac{R_2}{2\, \sqrt{2}}}(0,0)$, the function $g(s,t)=\|\Psi(s,t)-(p,q)\|^2$ has a local minimum at some point $(s_0,\, t_0)\in B_{R_2}(0,0)$, and then 
$$
\begin{cases} 
2(s_0-p)+2(\frac{r(s_0,t_0)-r(0,0)}{\|\nabla\rho(z_0)\|}-q)\, \frac{\frac{\partial r}{\partial s}(s_0,t_0)}{\|\nabla\rho(z_0)\|}=0,\\*[5pt]
2(\frac{r(s_0,t_0)-r(0,0)}{\|\nabla\rho(z_0)\|}-q)\, \frac{\frac{\partial r}{\partial t}(s_0,t_0)}{\|\nabla\rho(z_0)\|}=0.
\end{cases}
$$
	
And we can change $R_2$ by $R_3$ so that $\frac{\partial r}{\partial t}$ do not vanish at $B_{R_3}(0,0)$ and the other properties still hold.
		
%\bigskip
		
So, there is a new choice of the region $U$ such that $\|\nabla\rho\|$ is uniformly bounded from below on $U$, and then there are choices of $R',R''>0$ such that for every~$z_0\in U$ the corresponding map $\Psi$ composed with the rigid movement describe above map $\ka{1,\gamma}$-diffeomorphically $B_{R'}(z_0)$ in $B_{R''}(z_0)$.
		
%\bigskip
		
If $\Phi=\Psi^{-1}$, then $\Phi$ is a $\ka{1,\gamma}$-diffeomorphism such that $\Phi(0,0)=(0,0)$. Moreover $\Psi_1(s,t)=s$ and 
$$
\Psi_2(\{(s,t): r(s,t)=r(0,0)\})=0,
$$ 
so the level set  $\{r(s,t)=r(0,0)\}\cap B_{R'}(0,0)$ coincides wit the set $\{t=\Phi_2(s,0)\}$. 
		
Taking a convenient rectangle inside the balls, we have the function $\varphi$ in the statement of the theorem.
		
Finally, since essentially $\varphi(s)=\Phi_2(s,0)$, we have that 
$$
\varphi'(0)=\frac{\partial\Phi_2}{\partial p}(0,0),
$$ 
and since 
$$
J\Phi(0,0)=(J\Psi(0,0))^{-1}
$$ 
we have $\varphi'(0)=0$.
		
%\bigskip
		
Using the rigid movement above, we can take the construction backward to~$\C$ and get the theorem. 
\end{proof}

%	\bigskip
%	
%	\bigskip   

%	\bigskip 
%	\bigskip

\subsection{Decomposition of the singularities}

Set, for any $z\in\C$, $d(z)=d(z,\, \partial\Omega)$ and $\delta(z)=\max\{d(z),\, \frac{R_0}{2}\}$ ($R_0$ as in Lemma \ref{geom}).

\begin{prop}\label{desing} Let $\Omega\subset\C$ be a bounded domain with boundary of class $\ka{1, \gamma}$ and let $W$ denote $\Omega$ or $\C\setminus\bar\Omega$.
	 
If $f\in\ka{\gamma}(\bar W)\cap L^p(W)$ with $p>1$ and we identify $f$ with its extension by $0$ outside $W$, then there exists 
$$
\bar B[f](z)=\lim_{\epsilon\to0}\bar B_\epsilon[f](z)
$$ 
for every $z\in\C$.
		
%\bigskip
		
Moreover
$$
\bar B[f](z)=Q[f](z)+L[f](z)+f(z)\, \Theta_W^{\frac{R_0}{2}}(z),
$$ 
where
\begin{align*}
Q[f](z)&=\int_{\C\setminus B_{\delta(z)}(z)}f(\zeta)\, \frac{1}{(\bar\zeta-\bar z)^2}\, dm(\zeta),\\*[5pt]
L[f](z)&=\int_{B_{\delta(z)}(z)}\frac{f(\zeta)-f(z)}{|\zeta-z|^\gamma}\, 
\frac{|\zeta-z|^\gamma}{(\bar\zeta-\bar z)^2}\, dm(\zeta),
\end{align*} 
and		
$$
\Theta_W^{\frac{R_0}{2}}(z)=\begin{cases} 
0 &\text{if } d(z)>0,\\
\{\int_{B_{\frac{R_0}{2}}(z)\cap\bar W\cap\{\kappa_z>0\}}-\int_{B_{\frac{R_0}{2}}(z)\cap W^c\cap\{\kappa_z<0\}}\}\frac{ dm(\zeta)}{(\bar\zeta-\bar z)^2}\,  &\text{if } d(z)=0. 
\end{cases}
$$ 

\end{prop}

\medskip
\begin{rmk} The term $\Theta_W^{\frac{R_0}{2}}(z)$ is an intrinsec geometric object. Sometimes we will also use it in the form 
	$$
	\int_{B_{\frac{R_0}{2}}(z)\cap W}\frac{ dm(\zeta)}{(\bar\zeta-\bar z)^2}.
	$$
	
	\end{rmk}

\medskip

\begin{rmk} In order to avoid notation, from now on we will use in the proofs of results the notation $b_z(\zeta)$ instead of $\frac{ dm(\zeta)}{(\bar\zeta-\bar z)^2}$.
\end{rmk}

\medskip

\begin{rmk} Since the Beurling transform is a classical Calder\'on--Zygmund operator, the existence of the principal value is well known for functions in many different classes (See \cite[Corollary 5.8]{Duo}), nevertheless, we need here the existence of the princial value at each point in the plane, for functions in the aforementioned class. 
\end{rmk}

\begin{proof} Let $z\in\C$. 
	
	%\bigskip
	
	\begin{itemize}
		
		\item  If $d(z)>0$, then for $0<\epsilon<d(z)$ we have 
		$$
		B_\epsilon[f](z)=\int_{\C\setminus B_{\epsilon}(z)}f\, b_z=\int_{\C\setminus B_{d(z)}(z)}f\, b_z
		+\int_{B_{d(z)}(z)\setminus B_\epsilon(z)}f\, b_z.
		$$
		
		%\bigskip
		
		\begin{itemize}
			\item In the case of $z\notin W$, the second term is $0$ and we have 
			$$
			B[f](z)=\lim_{\epsilon\to0} B_\epsilon[f](z)=\int_{\C\setminus B_{d(z)}(z)}f\, b_z.
			$$
			
			%\bigskip
			
			\item In the case of $z\in W$, we have   
			\begin{equation*}
				\begin{split}
					B_\epsilon[f](z)&=\!\int_{\C\setminus B_{d(z)}(z)}f\, b_z\!+\!\int_{B_{d(z)}(z)\setminus B_\epsilon(z)}(f-f(z))\, b_z\\*[5pt]
					&=\!\int_{\C\setminus B_{d(z)}(z)}f\, b_z\!+\!\int_{B_{d(z)}(z)\setminus B_\epsilon(z)}\!\frac{f(\zeta)-f(z)}{|\zeta-z|^\gamma}\, 
					\frac{|\zeta-z|^\gamma}{(\bar\zeta-\bar z)^2}\, dm(\zeta),\!\!\!
				\end{split}
			\end{equation*} 
			and in the second term  
			has a weakely singular kernel acts against an integrable function, so the limit exists   
			\begin{equation*}
				\begin{split}
					B[f](z)&=\lim_{\epsilon\to0} B_\epsilon[f](z)\\*[5pt]
					&=\int_{\C\setminus B_{d(z)}(z)}f\, b_z+\int_{B_{d(z)}(z)}\frac{f(\zeta)-f(z)}{|\zeta-z|^\gamma}\, 
					\frac{|\zeta-z|^\gamma}{(\bar\zeta-\bar z)^2}\, dm(\zeta).
				\end{split}
			\end{equation*} 
			
			%\bigskip
			
			If $d(z)\leq\frac{R_0}{2}$, then 
			$$
			\int_{\C\setminus B_{d(z)}(z)}f\, b_z=
			\int_{\C\setminus B_{\frac{R_0}{2}}(z)}f\, b_z+\int_{B_{\frac{R_0}{2}}(z)\setminus B_{d(z)}(z)}f\, b_z,
			$$ 
			and  applying again the cancellation lemma (below), we have that the second integral is  
			$$
			\int_{C_{d(z)}^{\frac{R_0}{2}}(z)}\frac{f(\zeta)-f(z)}{|\zeta-z|^\gamma}\, 
			\frac{|\zeta-z|^\gamma}{(\bar\zeta-\bar z)^2}\, dm(\zeta),
			$$ where $C_{d(z)}^{\frac{R_0}{2}}=B_{\frac{R_0}{2}}(z)\setminus B_{d(z)}(z)$.

			Implementing these facts in the formula for $B[f](z)$ above, we have the result. 
		\end{itemize}	
		
		%\bigskip
		
		\item If $d(z)=0$, then we for any $\epsilon<\frac{R_0}{2}$ we have 
		$$
		B_\epsilon[f](z)=\int_{\C\setminus B_{\frac{R_0}{2}}(z)}f\, b_z
		+\int_{C^{\frac{R_0}{2}}_\epsilon(z)}f\, b_z.
		$$
		
		The second term, is 
		\begin{equation*}
			\begin{split}
				\int_{C^{\frac{R_0}{2}}_\epsilon(z)\cap\bar W} f\, b_z&=\int_{C^{\frac{R_0}{2}}_\epsilon(z)\cap\bar W}(f-f(z))\, b_z+f(z)\, \int_{C^{\frac{R_0}{2}}_\epsilon(z)\cap\bar W}b_z\\*[5pt]
				&=\int_{C^{\frac{R_0}{2}}_\epsilon(z)\cap\bar W}\frac{f(\zeta)-f(z)}{|\zeta-z|^\gamma}\, 
				\frac{|\zeta-z|^\gamma}{(\bar\zeta-\bar z)^2}\, dm(\zeta)\\*[5pt]
				&\quad+f(z) \int_{C^{\frac{R_0}{2}}_\epsilon(z)\cap\bar W}b_z.
			\end{split}
		\end{equation*}
		
		From these two terms, the first one is weakly singular integral and has a limit
		$$
		\int_{B_{\frac{R_0}{2}}(z)\cap\bar W}\frac{f(\zeta)-f(z)}{|\zeta-z|^\gamma}\, 
		\frac{|\zeta-z|^\gamma}{(\bar\zeta-\bar z)^2}\, dm(\zeta).
		$$
		
		%\medskip
		
		The second term, according to the Beurling geometric lemma (below) has a limit as $\epsilon\to0$, that is
\begin{equation*}
		f(z)\biggl\{\int_{B_{\frac{R_0}{2}}(z)\cap\bar W\cap\{\kappa_z>0\}}-\int_{B_{\frac{R_0}{2}}(z)\cap W^c\cap\{\kappa_z<0\}}\biggr\}b_z.\qedhere
\end{equation*}
	\end{itemize}	
\end{proof}
%\bigskip
%\bigskip
%
%
%
%\bigskip
%\bigskip

The next fact is proved in Lemma 3 in \cite{MaOrVe} in a more general situation. For completeness we give a proof adapted to this case.

\begin{lem}[Beurling-Geometric lemma] 
If $\Omega\subset\C$ and $\partial \Omega$ is compact and $\ka{1,\gamma}$, then there exists 
$$
\lim_{\epsilon\to0}\int_{C^{\frac{R_0}{2}}_\epsilon(z)\cap\bar \Omega}\frac{ dm(\zeta)}{(\bar\zeta-\bar z)^2}=\biggl\{\int_{B_{\frac{R_0}{2}}(z)\cap\bar \Omega\cap\{\kappa_z>0\}}-\int_{B_{\frac{R_0}{2}}(z)\cap \Omega^c\cap\{\kappa_z<0\}}\biggr\}\frac{ dm(\zeta)}{(\bar\zeta-\bar z)^2}.
$$  
\end{lem}

\begin{proof} 
First of all, we have 
$$
\!\int_{C^{\frac{R_0}{2}}_\epsilon(z)\cap\bar \Omega} \!b_z\!=\!\!
\int_{C^{\frac{R_0}{2}}_\epsilon(z)\cap\{\kappa_z\le0\}\!}b_z+\biggl\{\int_{B_{\frac{R_0}{2}}(z)\cap\bar \Omega\cap\{\kappa_z>0\}}\!-\!\int_{B_{\frac{R_0}{2}}(z)\cap \Omega^c\cap\{\kappa_z<0\}}\biggr\}b_z,\!\!\!
$$ 
and 
$$
\int_{C^{\frac{R_0}{2}}_\epsilon(z)\cap\{\kappa_z\le0\}}b_z=\int_\epsilon^{\frac{R_0}{2}}\frac{dr}{r}\, \int_0^\pi e^{2i\theta}\, d\theta=0.
$$ 

From the geometric lemma, we have that 
\begin{multline*}
\int_{C^{\frac{R_0}{2}}_\epsilon(z)\cap\{\kappa_z\le0\}}b_z+
\biggl\{\int_{B_{\frac{R_0}{2}}(z)\cap\bar \Omega\cap\{\kappa_z>0\}}-\int_{B_{\frac{R_0}{2}}(z)\cap \Omega^c\cap\{\kappa_z<0\}}\biggr\}b_z\\*[5pt]
=\int_\epsilon^{\frac{R_0}{2}}\frac{dr}{r}\, \int_{I(r)} e^{2i\theta}\, d\theta=\int_\epsilon^{\frac{R_0}{2}}\vartheta_\Omega(z;\, r)\, \frac{dr}{r^{1-\gamma}},
\end{multline*} 
where 
$$
\vartheta_\Omega(z;\, r)=\frac{\cos\alpha(r)-1+i\sin\alpha(r)-1-\cos\beta(r)-i\sin\beta(r)}{r^\gamma}.
$$

%\medskip

By the geometric Lemma \ref{geom}, we have 
$$
|\vartheta(r)|\le 2M
$$ 
and then the limit exists and is equal to 
\begin{equation*}
\int_0^{\frac{R_0}{2}}\vartheta_\Omega(z;\, r)\, \frac{dr}{r^{1-\gamma}}.\qedhere
\end{equation*}
\end{proof}
%\bigskip
%\bigskip

\begin{lem}[Cancellation lemma]\label{can}  
If $z\in B_{R_0}(z)$ and $\epsilon<R_0-|z-w|$, then 
$$
\int_{B_{R_0}(w)\setminus B_\epsilon(z)}\frac{ dm(\zeta)}{(\bar\zeta-\bar z)^2}=0.
$$
\end{lem}
\begin{proof} 
In complex coordinates, the integral is 
$$
\frac{1}{2i}\, \int_{B_{R_0}(w)\setminus B_\epsilon(z)}
\frac{d\bar\zeta\wedge d\zeta}{(\bar\zeta-\bar z)^2}=\frac{i}{2}\biggl(\int_{\partial B_{R_0}(w)}-\int_{\partial 
	B_\epsilon(z)}\biggr)
\frac{d\zeta}{\bar\zeta-\bar z},
$$ 
and 
$$
\int_{\partial B_\epsilon(z)}
\frac{d\zeta}{\bar\zeta-\bar z}=\int_0^{2\pi}i\, e^{2i\theta}\, d\theta=0,
$$ 
and  
\begin{equation*}
\begin{split}
\int_{\partial B_{\R_0}(w)}\frac{d\zeta}{\bar\zeta-\bar z}
&=\int_0^{2\pi}\frac{iR_0\, e^{i\theta}\, d\theta}{R_0\, e^{-i\theta}+\overline{w-z}}\\*[5pt]
&=R_0\, \int_0^{2\pi}\frac{e^{i\theta}\, i e^{i\theta}\, d\theta}{R_0+\overline{(w-z)}\, e^{i\theta}}
=R_0\, \int_{\partial B_1(0)}\frac{\tau\, d\tau}{R_0+\overline{(w-z)}\, \tau}=0,
\end{split}
\end{equation*} 
because 
$\frac{\tau}{R_0+\overline{(w-z)}\, \tau}$ is a function  holomorphic in $\tau$, as $|\frac{R_0}{\overline{(w-z)}}|>1$.
\end{proof}

%\bigskip\bigskip

%\bigskip
%\bigskip

\subsection{The jump formula for $\bar B$}

Now we prove the jump formula 
$$
\bar B[g](z)=\frac{1}{2}\biggl\{\lim_{w\to z;\, w\in\Omega} \chi_\Omega\,  \bar B[g](w)+\lim_{w\to z;\, w\in\C\setminus\Omega} \chi_{\C\setminus \bar\Omega}\, \bar B[g](w)\biggr\}
$$ 
that is indentity \eqref{jump} in Theorem \ref{Mth3}. 

%\bigskip

\begin{rmk} Jump formulas of this type, for Calder\'on--Zygmund operators in potential theory appear in \cite{HoMiTa}. For the special case of the conjugate Beurling transform we give a simple proof to keep the exposition more self-contained.
\end{rmk}

%	\bigskip
%	\bigskip
	
As we will see in section 4.4 and 4.5 we have that $\chi_\Omega\, \bar B[g]\in \operatorname{Lip}_{\gamma}(\bar\Omega)$ and  $\chi_{\C\setminus \bar\Omega}\, \bar  B[g]\in  \operatorname{Lip}_{\gamma}(\C\setminus\Omega)$, and then the limits $\lim_{w\to z;\, w\in\Omega} \bar \chi_\Omega\, B[g](w)$ and $\lim_{w\to z;\, w\in\C\setminus\Omega} \chi_{\C\setminus \bar\Omega}\, \bar B[g](w)$ both exist and we can choose $w=z\pm\lambda\eta(z)$, where $\eta(z)$ is the unit vector, normal exterior to $\partial\Omega$ at $z$.
	
%\medskip
	
Also, if $g_\pm$ are the Lipschitz extensions of $g$ to $\Omega^c$ and $\bar\Omega$, respectively, we have, for $w\in B_{\frac{R_0}{8}}(z)$, the following facts	
\begin{itemize}
\item If $w\in\Omega$, then 
\begin{equation*}
\begin{split}
\bar B[g_-](w)&=\int_\Omega (g_--g_-(w))\, b_w+g_-(w)\, \int_\Omega b_w\\*[5pt]
&=\int_{\Omega\setminus B_{2\, |z-w|}(z)} (g_--g_-(w))\, b_w\\*[5pt]
&\quad+\int_{\Omega\cap B_{2\, |z-w|}(z)} (g_--g_-(w))\, b_w+g_-(w)\, \int_\Omega b_w\\*[5pt]
&=(I)(w)+(II)(w)+g_-(w)\, (III)(w).
\end{split}
\end{equation*}

%\bigskip
		
For the integral $(I)(w)$, we have immediately that for any fixed $\zeta$, 
$$
\chi_{\Omega\setminus B_{2\, |z-w|}(z)}(\zeta)\,  (g_-(\zeta)-g_-(w))\, b_w(\zeta)\to_{w\to z}\chi_\Omega(\zeta)\,  (g_-(\zeta)-g_-(z))\, b_z(\zeta),
$$ 
and also 
\begin{equation*}
%\begin{split}
|\chi_{\Omega\setminus B_{2\, |z-w|}(z)}(\zeta)\,  (g_-(\zeta)-g_-(w))\, b_w(\zeta)|
\le \|g_-\|_{\Lip{(\gamma,\, \bar\Omega)}}\frac{2}{|\zeta-z|^{2-\gamma}},
%\end{split}
\end{equation*} 
and by the dominated convergence theorem, 
$$
(I)(w)\to_{w\to z}\int_\Omega (g_--g_-(z))\, b_z.
$$ 
	
%\medskip
	
The term $(II)(w)$ is controlled by $$\|g_-\|_{\Lip{(\gamma,\, \bar\Omega)}}\, \int_{\Omega\cap B_{2\, |z-w|}(z)} \frac{1}{|\zeta-w|^{2-\gamma}}\, dm(\zeta)$$ and the last integral is bounded by
$$\int_{\Omega\cap B_{3\, |z-w|}(w)} \frac{1}{|\zeta-w|^{2-\gamma}}\, dm(\zeta)
\le2\pi\, \int_0^{3\, |z-w|}\frac{1}{r^{1-\gamma}}\, dr=2\pi\, 3^\gamma\, |z-w|^\gamma,
$$ and then
$$
(II)(w)\to_{w\to z}0.
$$
	
%	\bigskip
%	\bigskip
	
Also, if $\varrho\in\mathcal{D}(\C)$ such that $\varrho\equiv1$ in a ball containing $\bar\Omega$, then
\begin{equation*}
\begin{split}
\bar B[g_+](w)&=\bar B[(1-\varrho)\, g_+](w)+\bar B[\varrho\, g_+](w)\\*[5pt]
&=(IV)(w)+\int_{\C\setminus(\Omega\cup B_{2\, |z-w|}(z))} \varrho\, ( g_+-g_+(z))\, b_w\\*[5pt]
&\quad+\int_{\Omega^c\cap B_{2\, |z-w|}(z)} \varrho\, (g_+- g_+(z))\, b_w+g_+(z)\, \int_{\Omega^c} \varrho\, b_w\\*[5pt]
&=(IV)(w)+(V)(w)+(VI)(w)+g_+(z)\, (VII)(w).
\end{split}
\end{equation*} 

%\bigskip
	
It is immediate that 
$$
\lim_{w\to z}(IV)(w)=\bar B[(1-\varrho)\, g_+](z).
$$
	
%\medskip
	
By arguments similar to those used for $(I)(w)$, we have that 
$$
\lim_{w\to z}(V)(w)=\int_{\bar\Omega^c} \varrho\, (g_+-g_+(z))\, b_z,
$$ 
and, analogously to $(II)(w)$, we have that 
$$
\lim_{w\to z}(VI)(w)=0.
$$

\begin{figure}[ht]
\centering	
\includegraphics[width=0.5\linewidth]{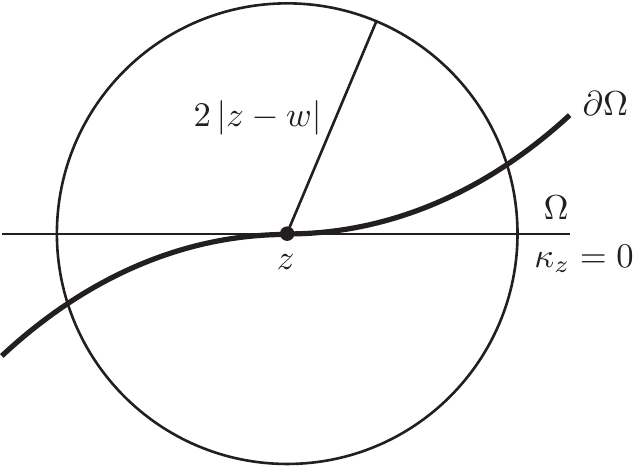}
\caption{}
\label{figura}
\end{figure}

%\bigskip
	
So, as all limits exist, we have 
\begin{equation*}
\begin{split}
\lim_{w\to z;\, w\in\Omega}& \chi_\Omega\,  \bar B[g](w)\\*[5pt]
&=\int_\Omega (g_--g_-(z))\, b_z+\bar B[(1-\varrho)\, g_+](z)\\*[5pt]
&\quad+\int_{\bar\Omega^c} \varrho\, (g_+-g_+(z))\, b_z+g_-(z) \lim_{w\to z}\!\int_\Omega b_w\!+\!g_+(z) \lim_{w\to z}\!\int_{\bar\Omega^c} \varrho\, b_w\!\\*[5pt]
&=\bar B[g_-](z)+\bar B[g_+](z)+g_-(z)\biggl\{\lim_{w\to z}\int_\Omega b_w-\int_\Omega b_z\biggr\}\\*[5pt]
&\quad+g_+(z)\biggl \{\lim_{w\to z}\int_{\bar\Omega^c} \varrho\, b_w-\int_{\bar\Omega^c} \varrho\, b_z\biggr\}\\*[5pt]
&=\bar B[g](z)+g_-(z)\biggl \{\lim_{w\to z}\int_\Omega b_w-\int_\Omega b_z\biggr\}\\*[5pt]
&\quad+g_+(z)\biggl\{\lim_{w\to z}\int_{\bar\Omega^c} \varrho\, b_w-\int_{\bar\Omega^c} \varrho\, b_z\biggr\}.
\end{split}
\end{equation*} 
%	\bigskip
%	
%	\bigskip
	
\item If $w\notin\bar\Omega$, then, in a similar way, we have
 \begin{equation*}
\begin{split}
 \lim_{w\to z;\, w\in\bar\Omega^c} \chi_{\bar\Omega^c}\,  \bar B[g](w)&=\bar B[g](z)+g_-(z)\biggl \{\lim_{w\to z}\int_{\Omega} b_w-\int_\Omega b_z\biggr\}\\*[5pt]
&\quad+g_+(z)\biggl\{\lim_{w\to z}\int_{\bar\Omega^c} \varrho\, b_w-\int_{\bar\Omega^c} \varrho\, b_z\biggr\}.
\end{split}
\end{equation*} 
%	\bigskip
%	
%	\bigskip
	
\item Then, for any $z\in\partial\Omega$ we have 
 \begin{equation*}
\begin{split}
&\frac{1}{2}\biggl \{\lim_{w\to z;\, w\in\Omega} \chi_\Omega\,  \bar B[g](w)+\lim_{w\to z;\, w\in\C\setminus\Omega} \chi_{\C\setminus \bar\Omega}\, \bar B[g](w)\biggr\}=\bar B[g](z)\\*[5pt]
&\quad+\frac{1}{2}\, g_-(z)\biggl \{\lim_{w\to z;\, w\in\Omega}\int_\Omega b_w-\int_\Omega\, b_z\biggr\}\\*[5pt]
&\quad+\frac{1}{2}\, g_+(z)\biggl\{\lim_{w\to z;\, w\in\Omega}\int_{\bar\Omega^c}\varrho\,  b_w-\int_{\bar\Omega^c}\varrho\,  b_z\biggr\}\\*[5pt]
&\quad+\frac{1}{2}\, g_-(z)\biggl \{\lim_{w\to z;\, z\notin\bar\Omega}\int_{\Omega} b_w-\int_\Omega b_z\biggr\}\\*[5pt]
&\quad+
	\frac{1}{2}\, g_+(z)\biggl\{\lim_{w\to z;\, z\notin\bar\Omega}\int_{\bar\Omega^c}\varrho\,  b_w-\int_{\bar\Omega^c}\varrho\,  b_z\biggr\}\\*[5pt]
&=\bar B[g](z)+g_-(z)\biggl(\frac{1}{2}\biggl \{\lim_{w\to z;\, w\in\Omega}\int_\Omega b_w+\lim_{w\to z;\, z\notin\bar\Omega}\int_{\Omega} b_w\biggr\}-\int_\Omega b_z\biggr)\\*[5pt]
&\quad+g_+(z)\biggl (\frac{1}{2}\biggl \{\lim_{w\to z;\, w\in\Omega}\int_{\bar\Omega^c}\varrho\,  b_w+\lim_{w\to z;\, z\notin\bar\Omega}\int_{\bar\Omega^c}\varrho\,  b_w\biggr\}-\int_{\bar\Omega^c}\varrho\,  b_z\biggr).
\end{split}
\end{equation*} 
	
And the lemma below finishes the proof of equation \eqref{jump}.
\end{itemize}

%\bigskip
%\bigskip

\begin{lem} Let $W\subset\C$ a domain with compact $\ka{1,\gamma}$ boundary.
	
Let $z\in\partial W$ and $\eta=\eta(z)$ be the normal exterior vector at $z$.  Consider the points $w=z\pm\lambda\, \eta$. 
	
For $h\in\mathcal{D}(\C)$ we have 
$$
\frac{1}{2}\biggl \{\lim_{\lambda\to 0;\, w=z-\lambda\, \eta}\int_{W}h\,  b_w+\lim_{\lambda\to 0;\, w=z+\lambda\, \eta}\int_{W}h\,  b_w\biggr\}=\int_{W}h\,  b_z.
$$
\end{lem}

\begin{proof} 
First of all, for $w\notin\partial W$, using Stokes theorem we have  
 \begin{equation*}
\begin{split}
\int_{W}h\,  b_w&=\frac{1}{2i}\, \int_{W}h(\zeta)\, \frac{d\bar\zeta\wedge d\zeta}{(\bar\zeta-\bar w)^2}\\*[5pt]
&=\frac{1}{2\, i}\biggl \{\int_{W}\frac{\partial h}{\partial\bar\zeta}(\zeta)\, \frac{d\bar\zeta\wedge d\zeta}{\bar\zeta-\bar w}-\int_{\partial W}h(\zeta)\, \frac{d\zeta}{\bar\zeta-\bar w}\biggr\}=(*).
\end{split}
\end{equation*} 
	
The integral 
$$
\int_{\partial W}h(\zeta)\, \frac{d\zeta}{\bar\zeta-\bar w}=\int_{\partial W\setminus B_{2\, |z-w|}(z)}h(\zeta)\, \frac{d\zeta}{\bar\zeta-\bar w}+\int_{\partial W\cap B_{2\, |z-w|}(z)}h(\zeta)\, \frac{d\zeta}{\bar\zeta-\bar w},
$$ 
and 	
 \begin{equation*}
\begin{split}
\int_{\partial W\cap B_{2\, |z-w|}(z)} \frac{d\zeta}{\bar\zeta-\bar w}&=\int_{\{\kappa_z=0\}\cap B_{2\, |z-w|}(z)}\frac{d\zeta}{\bar\zeta-\bar w}\\*[5pt]
&\quad+\biggl\{\int_{\partial W\cap B_{2\, |z-w|}(z)} \frac{d\zeta}{\bar\zeta-\bar w}-\int_{\{\kappa_z=0\}\cap B_{2\, |z-w|}(z)} \frac{d\zeta}{\bar\zeta-\bar w}\biggr\}\\*[5pt]
&=(I)+(II),
\end{split}
\end{equation*} 
and by the Stokes formula the second term is 
 \begin{equation*}
\begin{split}
(II)&=-\biggl\{\int_{\{\kappa_z<0\}\cap W^c\cap B_{2\, |z-w|}(z)}+\int_{\{\kappa_z>0\}\cap W\cap B_{2\, |z-w|}(z)}\biggr\}\,  \frac{d\bar\zeta\wedge d\zeta}{(\bar\zeta-\bar w)^2}\\*[5pt]
&\quad-\biggl\{\int_{\{\kappa_z<0\}\cap W^c\cap \partial B_{2\, |z-w|}(z)}+\int_{\{\kappa_z>0\}\cap W\cap \partial B_{2\, |z-w|}(z)}\biggr\} \frac{d\zeta}{\bar\zeta-\bar w}.
\end{split}
\end{equation*} 

%\bigskip  
	                                                                     
All the integrals are well defined for  $w=z\pm\lambda\, \eta$, and the geometric cancellation lemma and the facts that both 
$\{\kappa_z>0\}\cap W\cap \partial B_{2\, |z-w|}(z)$ and $\{\kappa_z<0\}\cap W\cap \partial B_{2\, |z-w|}(z)$ tend to $0$ as $r$ tends to $0$, imply that 
$$
\lim_{w\to z}(II)=0.
$$ 
	
After a translation and a rotation, 
 \begin{equation*}
(I)=\int_{-2\, |z-w|}^{2\, |z-w|}\frac{dx}{x-i\lambda},
\end{equation*} 
and then, by a direct computation 
$$
\lim_{w\to z}(I)=i\, \pi.
$$
	
%\bigskip
	
Then, 
 \begin{multline*}
\frac{1}{2}\biggl \{\lim_{\lambda\to 0;\, w=z-\lambda\, \eta}\int_{W}h\,  b_w+\lim_{\lambda\to 0;\, w=z+\lambda\, \eta}\int_{W}h\,  b_w\biggr\}\\*[5pt]
=\frac{1}{2\, i}\biggl \{\int_{W}\frac{\partial h}{\partial\bar\zeta}(\zeta)\, \frac{d\bar\zeta\wedge d\zeta}{\bar\zeta-\bar w}-\frac{1}{2}\, \biggl(\lim_{\lambda\to 0}\int_{\partial W\setminus B_{2\, \lambda}(z)}h(\zeta)\, \frac{d\zeta}{\bar\zeta-\overline{z-i\, \lambda} }\\*[5pt]
+\lim_{\lambda\to 0}\int_{\partial W\setminus B_{2\, \lambda}(z)}h(\zeta)\, \frac{d\zeta}{\bar\zeta-\overline{z+i\, \lambda}}\biggr)\biggr\}=(**) 
\end{multline*} 
because in the second term 
$$
\lim_{w\to z}(I)=-i\, \pi.
$$ 
Then 
 \begin{equation*}
\begin{split}
(**)&=\frac{1}{2\, i}\biggl \{\int_{W}\frac{\partial h}{\partial\bar\zeta}(\zeta)\, \frac{d\bar\zeta\wedge d\zeta}{\bar\zeta-\bar w}\\*[5pt]
&\hspace*{1.75cm}-\frac{1}{2}\, \lim_{\lambda\to 0}\int_{\partial W\setminus B_{2\, \lambda}(z)}h(\zeta)\biggl(\frac{1}{\bar\zeta-\overline{z-i\, \lambda} }+\frac{1}{\bar\zeta-\overline{z+i\, \lambda}}\biggr)\, d\zeta\biggr\}\\*[5pt]
&=\frac{1}{2\, i}\biggl \{\int_{W}\frac{\partial h}{\partial\bar\zeta}(\zeta)\, \frac{d\bar\zeta\wedge d\zeta}{\bar\zeta-\bar w}-\frac{1}{2}\, \lim_{\lambda\to 0}\int_{\partial W\setminus B_{2\, \lambda}(z)}h(\zeta)\, \frac{2\, (\bar\zeta-\bar z)}{(\bar\zeta-\bar z)^2+\lambda^2}\, d\zeta\biggr\}.
\end{split}
\end{equation*}
	
%\bigskip
%\bigskip
	
On the other hand, 
$$
\int_{W}h\,  b_z=\lim_{\lambda\to0}\int_{W\setminus B_{2\, \lambda}(z)}\frac{h(\zeta)}{(\bar\zeta-\bar z)^2}\,  dm(z)=\lim_{\lambda\to0} I_\lambda,
$$ 
and by Cauchy--Green's formula (cf.\ \cite{Hor1}) 
 \begin{equation*}
\begin{split}
I_\lambda&=\frac{1}{2\, i}\biggl \{-\int_{\partial W\setminus B_{2\, \lambda}(z)}h(\zeta)\, \frac{d\zeta}{\bar\zeta-\bar z}-\int_{\partial B_{2\, \lambda}(z)\cap W}h(\zeta)\, \frac{d\zeta}{\bar\zeta-\bar z}\\*[5pt]
&\hspace*{5cm}+\int_{W\setminus B_{2\, \lambda}(z)}\frac{\partial h}{\partial\bar\zeta}(\zeta)\, \frac{d\bar\zeta\wedge d\zeta}{\bar\zeta-\bar z}\biggr\},
\end{split}
\end{equation*}
	and as in the previous developments 
$$I_\lambda\to_{\lambda\to0}\frac{1}{2\, i}\biggl \{-\text{p.\ v.\ } \int_{\partial W}h(\zeta)\, \frac{d\zeta}{\bar\zeta-\bar z}+ \int_{W}\frac{\partial h}{\partial\bar\zeta}(\zeta)\, \frac{d\bar\zeta\wedge d\zeta}{\bar\zeta-\bar z}\biggr\},
$$ where the fact that 
$$
\int_0^\pi e^{2\, i\, \theta}\, d\theta=0
$$ 
is used. 	
\end{proof}

%\bigskip
%\bigskip
%
%
%\bigskip
%\bigskip

\subsection{Uniform estimates} 

Let us assume WLOG that $0\in\Omega$ and also that $R_0$ in Lemma~\ref{geom} is not larger than $1$. 

The parameter 
$$
\vartriangle(z)=\max\{|z|^2,\, d(z)^2\},
$$ 
for $z\in(\overline{\Omega\cup U_{R_0}})^c$ plays an important role in the estimates.

%\bigskip

We have

\begin{prop}\label{thm1} 
Let $\phi\in\ka{\infty}(\Omega)\cap \Lip{(\gamma,\bar\Omega)}$.
Let $\psi\in\ka{\infty}(\C\setminus\bar\Omega)\cap \Lip{(\gamma,\C\setminus\Omega)}$, satisfying that for a fixed constant $C(\psi)$, 
$$
|\psi(z)|\le\frac{C(\psi)}{\delta(z)^2},
$$ 
for $z\in(\Omega\cup U_{R_0})^c$.
	
%\medskip

Then there exists a constant $K=K(\gamma,\, \Omega,\, R_0)$ such that 
\begin{enumerate}
\item For $f=\phi$ or $\psi$, we have
$$
\|\chi_\Omega\, B[f]\|_\infty\le K\, \|f\|_\gamma.
$$
		
%\medskip
		
\item For $f$ as above and $z\in \Omega^c\cap U_{R_0}$, we have 
$$
|\chi_{\C\setminus\bar\Omega}\, B[f](z)|\le K\, \|f\|_\gamma.
$$

%\medskip
		
\item For $z\in\C\setminus(\overline{\Omega\cup U_{R_0}})^c$, we have 
$$
|(\chi_{\C\setminus\bar\Omega}\, B[\phi])(z)|\le K\, \frac{\|\phi\|_\gamma}{\delta(z)^2}.
$$
		
%\medskip
		
\item For $z\in\C\setminus(\overline{\Omega\cup U_{R_0}})^c$, we have 
$$
|(\chi_{\C\setminus\bar\Omega}\, B[\psi])(z)|\le K\, (1+\|\psi\|_\gamma)\, C(\psi)\biggl\{\frac{1}{\delta(z)^2}+ \frac{1}{\vartriangle(z)^2}\, (1+\ln(\vartriangle(z))\biggr\}.
$$
\end{enumerate} 	
\end{prop}

%\bigskip\bigskip

\begin{rmk}
In particular, if $f=\chi_\Omega$ then  
$$
\|\chi_\Omega\, B[\chi_\Omega]\|_\infty\le K
$$ 
and 
$$
|\chi_{\C\setminus\bar\Omega}\, B[\chi_\Omega](z)|\le K\biggl( \chi_{U_{R_0}}(z)+\frac{\, \chi_{\C\setminus(U_{R_0}\cup\Omega)}(z)}{\delta(z)^2}\biggr).
$$	
\end{rmk}

%\bigskip\bigskip

\subsubsection{Proof of Proposition \ref{thm1}:}
 After Proposition \ref{desing}, we have to estimate $Q[f](z)$ and $L[f](z)$ in different situations depending on the support of $f$ and the position of the point $z$.

%\medskip

For doing so, we have

\pagebreak

{\bf 1) Estimates for $L$:}

\begin{prop} 
If $z\in\Omega\setminus U_{R_0}$, there exists a constant $C=C(\gamma,\, \Omega)>0$ such that
$$
|L[\phi](z)|\le C\, \|\phi\|_\gamma
$$ 
and
$$
L[\psi](z)=0.
$$
\end{prop}

\begin{proof} 
We have now that $z\in\Omega$, $\delta(z)=d(z)>R_0$.
	
For the first part, 
$$
|L[\phi](z)|\le \|\phi\|_\gamma\, 2\pi\, \int_0^{d(z)}r^{\gamma-1}\, dr=\|\phi\|_\gamma\, 2\pi\,\frac{d(z)^\gamma}{\gamma}\le C\, \|\phi\|_\gamma .
$$
	
%\medskip
For the second part, we have $L[\psi](z)=0$. 	
\end{proof}
%\bigskip
%\bigskip

\begin{prop}\label{pu1} 
If $z\in \overline{U_{R_0}}$, there exists a constant $C=C(\gamma,\, \Omega,\, R_0)>0$ such that both for $f=\phi,\, \psi$, we have 
$$
|L[f](z)|\le C\, \|f\|_\gamma.
$$	
\end{prop}

\begin{proof} 
For $z\in\overline{U_{R_0}}$, let $\tau\in\partial\Omega$ such that $d(z)=|z-\tau|$ and  $\kappa_\tau$ defining the half space determined by the tangent line to $\partial\Omega$ across $\tau$ and 
containing the inward normal vector to $\partial\Omega$ in $\tau$. Then 
\begin{equation*}
\begin{split}
&L[f](z)=\biggl(\int_{B_{\delta(z)}(z)\cap\{\kappa_\tau\le0\}}+\int_{B_{\delta(z)}(z)\cap\{\kappa_\tau\geq0\}}\biggr)
\frac{f(\zeta)-f(z)}{|\zeta-z|^\gamma}\, 
\frac{|\zeta-z|^\gamma}{(\bar\zeta-\bar z)^2}\, dm(\zeta)\\*[5pt]
&\!=\!\biggl(\int_{B_{\delta(z)}(z)\cap\{\kappa_\tau\le0\}\cap\{\rho\le0\}}
\!+\!\int_{B_{\delta(z)}(z)\cap\{\kappa_\tau\le0\}\cap\{\rho>0\}}\biggr)\, \frac{f(\zeta)\!-\!f(z)}{|\zeta-z|^\gamma}\, 
\frac{|\zeta-z|^\gamma}{(\bar\zeta-\bar z)^2}\, dm(\zeta)\\*[5pt]
&\quad\!\!+\!\biggl(\int_{B_{\delta(z)}(z)\cap\{\kappa_\tau\geq0\}\cap\{\rho\le0\}\!}\!+\!
\int_{B_{\delta(z)}(z)\cap\{\kappa_\tau\geq0\}\cap\{\rho>0\}}\biggr)
\frac{f(\zeta)\!-\!f(z)}{|\zeta-z|^\gamma}\, \frac{|\zeta-z|^\gamma}{(\bar\zeta\!-\!\bar z)^2}\, dm(\zeta).\!\!\!
\end{split}
\end{equation*}	

%\medskip
	
If $f=\phi$, in the case of $z\in\Omega$, we have 
\begin{equation*}
\begin{split}
L[\phi](z)&=\int_{B_{\delta(z)}(z)\cap\{\kappa_\tau\le0\}\cap\{\rho\le0\}}\, 
	\frac{\phi(\zeta)-\phi(z)}{|\zeta-z|^\gamma}\, 
	\frac{|\zeta-z|^\gamma}{(\bar\zeta-\bar z)^2}\, dm(\zeta)\\*[5pt]
&\quad-\phi(z)\, \int_{B_{\delta(z)}(z)\cap\{\kappa_\tau\le0\}\cap\{\rho>0\}}
	\frac{1}{(\bar\zeta-\bar z)^2}\, dm(\zeta)\\*[5pt]
&\quad+\int_{B_{\delta(z)}(z)\cap\{\kappa_\tau\geq0\}\cap\{\rho\le0\}}\, \frac{\phi(\zeta)-\phi(z)}{|\zeta-z|^\gamma}\, 
	\frac{|\zeta-z|^\gamma}{(\bar\zeta-\bar z)^2}\, dm(\zeta)\\*[5pt]
&\quad-\phi(z)\, 
	\int_{B_{\delta(z)}(z)\cap\{\kappa_\tau\geq0\}\cap\{\rho>0\}} 
\frac{1}{(\bar\zeta-\bar z)^2}\, dm(\zeta)\\*[5pt]
&=(I)+(II)+(III)+(IV).
\end{split}
\end{equation*}
	
The integral 
$$
|(I)|\le C\, \|\phi\|_\gamma
$$ 
as in the previous proposition.
	
%\medskip
	
Concerning the terms $(II)$ and $(III)$, the integrals there are extended to the subsets of $B_{\delta(z)}(z)$ located 
between the tangent line to $\partial\Omega$ at $\tau$ and $\partial\Omega$ itself. Also $z$ is not in the domain of integration. Then, using the Lemma \ref{Aux1} below, we have
$$
|(II)|\le\|\phi\|_\infty\, \frac{4\, \delta(z)^\gamma}{\gamma}\le C\, \|\phi\|_\infty,
$$ 
and the integral 
$$
|(III)|\le \|\phi\|_\gamma\, \frac{4\, \delta(z)^\gamma}{\gamma}\le C\, \|\phi\|_\gamma.
$$
	
%\medskip
	
Finally, the integral in $(IV)$ can be written as 
$$
\biggl(\int_{B_{\delta(z)}(z)\cap\{\kappa_\tau\geq0\}}-
\int_{B_{\delta(z)}(z)\cap\{\kappa_\tau\geq0\}\cap\{\rho\le0\}}\biggr)\, b_z,
$$ 
and using Lemmas \ref{Aux1} and \ref{Aux2} we have the estimate 
$$
|(IV)|\le\|\phi\|_\infty\biggl(\frac{\pi}{2} +\frac{4\, \delta(z)^\gamma}{\gamma}\biggr).
$$
	
Finally, in this case $\frac{R_0}{2}\le\delta(z)\le R_0$, and this concludes.
\end{proof}	

%\medskip
	
The case of $f=\psi$ is completely similar.

%\bigskip
%\bigskip

\begin{lem}\label{Aux1} 
Under the conditions and notation of Proposition \ref{pu1}, the integrals 
$$ 
\int_{B_{\delta(z)}(z)\cap\{\kappa_\tau\le0\}\cap\{\rho>0\}} b_z
$$ 
and 
$$
\int_{B_{\delta(z)}(z)\cap\{\kappa_\tau\geq0\}\cap\{\rho\le0\}} b_z
$$ 
are bounded by 
$$
\frac{4\, R_0^\gamma}{\gamma}.
$$	
\end{lem}

\begin{proof}
Concerning the terms $(II)$ and $(III)$, the integrals are extended to the subsets of $B_{\delta(z)}(z)$ located 
between the tangent line to $\partial\Omega$ at $w_0$ and $\partial\Omega$ itself (we change the notation and use $w_0$ instead of $\tau$). 
	
Then, we can take coordinates centered at $w_0$ given by the frame $\eta(w_0)=\frac{2\, \bar\partial\rho(w_0)}{\|\nabla\rho(w_0)\|}$,
$\tau(w_0)=i\, \eta(w_0)$, so any point $z\in\C$ is of the form 
$$
z=w_0+\lambda \eta(w_0)+\mu\, \tau(w_0)
$$
and 
$$
\kappa_{w_0}(z)=\lambda\, \|\nabla\rho(w_0)\|
$$ 
and 
$$
\rho(z)=\kappa_{w_0}(z)
+\omega(\sqrt{\lambda^2+\mu^2})\, \sqrt{\lambda^2+\mu^2},
$$ 
where $\omega$ is the modulus of continuity of the derivatives of $\rho$, 
at $w_0$.
	
%\medskip
	
Then the isometric map
$$
P\colon \C\rightarrow\C
$$ 
given by 
$$
\phi(\mu+i\lambda)=w_0+\lambda \eta(w_0)+\mu\, \tau(w_0)
$$ 
transforms $0$ in $w_0$, $-i\, d(z_0)$  in $z_0$, the real and the imaginary axis in the lines across $w_0$ directed by $\eta(w_0)$ and $\tau(w_0)$ respectively. Also 
$$
P(B_{\delta(z)}(-i\, d(z)))=B_{\delta(z)}(z).
$$
	
Then, if
$$
A_1\!=\!\{\mu+i\, \lambda\in\C: |\mu|\!\le\!\sqrt{R_0^2-d(z)^2},\, \lambda\, \|\nabla\rho(w_0)\|+
\omega(\sqrt{\lambda^2+\mu^2})\, \sqrt{\lambda^2+\mu^2}\le0\}\!\!
$$ 
and 
$$
A_2\!=\!\{\mu+i\, \lambda\!\in\!\C: |\mu|\!\le\!\sqrt{R_0^2-d(z)^2},\,  \lambda\, \|\nabla\rho(w_0)\|+
\omega(\sqrt{\lambda^2+\mu^2}) \, \sqrt{\lambda^2+\mu^2}>0\},\!\!
$$
we have  that 
\begin{equation*}
\begin{split}	
(II)&=\int_{B_{\delta(z)}(-i\, d(z))\cap\{\lambda>0\}\cap A_1}\, 
\frac{1}{(\mu-i\, \lambda+i\, d(z))^2}\, dm(\mu,\lambda)\\*[5pt]
&=\int_{B_{\delta(z)}(-i\, d(z))\cap\{\lambda>0\}\cap A_1}\, 
\frac{1}{(\mu-i\, \lambda+i\,d(z))^2}\,dm(\mu,\lambda)\\*[5pt]
&=\int_{-\sqrt{R_0^2-d(z)^2}}^{\sqrt{R_0^2-d(z)^2}}\, d\mu\, 
\int_0^{\varphi_{+}(\mu)}\, \frac{1}{(\mu-i\,\lambda+i\,d(z))^2}\, d\lambda\\*[5pt]
&=\frac{1}{i}\, \int_{-\sqrt{R_0^2-d(z)^2}}^{\sqrt{R_0^2-d(z)^2}}\, 
\frac{\varphi_{+}(\mu)}{\mu-i\,\varphi_{+}(\mu)+i\,d(z)}\, d\mu.
\end{split}
\end{equation*}	

%\medskip
	
Then 
\begin{equation*}
\begin{split}	
|(II)|&\le\int_{-\sqrt{R_0^2-d(z)^2}}^{\sqrt{R_0^2-d(z)^2}}\, 
\frac{|\varphi_{+}(\mu)|}{|\mu-i\,(\varphi_{+}(\mu)+\,d(z))|}\, d\mu\\[5pt]
&\le C\, \int_{-\sqrt{R_0^2-d(z)^2}}^{\sqrt{R_0^2-d(z)^2}}\, 
\frac{|\mu|^\gamma}{\sqrt{\mu^2+ (\varphi_{+}(\mu)+ d(z))^2}}\, d\mu\\[5pt]
&\le 2C\, \int_{0}^{\sqrt{R_0^2-d(z)^2}}\, 
\frac{1}{\mu^{1-\gamma}}\, d\mu=\frac{2C}{\gamma}\biggl(\sqrt{R_0^2-d(z)^2}\biggr)^\gamma.\qedhere
\end{split}
\end{equation*}	
\end{proof}

%\bigskip
%\bigskip

\begin{lem}\label{Aux2} 
$$
\biggl|\int_{B_{\delta(z)}(z)\cap\{\kappa_\tau\geq0\}}b_z\biggr|\le\frac{\pi}{2}.
$$
\end{lem}

\begin{proof} 
It is clear that if $\delta(z)\geq\frac{R_0}{2}$ or if $z=\tau$, then 
$$
\int_{B_{\delta(z)}(z)\cap\{\kappa_\tau\geq0\}}b_z=0
$$ 
by the cancellation property. So we assume that $0<d(z)<\frac{R_0}{2}$ and then our integral is 
$$
\int_{B_{\frac{R_0}{2}}(z)\cap\{\kappa_\tau\geq0\}}b_z.
$$ 
Also by the cancellation property, and after a rigid movement, our integral is  
\begin{equation*}
\begin{split}	
\int_{B_{\delta(z)}(0)\cap\{\Im(\zeta)\geq\alpha\}}b_z&=\int_{\alpha}^{\frac{R_0}{2}}\frac{dr}{r}\,
	\int_{\arcsin(\frac{\alpha}{r})}^{\pi-\arcsin(\frac{{\alpha}}{r})} e^{2i\theta}\,d\theta\\*[5pt]
&\!=\!\int_{\alpha}^{\frac{R_0}{2}}\frac{dr}{r}\,\frac{1}{2i}\,\left(e^{-2i \arcsin(\frac{\alpha}{r})}-
	e^{2i \arcsin(\frac{\alpha}{r})}\right)\\[5pt]
&\!=\!\int_{\alpha}^{R_0}\!\frac{dr}{r}\frac{e^{-i \arcsin(\frac{\alpha}{r})}\!-\!
e^{i \arcsin(\frac{\alpha}{r})}}{2i}\left(e^{-i \arcsin(\frac{\alpha}{r})}\!+\!
e^{i \arcsin(\frac{\alpha}{r})}\right)\!\!\!\\[5pt]
&\!=\!-2\,\int_{\alpha}^{\frac{R_0}{2}}\frac{dr}{r}\,\frac{\alpha}{r}\,
\Re(e^{i \arcsin(\frac{\alpha}{r})})\\[5pt]
&\!=\!-2\,\int_{\alpha}^{\frac{R_0}{2}}\frac{dr}{r}\,\frac{\alpha}{r}\,
\sqrt{1-\biggl(\frac{\alpha}{r}\biggr)^2}=(*),
\end{split}
\end{equation*} 
so we have, for $t=\frac{\alpha}{r}$, that 
\begin{equation*}
\begin{split}	
(*)&=-2\,\int_\frac{2\alpha}{R_0}^1\sqrt{1-t^2}\,dt=-(\arcsin(t)
+t\,\sqrt{1-t^2})]_\frac{2\alpha}{R_0}^1\\
&=\arcsin\biggl(\frac{2\alpha}{R_0}\biggr)
+\frac{2\alpha}{R_0}\,\sqrt{1-\biggl(\frac{2\alpha}{R_0}\biggr)^2}-\arcsin 1\le\frac{\pi}{2}.\qedhere
\end{split}
\end{equation*} 
\end{proof}

%\bigskip
%\bigskip

\begin{prop} 
If $z\in\C\setminus\overline{U_{R_0}\cup\Omega}$, then 
$$
L[\phi](z)=0,
$$ 
and there exists a constant $C_3=C_3(\gamma,\, \Omega,\, R_0)$ such that 
\begin{multline*}
|L[\psi](z)|\le C_3\, (1+\|\psi\|_\gamma) \, C(\psi)\biggl\{\frac{1}{\max\{R_0^2,\,d(z)^2\}}\\*[5pt]
+\frac{1}{\max\{d(z)),\,  |z|\}^2}\, \{1+\ln(\max\{d(z),\,  |z|\})\biggr\}.
\end{multline*} 	
	
%\medskip
	
Here we are assuming that $\|\psi\|_\gamma>0$ and also that $0\in\Omega\setminus U_{R_0}$.	
\end{prop}

\begin{proof} 
The case of $\phi$ is immediate.
	
%\medskip

In the case of $\psi$, for any choice of a positive $\alpha(z)<\delta(z)$, we have 
\begin{equation*}
\begin{split}	
L[\psi](z)&=\int_{B_{\delta(z)}(z)}\frac{\psi(\zeta)-\psi(z)}{|\zeta-z|^\gamma}\, 
	\frac{|\zeta-z|^\gamma}{(\bar\zeta-\bar z)^2}\, dm(\zeta)\\*[5pt]
&=\biggl(\int_{B_{\alpha(z)}(z)}+\int_{B_{\delta(z)}(z)\setminus B_{\alpha(z)}(z)}\biggr)\frac{\psi(\zeta)-\psi(z)}{|\zeta-z|^\gamma}\, 
	\frac{|\zeta-z|^\gamma}{(\bar\zeta-\bar z)^2}\, dm(\zeta)\\*[5pt]
&=(I)_{\alpha(z)}+(II)_{\alpha(z)}.
\end{split}
\end{equation*}
	
%\medskip
	
Since 
$$
|\psi(\zeta)|\le|\psi(z)|+\|\psi\|_\gamma\, |\zeta-z|^\gamma
$$ 
then for $\zeta\in B_{\alpha(z)}(z)$ we have
 $$
 |\psi(\zeta)|\le|\psi(z)|+\|\psi\|_\gamma\, \alpha(z)^\gamma,
 $$ 
 and if $\alpha(z)\le(\max\{\frac{|\psi(z)|}{\|\psi\|_\gamma},\, \frac{1}{\max\{R_0^2,\,d(z)^2\}}\})^\frac{1}{\gamma} $, we have
$$
|\psi(\zeta)|\le 2\, |\psi(z)|
$$ 
and then
\begin{equation*}
\begin{split}
|(I)_{\alpha(z)}|&\le \|\psi\|_\gamma\, \int_{B_{\alpha(z)}(z)} \frac{dm(\zeta)}{|\zeta-z|^{2-\gamma}}\\*[5pt]
&=2\pi\, \|\psi\|_\gamma\, \int_0^{\alpha(z)} \frac{dr}{r^{1-\gamma}}\\*[5pt]
&\le 2\pi\, \|\psi\|_\gamma\, \frac{\alpha(z)^{\gamma}}{\gamma} 
=\frac{2\pi}{\gamma}\, |\psi(z)|.
\end{split}
\end{equation*}
	
%\medskip
	
The term 
\begin{equation*}
\begin{split}
(II)_{\alpha(z)}&=\biggl(\int_{B_{d(z)}(z)\cap U_{R_0}}+
	\int_{(B_{d(z)}(z)\setminus B_{\alpha(z)}(z))\setminus U_{R_0}}\biggr)
\frac{\psi(\zeta)-\psi(z)}{|\zeta-z|^\gamma}\, 
\frac{|\zeta-z|^\gamma}{(\bar\zeta-\bar z)^2}\, dm(\zeta)\!\\
&=(III)+(IV)_{\alpha(z)}.
\end{split}
\end{equation*}
	
%\medskip
	
The term 
\begin{equation*}
\begin{split}
(III)&=\biggl(\int_{B_{d(z)}(z)\cap U_{\frac{R_0}{2}}}+
\int_{B_{d(z)}(z)\cap (U_{R_0}\setminus U_{\frac{R_0}{2}})}\biggr)
\frac{\psi(\zeta)-\psi(z)}{|\zeta-z|\gamma}\, 
\frac{|\zeta-z|^\gamma}{(\bar\zeta-\bar z)^2}\, dm(\zeta)\\*[5pt]
&=(V)+(VI).
\end{split}
\end{equation*}
	
%\bigskip
	
Since for $\zeta\in U_{\frac{R_0}{2}}$ we have $|\zeta-z|\geq\frac{d(z)}{2}$, then
$$
|(V)|\le\|\psi\|_\infty\, \frac{4\, m(\bar U_{R_0})}{\max\{R_0^2,\,d(z)^2\}}.
$$
	
%\medskip
	
The integral of the term $(VI)$, extends to $\zeta\in U_{R_0}\setminus U_{\frac{R_0}{2}}$. Then, if $d(z)\geq\frac{3R_0}{2}$ we have 
$$
|\zeta-z|\geq||\zeta-\varsigma|-|\varsigma-z||\geq d(z)-R_0\geq\max\biggl\{\frac{d(z)}{3},\, \frac{R_0}{2}\biggr\}
$$ 
so 
$$
|(VI)|\le\|\psi\|_\infty\, \frac{36\, m(\bar U_{R_0})}{\max\{R_0^2,\,d(z)^2\}}.
$$
	
%\medskip
	
If $R_0\le d(z)<\frac{3R_0}{2}$, then $d(z)-R_0\le|\zeta-z|\le d(z)$ and then 
$$
|(VI)|\le\|\psi\|_\gamma\, 2\pi\, \int_{d(z)-R_0}^{d(z)}\frac{dr}{r^{1-\gamma}}\le\|\psi\|_\gamma\, \frac{2\pi}{\gamma}\, d(z)^\gamma\le\frac{\|\psi\|_\gamma\, \frac{9\pi}{8\gamma}\, (\frac{3}{4})^\gamma\, R_0^{2+\gamma}}{\max\{R_0^2,\,d(z)^2\}}.
$$  
	
%\bigskip
	
The term 
$$
(IV)_{\alpha(z)}=\biggl(\int_{C^{\frac{d(z)}{2}}_{\alpha(z)}(z)}+
	\int_{C^{d(z)}_{\frac{d(z)}{2}}(z)\setminus U_{R_0}}\biggr)\psi\, b_z=(VII)_{\alpha(z)}
+(VIII),
$$ 
and we can use the estimate 
$$
|\psi(z)|\le\frac{C(\psi)}{\max\{R_0^2,\,d(z)^2\}}.
$$
	
In the case of $\frac{d(z)}{2}>\alpha(z)$. Otherwise the first term is equal to $0$. Since for $\zeta\in B_{\frac{d(z)}{2}}(z)$, we have $d(\zeta)>\frac{d(z)}{2}$, we have 
\begin{equation*}
\begin{split}
|(VII)_{\alpha(z)}|&\le\, C(\psi)\, 
\int_{C^{\frac{d(z)}{2}}_{\alpha(z)}(z)}\frac{dm(\zeta)}{d(\zeta)^2\,|\zeta-z|^2}\\*[5pt]
&\le\, C(\psi)\, \frac{4}{d(z)^2}\, \int_{C^{\frac{d(z)}{2}}_{\alpha(z)}(z)}\frac{dm(\zeta)}{|\zeta-z|^2}=C(\psi)\, 
\frac{8\pi}{d(z)^2}\, \int^{\frac{d(z)}{2}}_{\alpha(z)}\frac{dr}{r}\\*[5pt]
&=C(\psi)\, \frac{8\pi}{d(z)^2}\, \ln(\frac{d(z)}{2\, \alpha(z)})=(*),
\end{split}
\end{equation*}
and if $\frac{|\psi(z)|}{\|\psi\|_\gamma}\geq\frac{1}{\max\{R_0^2,\,d(z)^2\}},$ then 
$$
(*)\le\, C(\psi)\, \frac{8\pi}{d(z)^2}\, \{\ln(d(z)\, \|\psi\|_\gamma^{\frac{1}{\gamma}})+\, \ln(2\, |\psi(z)|^{\frac{1}{\gamma}})\},
$$ 
and in the oposite case we have
$$
(*)\le\, C(\psi)\, \frac{8\pi}{d(z)^2}\, \{\ln(d(z))+\ln(2\, \max\{R_0^2,\,d(z)^2\}^{\frac{1}{\gamma}}\}).
$$
	
%\medskip
	
For the term $(VIII)$, we consider 
$$
C^{d(z)}_{\frac{d(z)}{2}}(z)\setminus U_{R_0}=E_1(z)\cup E_2(z),
$$ 
where $E_1(z)$ and $E_2(z)$ are the intersection of the domain of integration with the set 
$$
\{\zeta\in\C:\, d(\zeta)\geq
	|\zeta-z|\}
$$ 
or its complement, respectively. Then 
$$
(VIII)=\biggl(\int_{E_1(z)}+\int_{E_2(z)}\biggr)\, \psi b_z=(IX)+(X),
$$ 
and the term 
$$
|(IX)|\le\, C(\psi)\, \int_{\C\setminus B_{\frac{d(z)}{2}}(z)}\frac{dm(\zeta)}{|\zeta-z|^4}\le\, C(\psi)\,  2\pi\, \int_{\frac{d(z)}{2}}^\infty \frac{dr}{r^3}\le C(\psi)\, \frac{8\pi}{d(z)^2}.
$$
	
%\medskip
	
Let 
$$
M=\max\{|w|:w\in\overline{(U_{R_0}\cup\Omega)}\}.
$$ 
The term 
$$
(X)=\biggl(\int_{E_2(z)\cap B_{2M}(0)}+\int_{E_2(z)\setminus B_{2M}(0)}\biggr)\, \psi b_z
\newline	=(XI)+(XII)
$$ 
since 
we have that if $|\zeta|\geq 2M$, then $d(\zeta)\geq\frac{|\zeta|}{2}$ and then 
$$
|(XII)|\le\, C(\psi)\,  4\int_{C^{d(z)}_{\frac{d(z)}{2}}(z)\setminus B_{2M}(0)}\frac{dm(\zeta)}{|\zeta|^2\,|\zeta-z|^2}.
$$
	
%\medskip
	
If $|z|\geq 5M$, then $d(z)\geq 3M$ and for $r<3M$, the balls $B_{2M}(0)$ and $B_{\frac{d(z)}{2}}(z)$ are 
mutually disjoint and we consider then the decomposition in disjoint sets 
$$
C^{d(z)}_{\frac{d(z)}{2}}(z)\setminus B_{2M}(0))=A_1\cup A_2\cup A_3,
$$ 
where
\begin{align*}
A_1&=C^{d(z)}_{\frac{d(z)}{2}}(z)\setminus B_{2|z|}(0),\\*[5pt]
A_2&=\{\zeta\in B_{2|z|}(0)\setminus B_{2M}(0):\, |\zeta|\le |\zeta-z|\},
\intertext{and}
A_3&=\{\zeta\in B_{2|z|}(0)\setminus B_{\frac{d(z)}{2}}(z):\, \, |\zeta|\geq |\zeta-z|\}.
\end{align*}
	
If $|\zeta|>2|z|$, then $\frac{1}{2}\, |\zeta|\le |\zeta-z|\le\frac{1}{2}\, |\zeta|$, so the integral over
$A_1$ is bounded by 
$$
\int_{C^{d(z)}_{\frac{d(z)}{2}}(z)\setminus B_{2|z|}(0)}\frac{4\, dm(\zeta)}{|\zeta|^4}=\frac{4\pi}{|z|^2}.
$$
	
%\medskip
	
Since for $\zeta\in A_2$, $|\zeta-z|\geq\frac{|z|}{2}$, the integral over $A_2$ is bounded by 
$$
\frac{4}{|z|^2}\, 
\int_{C_{2M}^{2|z|}(0)}\frac{dm(\zeta)}{|\zeta|^2}\le\frac{8\pi}{|z|^2}\, 
\int_{2M}^{2|z|}\frac{dr}{r}=\frac{8\pi}{|z|^2}\, \ln\biggl(\frac{|z|}{M}\biggr).
$$
	
%\medskip
	
Finally, for $\zeta\in A_3$ we have that $|\zeta|\geq\frac{|z|}{2}$ an the integral in $A_3$ is bounded by 
\begin{equation*}
\begin{split}
\frac{4}{|z|^2}\,\int_{B_{2|z|}(0)\setminus B_{\frac{d(z)}{2}}(z)}\frac{dm(\zeta)}{|\zeta-z|^2}&\le
	\frac{4}{|z|^2}\,\int_{B_{4|z|}(0)\setminus B_{\frac{d(z)}{2}}(z)}\frac{dm(\zeta)}{|\zeta-z|^2}\\*[5pt]
&=\frac{8\pi}{|z|^2}\, \ln\biggl(\frac{8|z|}{d(z)}\biggr).
\end{split}
\end{equation*}
	
%\bigskip
	
If $|z|\leq 5M$, then the set 
$$
C^{d(z)}_{\frac{d(z)}{2}}(z)\setminus B_{2M}(0)
$$ 
decomposes in a disjoint union of the resulting intersection with the set 
$$
\{|\zeta|>|\zeta-z|\}
$$ 
and its complement, namely $G_1$ and $G_2$. Then 
$$
\int_{G_1}\frac{dm(\zeta)}{|\zeta|^2\,|\zeta-z|^2}\le
\int_{\C\setminus B_{\frac{d(z)}{2}}(z))}\frac{dm(\zeta)}{|\zeta-z|^4}=\frac{16\pi}{d(z)^2},
$$ 
and 
$$
\int_{G_2}\frac{dm(\zeta)}{|\zeta|^2\,|\zeta-z|^2}\le
\int_{\C\setminus B_{2M}(0)}\frac{dm(\zeta)}{|\zeta|^4}=\frac{16\pi}{4M^2}\le\frac{16\pi\, 25}{|z|^2}.
$$

%\bigskip
	
The integral 
\begin{equation*}
|(XI)|\!\le\!\!\int_{\!E_2(z)\cap B_{2M}(0)}\! |\psi b_z|
\le C(\psi) \!\int_{\!B_{8M}(0)} \frac{dm(\zeta)}{R_0^2|\zeta-z|^2}\le C(\psi) \frac{64\pi M^2}{R_0^2} \frac{4}{d(z)^2}.\!\qedhere
\end{equation*}
\end{proof}

%\bigskip
%\bigskip

%\bigskip
%\bigskip

\medskip

2) {\bf Estimates for $Q$:}

\begin{prop} 
There exists a constant $C_4=C_4(\Omega,\, R_0)$ such that if $z\in\Omega\cup U_{R_0}$, we have 
$$
|Q[\phi](z)|\le C_4\, \|\phi\|_\infty
$$ 
and 
$$
|Q[\psi](z)|\le C_4\, (\|\psi\|_\infty+C(\psi)).
$$	
\end{prop}

\begin{proof} 
A) In the case of $z\in\Omega\setminus U_{R_0}$, we have $\delta(z)=d(z)\geq R_0$ and 
then 
$$
Q[\phi](z)=\int_{\Omega\setminus B_{d(z)}(z)}\phi\, b_z
$$ 
so
$$
|Q[\phi](z)|\le\|\phi\|_\infty\, \frac{m(\overline{\Omega\cup U_{R_0}})}{R_0^2}.
$$
	
%\bigskip 

Also 
\begin{equation*}
\begin{split}
Q[\psi](z)&=\int_{\C\setminus B_{d(z)}(z)}\psi\, b_z=\int_{\C\setminus\Omega}\psi\, b_z\\*[5pt]
&=\int_{(\C\setminus\Omega)\cap U_{R_0}}\psi\, b_z+\int_{\C\setminus(\Omega\cup U_{R_0})}\psi\, b_z=(I)+(II),
\end{split}
\end{equation*}
 and $(I)$ has the same control as $Q[\phi](z)$.
 
 %\pagebreak
	
The integral 
$$
|(II)|\le\, C(\psi)\, \int_{\C\setminus(\Omega\cup U_{R_0})}\frac{dm(\zeta)}{\max\{R_0^2,\,d(\zeta)^2\}\, |\zeta-z|^2},
$$ 
and the set 
$$
\C\setminus(\Omega\cup U_{R_0})=E_1\cup E_2,
$$ 
where $E_1$ is the intersection of $\C\setminus(\Omega\cup U_{R_0})$ with 
$$
\{d(\zeta)\geq |\zeta-z|\}
$$ 
and $E_2$ with the complement, we have that the integral is 
$$
\int_{E_1}\frac{dm(\zeta)}{d(\zeta)^2\, |\zeta-z|^2}
	+\int_{E_2}\frac{dm(\zeta)}{d(\zeta)^2\, |\zeta-z|^2}=(III)+(IV),
$$ 
and 
$$
(III)\le\int_{E_1}\frac{dm(\zeta)}{|\zeta-z|^4}\le\int_{\C\setminus B_{R_0}(z)}\frac{dm(\zeta)}{|\zeta-z|^4}\le \frac{2\pi}{3\, R_0^2}.
$$
		
%\medskip
	
For the integral $(IV)$, we have that in our domain, $d(\zeta)\geq R_0$ and 
$$
|\zeta-z|\le|\zeta-\tau|+|\tau-z|\le d(\zeta)+\dmt(\Omega),
$$ 
where $|\zeta-\tau|=d(\zeta)$. This implies that $E_2\subset B_{R_0+3\dmt(\Omega)}(0)$ and so 
$$
|(IV)|\le \frac{m(B_{R_0+3\dmt(\Omega)}(0))}{R_0^4}.
$$
	
%	\bigskip
%	\bigskip
	
B) In the case of $z\in U_{R_0}$, we have 
\begin{equation*}
\begin{split}
Q[f](z)&=\int_{\C\setminus B_{R_0}(z)}f\, b_z=\biggl(\int_{\C\setminus (B_{R_0}(z)\cup U_{R_0}\cup\Omega)}+\int_{(U_{R_0}\cup\Omega)\setminus B_{R_0}(z)}\biggr)\,f\, b_z\\*[5pt]
&=(I)+(II).
\end{split}
\end{equation*}

%\medskip
	
In the case $f=\phi$, the integral 
$$
(I)=0.
$$ 
	
%\medskip
	
The integral 
$$
|(II)|\le\|\phi\|_\infty\, \frac{m(\overline{\Omega\cup U_{R_0}})}{R_0^2}
$$ 
in all cases.
	
%\medskip
	
In the case of $f=\psi$, 
$$
|(I)|\le\, C(\psi)\, \int_{\C\setminus(\Omega\cup U_{R_0})}\frac{dm(\zeta)}{\max\{R_0^2,\,d(\zeta)^2\}\, |\zeta-z|^2}\le\, C(\psi)\, \frac{m(\Omega\cup U_{R_0})}{R_0^4} 
$$
and $(II)$ in A).
\end{proof}

%\bigskip
%\bigskip

\begin{prop} 
There exists a constant $C_5=C_5(\Omega,\, R_0)$ such that if  $z\in\C\setminus\overline{(\Omega\cup U_{R_0})}$, then 
$$
|Q[\phi](z)|\le\, C_5\, \|\phi\|_\infty\, \frac{1}{\max\{R_0^2,\,d(z)^2\}}
$$ 
and 
$$
|Q[\psi](z)|\le\, C_5\biggl \{\|\psi\|_\infty\, \frac{1}{\delta(z)^2}+C(\psi)\, \frac{1}{\max\{|z|^2,\,d(z)^2\}}\, (1+\ln\max\{|z|,\,d(z)\})\biggr\}.
$$	
\end{prop}

\begin{proof} 
If $z\in\C\setminus(\Omega\cup U_{R_0})$, then 
$$
Q[\phi](z)=\int_{\Omega\setminus B_{\delta(z)}(z)}\phi\, b_z
$$ 
and is bounded by 
$$
\|\phi\|_\infty\, \frac{m(\overline{\Omega\cup U_{R_0}})}{\delta(z)^2},
$$ 
leading to the statement.
	
%\medskip
	
Also 
\begin{equation*}
\begin{split}
Q[\psi](z)&=\int_{\C\setminus(\Omega\cup B_{\delta(z)}(z))}\psi\, b_z\\*[5pt]
&=\int_{\C\setminus(\Omega\cup U_{R_0}\cup B_{\delta(z)}(z))}\psi\, b_z+\int_{ U_{R_0}\setminus(\Omega\cup B_{\delta(z)}(z))}\psi\, b_z=(I)+(II).
\end{split}
\end{equation*}

%\medskip
	
The integral 
$$
|(II)|\le\|\psi\|_\infty\, \frac{m(\overline{U_{R_0}})}{\delta(z)^2}.
$$
	
%\medskip
	
The term 
$$
|(I)|\le\, C(\psi)\, \int_{\C\setminus(\Omega\cup U_{R_0}\cup B_{\delta(z)}(z))}\frac{dm(\zeta)}{\max\{R_0^2,\,d(\zeta)^2\}\, |\zeta-z|^2},
$$ 
and we consider 
$$
\C\setminus(\Omega\cup B_{d(z)}(z))=E_1(z)\cup E_2(z),
$$ 
where $E_1(z)$ and $E_2(z)$ are the intersection of the domain of integration with the set 
$$
\{\zeta\in\C:\, d(\zeta)\geq |\zeta-z|\}
$$ 
or its complement. Then the integral in $(I)$ can be decomposed as  
$$
\biggl(\int_{E_1(z)}+\int_{E_2(z)}\biggr)\, \varphi b_z=(III)+(IV),
$$ 
and the term 
$$
|(III)|\le\int_{\C\setminus B_{d(z)}(z)}\frac{dm(\zeta)}{|\zeta-z|^4}\le 2\pi\, \int_{d(z)}
	^\infty \frac{dr}{r^3}\le\frac{\pi}{d(z)^2}.
$$
	
%\medskip
	
The term 
$$
(IV)=\biggl(\int_{E_2(z)\cap B_{2M}(0)}+\int_{E_2(z)\setminus B_{2M}(0)}\biggr)\, \varphi b_z
=(V)+(VI)
$$ 
since in general $d(\zeta)\geq|\, |\zeta|-M|$,  we have that if $|\zeta|\geq 2M$, then $d(\zeta)\geq\frac{|\zeta|}{2}$ and then
$$
|(VI)|\le 4\int_{\C\setminus(B_{2M}(0)\cup B_{\frac{d(z)}{2}}(z))}\frac{dm(\zeta)}{|\zeta|^2\,|\zeta-z|^2}.
$$
	
%\medskip
	
If $|z|\geq 5M$, then $d(z)\geq 3M$ and  the balls $B_{2M}(0)$ and $B_{\frac{d(z)}{2}}(z)$ are 
mutually disjoint and we consider then the decomposition in disjoint sets 
$$
\C\setminus(B_{2M}(0)\cup B_{\frac{d(z)}{2}}(z))=A_1\cup A_2\cup A_3,
$$ 
where
\begin{align*}
A_1&=\C\setminus B_{2|z|}(0),\\*[5pt]
A_2&=\{\zeta\in B_{2|z|}(0)\setminus B_{2M}(0):\, |\zeta|\le |\zeta-z|\},
\intertext{and} 
A_3&=\{\zeta\in B_{2|z|}(0)\setminus B_{\frac{d(z)}{2}}(z):\, \, |\zeta|\geq |\zeta-z|\}.
\end{align*}

If $|\zeta|>2|z|$, then $\frac{1}{2}\, |\zeta|\le |\zeta-z|\le\frac{3}{2}\, |\zeta|$, so the integral over
$A_1$ is bounded by 
$$
\int_{\C\setminus(B_{2|z|}(0)}\frac{4\, dm(\zeta)}{|\zeta|^4}=\frac{4\pi}{|z|^2}.
$$
	
%\medskip
	
Since for $\zeta\in A_2$, $|\zeta-z|\geq\frac{|z|}{2}$, the integral over $A_2$ is bounded by 
$$
\frac{4}{|z|^2}\, \int_{C_{2M}^{2|z|}(0)}\frac{dm(\zeta)}{|\zeta|^2}
\le\frac{8\pi}{|z|^2}\, \int_{2M}^{2|z|}\frac{dr}{r}=\frac{8\pi}{|z|^2}\, \ln\biggl(\frac{|z|}{M}\biggr).
$$
	
\medskip
	
Finally, for $\zeta\in A_3$ we have that $|\zeta|\geq\frac{|z|}{2}$ an the integral in $A_3$ is bounded by 
\begin{equation*}
\begin{split}
\frac{4}{|z|^2}\,\int_{B_{2|z|}(0)\setminus B_{\frac{d(z)}{2}}(z)}\frac{dm(\zeta)}{|\zeta-z|^2}&\le
\frac{4}{|z|^2}\,\int_{B_{4|z|}(0)\setminus B_{\frac{d(z)}{2}}(z)}\frac{dm(\zeta)}{|\zeta-z|^2}\\*[5pt]
&=\frac{8\pi}{|z|^2}\, \ln\biggl(\frac{8|z|}{d(z)}\biggr).
\end{split}
\end{equation*}
	
%\bigskip
	
If $|z|\leq 5M$, then the set 
$$
\C\setminus(B_{2M}(0)\cup B_{\frac{d(z)}{2}}(z))
$$ 
decomposes in a disjoint union of the resulting intersection with the set 
$$
\{|\zeta|>|\zeta-z|\}
$$ 
and its complement, namely $G_1$ and $G_2$. Then 
$$
\int_{G_1}\frac{dm(\zeta)}{|\zeta|^2\,|\zeta-z|^2}\le
\int_{\C\setminus B_{\frac{d(z)}{2}}(z))}\frac{dm(\zeta)}{|\zeta-z|^4}=\frac{16\pi}{d(z)^2},
$$ 
and 
$$
\int_{G_2}\frac{dm(\zeta)}{|\zeta|^2\,|\zeta-z|^2}\le
\int_{\C\setminus B_{2M}(0)}\frac{dm(\zeta)}{|\zeta|^4}=\frac{16\pi}{4M^2}\le\frac{16\pi\, 25}{|z|^2}.
$$
	
%\bigskip

The integral 
\begin{equation*}
|(V)|\le\int_{E_2(z)\cap B_{2M}(0)} |\varphi b_z|
\le\int_{B_{8M}(0)} \frac{dm(\zeta)}{R_0^2|\zeta-z|^2}\le\frac{64\pi M^2}{R_0^2}\, \frac{4}{d(z)^2}.\qedhere
\end{equation*}
%	\bigskip
%	\bigskip
%	
%	
%	
%	
\end{proof}
%\bigskip
%\bigskip
%
%\bigskip
%\bigskip

\subsection{Lipschitz estimates} 

We assume the same hypotheses, definitions and notation of the previous sections and subsections.

%\medskip

\begin{prop}\label{thm2} 
There exists a constant $K=K(\gamma,\, \Omega,\, R_0)$ such that for $\phi$ and~$\psi$ satisfying the conditions of Proposition \ref{thm1} we have,  for $f=\phi$ or $\psi$, that,   
if both $z,w\in\Omega$, or both $z,w\in\C\setminus\bar\Omega$, then 
$$
\frac{| B[f](z)- B[f](w)|}{|z-w|^\gamma}\le K_0\, \|f\|_\gamma.
$$ 
\end{prop}
%
%\bigskip
%\bigskip

\begin{rmk} As we sill see,
the theorem implies that $B[f]$ has a Lipschitz extension to $\bar\Omega$ and also to $\C\setminus\Omega$ but these extensions differ in a jump along $\partial\Omega$.
\end{rmk}

\subsubsection{Proof of Proposition \ref{thm2}:} 
We will estimate the quotients in the left hand side of the inequalities in the statement above, only in the case of $|z-w|\le\frac{R_0}{4}$. 

Also, most of the proof goes along for $\phi$ or $\psi$ with no distinction, so unless it be necessary, we will use $f$ for $\phi$ or $\psi$.

%\bigskip

The goal is to estimate $\frac{|B[f](z)-B[f](w)|}{|z-w|^\gamma}$ for $z,\, w\in\Omega$ and $z,\, w\in\C\setminus\bar\Omega$.

%\bigskip

Then we have that  
\begin{equation*}
\begin{split}
\frac{|B[f](z)-B[f](w)|}{|z-w|^\gamma}&\le \frac{|Q[f](z)-Q[f](w)|}{|z-w|^\gamma}+\frac{|L[f](z)-L[f](w)|}{|z-w|^\gamma}\\
&
=(I)+(II).
\end{split}
\end{equation*}

%\medskip

%\medskip

The term 
\begin{equation*}
\begin{split}
Q[f](z)-Q[f](w)&=\int_{\C\setminus B_{\delta(z)}(z)}f\, b_z-\int_{\C\setminus B_{\delta(w)}(w)}f\, b_w\\*[5pt]
&=\int_{\C\setminus (B_{\delta(z)}(z)\cup B_{\delta(w)}(w))}f\, (b_z-b_w)\\*[5pt]
&\quad+
\int_{B_{\delta(w)}(w)\setminus B_{\delta(z)}(z)}f\, b_z-\!\int_{B_{\delta(z)}(z)\setminus B_{\delta(w)}(w)}f\,b_w
\!=\!(1)+(2).\!
\end{split}
\end{equation*}

%\medskip

And we have, after a direct computation, that
\begin{equation}\label{id}
b_z(\zeta)-b_w(\zeta)=-(\bar z-\bar w)\, (\bar z+\bar w-2\, \bar\zeta)\, b_z(\zeta)\, b_w(\zeta).
\end{equation}

%\bigskip
%\bigskip

Then 
$$
(1)=-(\bar z-\bar w)\, \int_{\C\setminus (B_{\delta(z)}(z)\cup B_{\delta(w)}(w))}f(\zeta)\, (\bar z+\bar w-2\, \bar\zeta)\, b_z(\zeta)\, b_w(\zeta).
$$

%\medskip

Also 
\begin{equation*}
\begin{split}
L[f](z)-L[f](w)&=\int_{B_{\delta(z)}(z)}\frac{f(\zeta)-f(z)}{|\zeta-z|^\gamma}\, 
\frac{|\zeta-z|^\gamma}{(\bar\zeta-\bar z)^2}\, dm(\zeta)\\*[5pt]
&\quad-\int_{B_{\delta(w)}(w)}\frac{f(\zeta)-f(w)}{|\zeta-w|^\gamma}\, 
\frac{|\zeta-w|^\gamma}{(\bar\zeta-\bar w)^2}\, dm(\zeta)\\*[5pt]
&=\int_{B_{\delta(w)}(w)\cap B_{\delta(z)}(z)}\biggl[\frac{f(\zeta)-f(z)}{|\zeta-z|^\gamma}\, 
\frac{|\zeta-z|^\gamma}{(\bar\zeta-\bar z)^2}\\*[5pt]
&\hspace*{4cm}-\frac{f(\zeta)-f(w)}{|\zeta-w|^\gamma}\, 
\frac{|\zeta-w|^\gamma}{(\bar\zeta-\bar w)^2}\biggr]\, dm(\zeta)\\*[5pt]
&\quad+\int_{B_{\delta(z)}(z)\setminus B_{\delta(w)}(w)}\frac{f(\zeta)-f(z)}{|\zeta-z|^\gamma}\, 
\frac{|\zeta-z|^\gamma}{(\bar\zeta-\bar z)^2}\\*[5pt]
&\quad-
\int_{B_{\delta(w)}(w)\setminus B_{\delta(z)}(z)}\frac{f(\zeta)-f(w)}{|\zeta-w|^\gamma}\, 
\frac{|\zeta-w|^\gamma}{(\bar\zeta-\bar w)^2}\, dm(\zeta)\\*[5pt]
&=(3)+(4).
\end{split}
\end{equation*}

%\bigskip
%\bigskip

%\begin{enumerate}

\noindent
{\bf 1.} In general, the integral 
$$
(1)\!=\!-(\bar z-\bar w)\!\!\int_{\C\setminus(B_{\delta(z)}(z)\cup B_{\delta(w)}(w))\!} f(\zeta)
\frac{\bar w-\bar\tau+\bar z-\bar\tau-2(\bar\zeta-\bar\tau)}{(\bar\zeta\!-\!\bar\tau\!+\!\bar\tau-\bar z)^2
	(\bar\zeta\!-\!\bar\tau\!+\!\bar\tau\!-\!\bar w)^2}\,dm(\zeta),\!\!\!
$$ 
and since $\tau-z=\frac{z+w}{2}-z
=\frac{w-z}{2}\overset{\text{def}}{=}a$ and $\tau-w=-\frac{w-z}{2}=-a$, we have 
\begin{equation*}
\begin{split}
(1)&=2\,(\bar z-\bar w)\,\int_{\C\setminus (B_{\delta(z)}(z)\cup B_{\delta(w)}(w))} f(\zeta)\,
\frac{\bar\zeta-\bar\tau}{(\bar\zeta-\bar\tau+\bar a)^2\,(\bar\zeta-\bar\tau-\bar a)^2}\,dm(\zeta)\\*[5pt]
&\overset{\text{def}}{=}2\,(\bar z-\bar w)\,K[f](z,w).
\end{split}
\end{equation*}

Now, we have that, if 
$$
\lambda=\sqrt{\frac{\delta(z)^2-|a|^2+\delta(w)^2-|a|^2}{2}},
$$ 
then $B_{\lambda}(\tau)\subset B_{\delta(z)}(z)\cup B_{\delta(w)}(w)$, and, since $f$ is a bounded function, then 
\begin{equation*}
\begin{split}
|K[f](z,w)|&\le 
\int_{\C\setminus B_\lambda(\tau)}\biggl|f(\zeta)\,\frac{\bar\zeta-\bar\tau}{(\bar\zeta-\bar\tau+\bar a)^2\,
	(\bar\zeta-\bar\tau-\bar a)^2}\biggr|\,dm(\zeta)\\*[5pt]
&\le \|f\|_\infty\, \int_\lambda^\infty\int_0^{2\pi} \frac{r^2}{|r\,e^{-i\theta}+\bar a|^2\, |r\,e^{-i\theta}-\bar a|^2}\,dr\,d\theta\\*[5pt]
&=2\pi\|f\|_\infty\, \int_\lambda^\infty\frac{r^2}{|a|^4+r^4}\,dr=\frac{2\pi\|f\|_\infty}{|a|}\, \int_\frac{\lambda}{|a|}^\infty\frac{t^2}{1+t^4}\,dt\\*[5pt]
&\le \frac{2\pi\, \sqrt{2}\, \|f\|_\infty}{\sqrt{\delta(z)^2-|a|^2+\delta(w)^2-|a|^2}}\le\frac{4\pi\, \sqrt{2}\, \|f\|_\infty}{R_0},
\end{split}
\end{equation*}
by the lemma below, and the choice of $\upsilon_0$. 

%\bigskip
%\bigskip

\begin{lem} 
If $\upsilon_0\le\frac{R_0}{\sqrt{2}}$, then 
$$
\int_0^{2\pi} \frac{1}{|r\,e^{-i\theta}+\bar a|^2\,
|r\,e^{-i\theta}-\bar a|^2}\,d\theta=\frac{2\pi}{|a|^4+r^4}.
$$
\end{lem}

\begin{proof}
Since 
\begin{equation*}
\begin{split}
|r\,e^{-i\theta}+\bar a|^2\,	|r\,e^{-i\theta}-\bar a|^2&=|r^2\,e^{-2\, i\theta}-\bar a^2|^2\\*[5pt]
&=(r^2\,e^{-2\, i\theta}-\bar a^2)\, (r^2\,e^{2\, i\theta}-a^2)\\*[5pt]
&=r^4+|a|^4-(a^2\, r^2\,e^{-2\, i\theta}+\bar a^2\, r^2\,e^{2\, i\theta})\\*[5pt]
&=-\bar a^2\, r^2\,\varsigma^4+(r^4+|a|^4)\, \varsigma^2-a^2\, r^2,
\end{split}
\end{equation*}
for $\varsigma=e^{i\, \theta}$.
	
This is a polynomial in $\varsigma$, and the roots are $\pm\frac{r}{\bar a}$ and $\pm i\, \frac{a}{r}$, so if $r\neq|a|$, then 
\begin{equation*}
\begin{split}
\int_0^{2\pi} &\frac{1}{|(r^2\,e^{-2\, i\theta}-\bar a^2)\, (r^2\,e^{-2\, i\theta}-\bar a^2)|^2\,
	}\,d\theta\\*[5pt]
&=\frac{1}{i}\int_{\T} \frac{\varsigma}{-\bar a^2\, r^2\,\varsigma^4+(r^4+|a|^4)\, \varsigma^2-a^2\, r^2}\, d\varsigma\\*[5pt]
&=\frac{-1}{\bar a^2\, r^2\, i}\int_{\T} \frac{\varsigma}{(\varsigma-\frac{r}{\bar a})\, (\varsigma+\frac{r}{\bar a})\, (\varsigma-i\, \frac{a}{r})\, (\varsigma+i\, \frac{a}{r})}\, d\varsigma=(*).
\end{split}
\end{equation*}
	
	If $\upsilon_0\le\frac{R_0}{\sqrt{2}}$, then $\lambda\geq\sqrt{\frac{R_0^2}{4}-|a|^2}\geq |a|$. Then, by the residue theorem, we have 
\begin{equation*}
(*)=\frac{-2\pi}{\bar a^2\, r^2}\, \frac{-\bar a^2\, r^2}{|a|^4+r^4}=\frac{2\pi}{|a|^4+r^4}.\qedhere
\end{equation*}
\end{proof}
%\bigskip
%\bigskip
%
%\bigskip
%\bigskip

For the remaining terms, both in the decomposition of $Q$ and $L$,  the absolute and mutual positions of $z$ and $w$, specially related to $\partial\Omega$, will play an important role. 

For this purpose, we will consider the following (four) situations:

\begin{itemize}
\item[1)] If $B_{\delta(w)}(w)\cap B_{\delta(z)}(z)=\emptyset$, equivalent to the fact that $\delta(w)+\delta(z)\le|z-w|$ can never happen because then $\frac{R_0}{2}+\frac{R_0}{2}\le\frac{R_0}{\sqrt{2}}$. 
	
%\bigskip
	
\item[2)] The cases of $B_{\delta(w)}(w)\subset B_{\delta(z)}(z)$, or $B_{\delta(z)}(z)\subset B_{\delta(w)}(w)$, corresponding respectively to the facts that $\delta(w)+|z-w|\le\delta(z)$ or $\delta(z)+|z-w|\le\delta(w)$, so is 
$$
\frac{|z-w|}{\delta(z)+\delta(w)}\le 1-\frac{2\, \delta(w)}{\delta(z)+\delta(w)}
$$ 
or 
$$\frac{|z-w|}{\delta(z)+\delta(w)}\le 1-\frac{2\, \delta(z)}{\delta(z)+\delta(w)},
$$ 
respectively.
	
The situations are completely symmetric and in each case only one term in $(2)$ survives. 
	
%\bigskip
	
\item[3)] The case of the conditions 
$$
\begin{cases} 
B_{\delta(w)}(w)\cap B_{\delta(z)}(z)\neq\emptyset ,\\ 
B_{\delta(w)}(w)\cap B_{\delta(z)}(z)^c\neq \emptyset,
\end{cases}
$$ 
or reciprocally, are both satisfied.

In this case, we have 
$$
\begin{cases} \delta(w)+\delta(z)\geq|z-w|, \\
\delta(z)\le|z-w|+\delta(w),
\end{cases}
$$ 
or 
$$
1-\frac{2\, \delta(w)}{\delta(z)+\delta(w)}\le\frac{|z-w|}{\delta(z)+\delta(w)}\le 1,
$$ 
and, simultaneously, 
$$
1-\frac{2\, \delta(z)}{\delta(z)+\delta(w)}\le\frac{|z-w|}{\delta(z)+\delta(w)}\le 1.
$$ 
\end{itemize}

%\bigskip
%\bigskip
%
%
%\bigskip
%\bigskip

\noindent
{\bf 2.} For the term $(2)$, let us consider first the case 2). WLOG we assume that we are in the first situation of this case. Then 
$$
(2)=\int_{B_{\delta(w)}(w)\setminus B_{\delta(z)}(z)}f\, b_z.
$$

%\bigskip

\begin{itemize}

\item If $z\notin U_{\frac{R_0}{2}}$, we have $\delta(z)=d(z)$, and since $\delta(w)\geq|z-w|+\delta(z)$, then also $\delta(w)=d(w)$, and  
\begin{equation*}
\begin{split}
(2)&=\int_{B_{d(w)}(w)\setminus B_{d(z)}(z)}f\, b_z\\*[5pt]
&=\int_{B_{d(w)}(w)\setminus B_{d(z)}(z)}\frac{f(\zeta)-f(z)}{|\zeta-z|^\gamma}\, 
\frac{|\zeta-z|^\gamma}{(\bar\zeta-\bar z)^2}\, dm(\zeta),
\end{split}
\end{equation*}
by the cancellation Lemma \ref{can}, and since $B_{\delta(w)}(w)\subset\Omega$ or $B_{\delta(w)}(w)\subset\bar\Omega^c$, we have clearly that 
\begin{equation*}
\begin{split}
|(2)|&\le \|f\|_\gamma\, 2\pi\,  \int_{d(z)}^{d(w)}\frac{dr}{r^{1-\gamma}}=\|f\|_\gamma\, \frac{2\pi}{\gamma}\,  \{d(w)^\gamma-d(z)^\gamma\}\\*[5pt]
&\le\|f\|_\gamma\, 2\pi\, \frac{d(w)-d(z)}{d(z)^{1-\gamma}}\le\|f\|_\gamma\, \frac{2\pi\, 2^{1-\gamma}}{R_0^{1-\gamma}}\, |w-z|.
\end{split}
\end{equation*}

%\bigskip

\item If $z\in U_{\frac{R_0}{2}}$, we have that $\delta(z)=\frac{R_0}{2}$ and 
$$
\delta(w)\geq|z-w|+\delta(z)=|z-w|+\frac{R_0}{2}\geq\frac{ R_0}{2},
$$ 
then 
$$
(2)=\int_{B_{d(w)}(w)\setminus B_{\frac{R_0}{2}}(z)}\frac{f(\zeta)-f(z)}{|\zeta-z|^\gamma}\, 
\frac{|\zeta-z|^\gamma}{(\bar\zeta-\bar z)^2}\, dm(\zeta),
$$ 
so 
\begin{equation*}
\begin{split}
|(2)|&\le \|f\|_\gamma\, 2\pi\, \int_{\frac{R_0}{2}}^{d(w)+|z-w|}\frac{dr}{r^{1-\gamma}}\\*[5pt]
&=\|f\|_\gamma\, \frac{2\pi}{\gamma}\biggl \{(d(w)+|z-w|)^\gamma-\biggl(\frac{R_0}{2}\biggr)^\gamma\biggr\}\\*[5pt]
&\le\|f\|_\gamma\, \frac{2\pi\, 2^{1-\gamma}}{R_0^{1-\gamma}}\, \biggl(d(w)-\frac{R_0}{2}+|z-w|\biggr)\\*[5pt]
&\le\|f\|_\gamma\, \frac{2\pi\, 2^{1-\gamma}}{R_0^{1-\gamma}}\, (d(w)-d(z)+|z-w|),
\end{split}
\end{equation*}
having the same estimate as above. 
\end{itemize}

%\bigskip

In the situation of the case 3), 
$$
(2)=
\int_{B_{\delta(w)}(w)\setminus B_{\delta(z)}(z)}f\, b_z-\int_{B_{\delta(z)}(z)\setminus B_{\delta(w)}(w)}f\, b_w.
$$ 

\begin{itemize}
\item If $z,w\notin U_{\frac{R_0}{2}}$, since $|z-w|\le\frac{R_0}{4}$, we have that $z,w\in\Omega$ or $z,w\in\bar\Omega^c$ are the only possibilities, and also that 
$$
(2)=\int_{B_{d(w)}(w)\setminus B_{d(z)}(z)}f\, b_z-\int_{B_{d(z)}(z)\setminus B_{d(w)}(w)}f\, b_w.
$$

Since $B_{d(w)}(w)\subset B_{d(z)+|z-w|}(z)$, then the first term 
\begin{equation*}
\begin{split}
&\biggl|\int_{B_{d(w)}(w)\setminus B_{d(z)}(z)}f\, b_z\biggr|\\*[5pt]
&\quad\le
\|f\|_{L^\infty}\frac{1}{d(z)^2}\, m(C_{d(z)}^{d(z)+|z-w|}(z))\\*[5pt]
&\quad=\pi\, \|f\|_{L^\infty}\frac{d(z)+|z-w|+d(w)}{d(z)^2}\, (d(z)+|z-w|-d(w))\\*[5pt] 
&\quad\le\frac{4\pi\, \|f\|_{L^\infty}}{R_0}\, |z-w|\, \frac{d(z)+|z-w|+d(w)}{d(z)}\le\frac{16\pi\, \|f\|_{L^\infty}}{R_0}\, |z-w|,
\end{split}
\end{equation*}
 because $d(w)\le d(z)+|z-w|$.

%\medskip

The other term is similar.

%\bigskip

\item If $w\notin U_{\frac{R_0}{2}}$, but $z\in U_{\frac{R_0}{2}}$, we have that 
$$
(2)=\int_{B_{d(w)}(w)\setminus B_{\frac{R_0}{2}}(z)}f\, b_z-\int_{B_{\frac{R_0}{2}}(z)\setminus B_{d(w)}(w)}f\, b_w,
$$ 
and we consider several subcases, considering in each one $f=\phi$ and $f=\psi$ separately.  

\begin{itemize} 
\item If $z,w\in\Omega$, then $(2)$ is
\begin{multline*}
\int_{B_{d(w)}(w)\setminus B_{\frac{R_0}{2}}(z)}\phi\, b_z-\int_{B_{\frac{R_0}{2}}(z)\setminus B_{d(w)}(w)}\phi\, b_w\\
=(21)_\phi=(211)_\phi-(212)_\phi,
\end{multline*}
 or
 $$
 -\int_{B_{\frac{R_0}{2}}(z)}\psi\, b_w=(21)_\psi.
$$

%\medskip

Now, 
\begin{equation*}
\begin{split}
|(211)|_\phi &\le\|\phi\|_\infty\,\int_{B_{d(w)}(w)\setminus B_{\frac{R_0}{2}}(z)}\frac{1}{|\zeta-z|^2}\, dm\\*[5pt]
&\le\frac{4\, \|\phi\|_\infty}{R_0^2}\, m(B_{d(w)}(w)\setminus B_{\frac{R_0}{2}}(z))\\*[5pt]
&=\frac{4\pi\, \|\phi\|_\infty}{R_0^2}\, ((d(w)+|z-w|)^2-(\frac{R_0}{2})^2)\\*[5pt]
&\le\frac{16\, \pi\, \|\phi\|_\infty\, (\dmt(\Omega)+R_0)}{R_0^2}\, |w-z|.
\end{split}
\end{equation*}

%\bigskip

Also 
\begin{equation*}
\begin{split}
|(212)|_\phi&\le\frac{4\, \|\phi\|_\infty}{R_0^2}\, m(B_{\frac{R_0}{2}}(z)\setminus B_{d(w)}(w))\\*[5pt]
&\le\frac{4\, \|\phi\|_\infty}{R_0^2}\, m(C_{d(w)}^{\frac{R_0}{2}+|z-w|}(w))\\*[5pt]
&\le\frac{4\pi\, \|\phi\|_\infty}{R_0^2}\biggl(\biggl|\frac{R_0}{2}\!+\!|z\!-\!w|\!-\!d(w)\biggr|\biggr)\biggl (\frac{R_0}{2}\!+\!|z\!-\!w|\!+\!d(w)\biggr)\!\!\\*[5pt]
&\le\frac{16\, \pi\, \|\phi\|_\infty\, (\dmt(\Omega)+R_0)}{R_0^2}\, |w-z|.
\end{split}
\end{equation*}

%\bigskip

The term 
\begin{equation*}
\begin{split}
|(21)_\psi|&\le\|\psi\|_{\infty}\, \int_{B_{\frac{R_0}{2}}(z)\cap\Omega^c}|b_w|\le\|\psi\|_{\infty}\, \frac{4}{R_0^2}\,  m(C_{d(z)}^{\frac{R_0}{2}}(z))\\*[5pt]
&\le\|\psi\|_{\infty}\, \frac{2\pi}{R_0^2}\,  \pi\biggl(\biggl(\frac{R_0}{2}\biggr)^2-d(z)^2\biggr)\\*[5pt]
&\le\|\psi\|_{\infty}\, \frac{2\pi\, (\dmt(\Omega)+R_0)}{R_0^2}\, (d(w)-d(z)),
\end{split}
\end{equation*}
because $d(w)\geq\frac{R_0}{2}$.

%\bigskip

\item If $z,w\in\Omega^c$, then $(2)$ is
$$
-\int_{B_{\frac{R_0}{2}}(z)}\phi\, b_w=(22)_\phi,
$$ 
or 
$$
\int_{B_{d(w)}(w)\setminus B_{\frac{R_0}{2}}(z)}\psi\, b_z-\int_{B_{\frac{R_0}{2}}(z)\setminus B_{d(w)}(w)}\psi\, b_w=(22)_\psi.
$$	

%\medskip

The term $(22)_\phi$ is analogous to $(21)_\psi,$ so 
$$
|(22)_\phi|\le\|\phi\|_{\infty}\, \frac{2\pi\, (\dmt(\Omega)+R_0)}{R_0^2}\, (d(w)-d(z)).
$$

%\medskip

Also the term 
$$
(22)_\psi=(221)_\psi-(222)_\psi
$$ 
is analogous to $(21)_\phi$ and 
\begin{equation*}
\begin{split}
|(221)_\psi|&\le\frac{4\pi\, \|\phi\|_\infty}{R_0^2}\biggl (d(w)^2-\biggl(\frac{R_0}{2}\biggr)^2\biggr)\\*[5pt]
&\le\frac{4\pi\, \|\phi\|_\infty}{R_0^2}\biggl (d(w)-\frac{R_0}{2}\biggr)\biggl(d(w)+\frac{R_0}{2}\biggr)\\*[5pt]
&\le\frac{8\pi\, (\dmt(\Omega)+R_0)\, \|\phi\|_\infty}{R_0^2}\, (d(w)-d(z)).
\end{split}
\end{equation*}

%\medskip

The term $(222)_\psi$ is similar to the previous one. 
\end{itemize}

%\bigskip

\item If $z,w\in U_{\frac{R_0}{2}}$, we have that 
$$
(2)=
\int_{B_{\frac{R_0}{2}}(w)\setminus B_{\frac{R_0}{2}}(z)}f\, b_z-\int_{B_{\frac{R_0}{2}}(z)\setminus B_{\frac{R_0}{2}}(w)}f\, b_w,
$$ 
and then 
$$
|(2)|\le2\, \frac{4\, \|f\|_{L^\infty}}{R_0^2}\, m(B_{\frac{R_0}{2}}(w)\setminus B_{\frac{R_0}{2}}(z))
$$ 
and the previous procedure applies. 
\end{itemize}

%\bigskip
%\bigskip

\noindent
{\bf 3.} Now, the term 
$$
(3)\!=\!\int_{B_{\delta(w)}(w)\cap B_{\delta(z)}(z)\!}\biggl[\frac{f(\zeta)\!-\!f(z)}{|\zeta-z|^\gamma}\, 
\frac{|\zeta-z|^\gamma}{(\bar\zeta-\bar z)^2}-\frac{f(\zeta)\!-\!f(w)}{|\zeta-w|^\gamma}\, 
\frac{|\zeta-w|^\gamma}{(\bar\zeta-\bar w)^2}\biggr]\, dm(\zeta).\!\!\!
$$

%\medskip

%Since $|z-w|<\frac{R_0}{4}$, then $z,w\in %B_{\frac{\delta(w)}{2}}(w)\cap %B_{\frac{\delta(z)}{2}}(z)\subset B_{\delta(w)}(w)\cap %B_{\delta(z)}(z)$, and we can reduce the study to this %situation.

In such case, $d(z,\partial B_{\delta(w)}(w))=\delta(w)-|z-w|$ and $d(w,\partial B_{\delta(z)}(z))=\delta(z)-|z-w|$, so the maximum radius of a ball centered at $z$ and contained in $B_{\delta(z)}(z)\cap B_{\delta(w)}(w)$ is equal to $\min\{\delta(z),\,\delta(w)-|z-w|\}$. For $w$ we have an analogous expression. And since $\delta(z)\geq\frac{R_0}{2}\geq2\, |z-w|$, and analogously for $w$, we can consider, defining 
\begin{align*}
k^\gamma_z(\zeta)&=|\zeta-z|^\gamma,\\
\Delta^\gamma_z (\zeta)&=\frac{f(\zeta)-f(z)}{k_z^\gamma(\zeta)},
\end{align*} 
the decomposition 
\begin{equation*}
\begin{split}
(3)&=\int_{B_{\delta(w)}(w)\cap B_{\delta(z)}(z)}[\Delta^\gamma_z\, k^\gamma_z\, b_z- \Delta^\gamma_w\, k^\gamma_w\, b_w]\\*[5pt]
&=\int_{B_{\frac{|z-w|}{2}}(z)}\Delta^\gamma_z\, k^\gamma_z\, b_z-\int_{B_{\frac{|z-w|}{2}}(w)}\Delta^\gamma_w\, k^\gamma_w\, b_w\\*[5pt]
&\quad+\int_{(B_{\delta(w)}(w)\cap B_{\delta(z)}(z))\setminus B_{\frac{|z-w|}{2}}(z)}\Delta^\gamma_z\, k^\gamma_z\, b_z\\*[5pt]
&\quad-\int_{(B_{\delta(w)}(w)\cap B_{\delta(z)}(z))\setminus B_{\frac{|z-w|}{2}}(w)}\Delta^\gamma_w\, k^\gamma_w\, b_w\\*[5pt]
&=(31)+(32).
\end{split}
\end{equation*}

Now, 
$$
|(31)|\le 2\, \|f\|_\gamma\, 2\pi\, \int_0^{\frac{|z-w|}{2}}\frac{dr}{r^{1-\gamma}}=\frac{4\pi\, \|f\|_\gamma}{\gamma\, 2^\gamma}\, |z-w|^\gamma.
$$

%\medskip

Since for $\tau=\frac{z+w}{2}$ we have $B_{\frac{|z-w|}{2}}(z)\cup B_{\frac{|z-w|}{2}}(w)\subset B_{|z-w|}(\tau)$, we have that 
\begin{equation*}
\begin{split}
(32)&=\int_{B_{|z-w|}(\tau)\setminus B_{\frac{|z-w|}{2}}(z)}\Delta^\gamma_z\, k^\gamma_z\, b_z-\int_{B_{|z-w|}(\tau)\setminus B_{\frac{|z-w|}{2}}(w)}\Delta^\gamma_w\, k^\gamma_w\, b_w\\*[5pt]
&\quad+\int_{(B_{\delta(w)}(w)\cap B_{\delta(z)}(z))\setminus B_{|z-w|}(\tau)}[\Delta^\gamma_z\, k^\gamma_z\, b_z-\Delta^\gamma_w\, k^\gamma_w\, b_w]\\*[5pt]
&=(321)+(322),
\end{split}
\end{equation*}
and the term 
$$
|(321)|\le 2\,  \|f\|_\gamma\, \frac{4\pi}{\gamma}\biggl(\biggl(\frac{3}{2}\biggr)^\gamma-\biggl(\frac{1}{2}\biggr)^\gamma\biggr)\, |z-w|^\gamma.
$$

For the term $(322)$, we have that 
$$
\Delta^\gamma_z\, k^\gamma_z\, b_z-\Delta^\gamma_w\, k^\gamma_w\, b_w=(f(\zeta)-f(z))\, b_z(\zeta)
-(f(\zeta)-f(w))\, b_w(\zeta)=(*).
$$

%\medskip

\begin{itemize}

\item If $z\notin U_{\frac{R_0}{4}}$, then, as $|z-w|<\frac{R_0}{4}$, we have that $w,\tau=\frac{z+w}{2}\in B_{\frac{R_0}{4}}(z)$, 
%so all are in $\bar\Omega$ or all are in $\Omega^c$,
 and then, using the decomposition
$$
(*)\!=\!(f(\zeta)-f(\tau)) (b_z-b_w)(\zeta)+
(f(w)-f(\tau))\, b_w(\zeta)-(f(z)-f(\tau))\, b_z(\zeta),\!\!
$$
we have 
\begin{equation*}
\begin{split}
(322)&=\int_{(B_{\delta(w)}(w)\cap B_{\delta(z)}(z))\setminus B_{|z-w|}(\tau)}\Delta^\gamma_\tau\, k^\gamma_\tau\, (b_z-b_w)\\*[5pt]
&\quad+(f(w)-f(\tau))\, \int_{(B_{\delta(w)}(w)\cap B_{\delta(z)}(z))\setminus B_{|z-w|}(\tau)}b_w\\*[5pt]
&\quad-(f(z)-f(\tau))\, \int_{(B_{\delta(w)}(w)\cap B_{\delta(z)}(z))\setminus B_{|z-w|}(\tau)}b_z\\*[5pt]
&=(3221)
+(3222).
\end{split}
\end{equation*}

%\medskip

For the integral $(3321)$ we use the change of variables 
$$
\zeta\rightarrow s=\frac{\zeta-\tau}{a}=\phi(\zeta),
$$ 
where $a=\frac{z-w}{2}$, then $\phi(\tau)=0$ and $J\phi=|a|^{-2}$. 

We have 
\begin{equation*}
\begin{split}
k^\gamma_\tau(\zeta)\, (b_z-b_w)(\zeta)&=|\zeta-\tau|^\gamma\biggl \{\frac{1}{(\bar\zeta-\bar z)^2}-\frac{1}{(\bar\zeta-\bar w)^2}\biggr\}\\*[5pt]
&=|a\, s|^\gamma\biggl \{\frac{1}{(\bar a\, (\bar s+1))^2}-\frac{1}{(\bar a\, (\bar s-1))^2}\biggr\}\\*[5pt]
&=\frac{|a|^\gamma}{\bar a^2}\,  \frac{-|s|^\gamma\, 4\bar s}{(\bar s+1)^2\, (\bar s-1)^2}.
\end{split}
\end{equation*} 

And then, since $\phi((B_{\delta(w)}(w))=B_{\frac{\delta(w)}{|a|}}(-1)$, $\phi((B_{\delta(z)}(z))=B_{\frac{\delta(z)}{|a|}}(1)$ and $\phi((B_{|z-w|}(\tau))=B_{2}(0)$, we have that 
\begin{multline*}
(3221)=-\frac{|a|^{\gamma+2}}{\bar a^2}\\
\times \int_{(B_{\frac{\delta(w)}{|a|}}(-1)\cap B_{\frac{\delta(z)}{|a|}}(1))\setminus B_{2}(0)}\Delta^\gamma_\tau(\tau+a\, s)\, \frac{|s|^\gamma\, 4\bar s}{(\bar s+1)^2\, (\bar s-1)^2}\, dm(s),
\end{multline*}
so 
$$
|(3221)|\le|a|^\gamma\, \|f\|_\gamma\, \int_{\C\setminus B_{2}(0)}\frac{4|s|^{1+\gamma}}{(\bar s+1)^2\, (\bar s-1)^2}\, dm(s).
$$

%\medskip 

For the term $(3322)$, we have 

\begin{lem} 
If  $z\notin U_{\frac{R_0}{4}}$ and $|z-w|<\frac{R_0}{4}$, then 
$$
\biggl|\int_{(B_{\delta(w)}(w)\cap B_{\delta(z)}(z))\setminus B_{|z-w|}(\tau)}b_z\biggr|\le 16\pi .
$$
\end{lem}

\begin{proof} 
Using Stokes formula, 
\begin{equation*}
\begin{split}
&\int_{(B_{\delta(w)}(w)\cap B_{\delta(z)}(z))\setminus B_{|z-w|}(\tau)}b_w\\*[5pt]
&=\frac{1}{2i}\biggl \{\int_{(\partial B_{\delta(w)}(w))\cap B_{\delta(z)}(z)}\!+\!\int_{B_{\delta(w)}(w)\cap\partial(B_{\delta(z)}(z))}\!-\!\int_{\partial B_{|z-w|}(\tau)}\biggr\}\, \frac{d\zeta}{\bar\zeta-\bar z}\\*[5pt]
&=(I)+(II)-(III).
\end{split}
\end{equation*}

Then using the change $\zeta=\tau+2\, |a|\, e^{i\theta}$, we have 
\begin{equation*}
\begin{split}
(III)&=\int_{\partial B_{2\, |a|}(\tau)}\, \frac{d\zeta}{\bar\zeta-\bar\tau-\bar a}=\int_{0}^{2\pi}\, \frac{2i\, |a|\, e^{i\theta}\, d\theta}{2\, |a|\, e^{-i\theta}-\bar a}\\*[5pt]
&=2i\, \frac{|a|}{\bar a}\, \int_{0}^{2\pi}\, \frac{ e^{2i\theta}\, d\theta}{2\, \frac{|a|}{\bar a}-e^{i\theta}}=2\, \frac{|a|}{\bar a}\, \int_{S_1(0)}\, \frac{ s\, ds}{2\, \frac{|a|}{\bar a}-s}=0,
\end{split}
\end{equation*}
by the Cauchy formula.

%\medskip

The integral 
$$
|(II)|\le \frac{1}{\delta(z)}\, \int_{\partial B_{\delta(z)}(z)} |d\zeta|=2\pi.
$$

%\medskip

In the same way, in our situation,  
$$
|(I)|\le \frac{1}{\delta(z)-|z-w|}\, \int_{\partial B_{\delta(w)}(w)}|d\zeta|=2\pi\, \frac{\delta(w)}{\delta(w)-|z-w|}.
$$

If $d(w)\le\frac{R_0}{2}$, then 
$$
|(I)|\le 2\pi\, \frac{\frac{R_0}{2}}{\frac{R_0}{2}-\frac{R_0}{4}}=4\pi.
$$

If $d(w)>\frac{R_0}{2}$, then 
\begin{equation*}
|(I)|\le 2\pi\, \frac{d(w)}{d(w)-\frac{R_0}{4}}\le 2\pi\, \frac{d(w)}{\frac{d(w)}{2}}=4\pi.\qedhere
\end{equation*}
\end{proof}
%\bigskip
%\bigskip

The lemma implies that 
$$
|(3222)|\le \frac{32\pi}{2^\gamma}\, |z-w|^\gamma\, \|f\|_\gamma.
$$ 

%\bigskip

\item The situation is completely equivalent to the case of $w\notin U_{\frac{R_0}{4}}$.

%\bigskip

\item If $z,w\in U_{\frac{R_0}{4}}$, then $\delta(z)=\delta(w)=\frac{R_0}{2}$, and if $z,\, w\in W$ where $W=\Omega$ or $\C\setminus\bar\Omega$, taking $\tau'$  a point in $\bar W$ at minimum distance of $\tau$
\begin{equation*}
\begin{split}
(322)&=\int_{(B_{\frac{R_0}{2}}(w)\cap B_{\frac{R_0}{2}}(z))\setminus B_{|z-w|}(\tau)} [(f(\zeta)-f(z)) b_z(\zeta)\\*[5pt]
&\hspace*{5.5cm}
-(f(\zeta)-f(w)) b_w(\zeta)]\\*[5pt]
&=\int_{(B_{\frac{R_0}{2}}(w)\cap B_{\frac{R_0}{2}}(z))\setminus B_{|z-w|}(\tau)}(f(\zeta)-f(\tau'))\, [b_z(\zeta)-b_w(\zeta)]\\*[5pt]
&\quad+(f(\tau')-f(z))\, \int_{(B_{\frac{R_0}{2}}(w)\cap B_{\frac{R_0}{2}}(z))\setminus B_{|z-w|}(\tau)}b_z(\zeta)\\*[5pt]
&\quad-(f(\tau')-f(w))\, \int_{(B_{\frac{R_0}{2}}(w)\cap B_{\frac{R_0}{2}}(z))\setminus B_{|z-w|}(\tau)}b_z(\zeta)\\*[5pt]
&=(3223)+(3224).
\end{split} 
\end{equation*}

%\bigskip

Now, using the identity \eqref{id}, we have  
\begin{multline*}
(3223)=-(\bar z-\bar w)
\int_{(B_{\frac{R_0}{2}}(w)\cap B_{\frac{R_0}{2}}(z))\setminus B_{|z-w|}(\tau)}\biggl (\frac{f(\zeta)-f(\tau')}{|\zeta-\tau'|^\gamma}\biggr)\\*[5pt]
\times k^\gamma_{\tau'}(\zeta)\, (\bar z+\bar w-2\, \bar\zeta)\, b_z(\zeta)\, b_w(\zeta).
\end{multline*}

Since $\tau'\in \partial W$ then, 
$$
|\tau'-\tau|\le|\tau-z|=\frac{|z-w|}{2}=|a|,
$$ 
then 
$$
|z-\tau'|\le|z-\tau|+|\tau-\tau'|\le 2\, |a|
$$ 
and 
$$
|\zeta-\tau'|\le|\zeta-\tau|+|\tau-\tau'|\le|\zeta-\tau|+|a|,
$$ 
so 
$$
|\zeta-\tau'|^\gamma\le(|\zeta-\tau|+|a|)^\gamma\le\frac{\gamma\, |\zeta-\tau|}{|a|^{1-\gamma}},
$$ 
and then 
\begin{multline*}
|(3223)|\le \frac{\gamma\, |z-w|}{|a|^{1-\gamma}}\, \|f\|_\gamma \\
\times\int_{(B_{\frac{R_0}{2}}(w)\cap B_{\frac{R_0}{2}}(z))\setminus B_{|z-w|}(\tau)}\frac{|\zeta-\tau|^\gamma\, (|z-\zeta|+|\zeta-w|)}{|\zeta-z|^2\, |\zeta-w|^2}\, dm(\zeta)
\end{multline*}
and the estimate is identical to the case of $(322)$. 

%\medskip

Also, a repetition of the arguments shows that the term $(3224)$ has the same estimate as the term $(3222)$.
\end{itemize}

%\bigskip
%\bigskip

\noindent
{\bf 4.} For the term 
\begin{equation*}
\begin{split}
(4)&=\int_{B_{\delta(z)}(z)\setminus B_{\delta(w)}(w)}\frac{f(\zeta)-f(z)}{|\zeta-z|^\gamma}\, 
\frac{|\zeta-z|^\gamma}{(\bar\zeta-\bar z)^2}\\*[5pt]
&\quad-
\int_{B_{\delta(w)}(w)\setminus B_{\delta(z)}(z)}\frac{f(\zeta)-f(w)}{|\zeta-w|^\gamma}\, 
\frac{|\zeta-w|^\gamma}{(\bar\zeta-\bar w)^2}\, dm(\zeta)=(41)+(42).
\end{split}
\end{equation*}

%\medskip

The term 
$$
|(41)|\le \|f\|_\gamma\, \int_{B_{\delta(z)}(z)\setminus B_{\delta(w)}(w)}\frac{dm(\zeta)}{|\zeta-z|^{1-\gamma}},
$$ 
and since $|z-w|\le\frac{R_0}{4}<\delta(w),\, \delta(z)$, then we have that
$z\in B_{\delta(w)}(w)$. Moreover, if $|\zeta-z|<\delta(w)-|z-w|$, then $\delta(w)>|\zeta-z|+|z-w|\geq|\zeta-w|$. This implies that 
$$
B_{\delta(w)-|w-z|}(z)\subset B_{\delta(w)}(w)
$$ 
and then, since the $\delta(w)-|z-w|<\delta(z)$, otherwise the domain of integration is empty, the previous integral is bounded by 
\begin{equation*}
\begin{split}
\int_{C^{\delta(z)}_{\delta(w)-|z-w|}(z)}&\frac{dm(\zeta)}{|\zeta-z|^{1-\gamma}}\\*[5pt]
&=2\pi\, \int^{\delta(z)}_{\delta(w)-|z-w|}\frac{dr}{r^{1-\gamma}}\\*[5pt]
&=\frac{2\pi}{\gamma}\, \{\delta(z)^\gamma-(\delta(w)-|z-w|)^\gamma\}\\*[5pt]
&=2\pi\, \frac{\delta(z)-(\delta(w)-|z-w|)}{(\delta(w)-|z-w|+\lambda\, (\delta(z)-(\delta(w)-|z-w|))^{1-\gamma}}\\*[5pt]
&\le2\pi\, \frac{\delta(z)-\delta(w)+|z-w|)}{(\delta(w)-|z-w|)^{1-\gamma}}
\le2\pi\, \frac{\delta(z)-\delta(w)+|z-w|)}{(\frac{R_0}{4})^{1-\gamma}},
\end{split}
\end{equation*}
where

 $\lambda\in(0,1)$.

%\bigskip

The term $(42)$ is symmetric and analogous to $(41)$.

%\bigskip

With the previous arguments we have the conjugate Beurling transforms of $\phi$ and $\psi$ are H\"older at $\Omega$ or $(\bar\Omega)^c$. The following lemma completes the case of~$\partial\Omega$.

\begin{lem}
If $W=\Omega$ or $W=(\bar\Omega)^c$ and $f\in \Lip{(\gamma, W)}$ and $\|f\|_\gamma<+\infty$, then $f$ extends to a Lipschitz function on $\bar W$, with the same Lipschitz norm.  
\end{lem}

\begin{proof} 
Since $\partial W$ is compact, then $f$ is uniformly continuous in $\bar U_{\frac{R_0}{4}}$.

If $z,w\in\partial\Omega$ and $|z-w|<\frac{R_0}{4}$, for $U_{\frac{R_0}{4}}\cap W$ we have 
$$
 f(z)-f(w)=f(z)-f(z')+f(z')-f(w')+f(w')-f(w)=(1)+(2)+(3).
 $$  

If $|z'-z|,\, |w'-w|<\min\{\delta,\, \frac{|z-w|}{3}\}$, then 
$$
(1),\, (3)\le M\, |z-w|^\gamma.
$$ 

Also 
$$
|z'-w'|\le \frac{5}{3}\, |z-w|,
$$ 
so 
\begin{equation*}
(2)\le M\, |z-w|^\gamma.\qedhere
\end{equation*}
\end{proof}
%\bigskip
%\bigskip

Finally, we can apply this theorem to the last term and we have

\begin{cor} 
There exists a constant $C_1$, depending only on $R_0$, $\upsilon_0$, $\dmt(\Omega)$, $\gamma$, such that 
$$
|f(z)\, \Theta_\Omega^{\frac{R_0}{2}}(z)-f(w)\, \Theta_\Omega^{\frac{R_0}{2}}(w)|\le |z-w|^\gamma\, \|f\|_\gamma\, C_1.
$$
\end{cor}

%\bigskip
%\bigskip

%\end{enumerate}

\bigskip
\bigskip

%\medskip

\subparagraph{Aknowledgements:} The authors are grateful to Joan Verdera for some useful conversations and comments on the results of the paper. Both authors are  partially supported by grants 2017-SGR-0395 (Generalitat de Catalunya) and MDM-2014-044 (MICINN, Spain). The second named author is also partially supported by MTM-2016-75390 (MINECO, Spain).

\medskip

\newpage

\newpage

\bibliography{trajan}
\bibliographystyle{alpha}

\noindent
{\small
\begin{tabular}{@{}l}
J.\ M.\ Burgu\'es\\
J.\ Mateu\\
Departament de Matem\`atiques\\
Universitat Aut\`onoma de Barcelona\\
08193 Bellaterra, Barcelona, Catalonia\\
{\it E-mail:} {\tt josep@mat.uab.cat}\\
{\it E-mail:} {\tt mateu@mat.uab.cat}
\end{tabular}}

\end{document}